\documentclass[11pt,a4paper,leqno]{amsart}

\usepackage{silence}
\usepackage[T1]{fontenc}
\usepackage[utf8]{inputenc}
\usepackage[english]{babel}

\usepackage[dvipsnames]{xcolor}
\definecolor{cornellred}{rgb}{0.7, 0.11, 0.11}
\definecolor{britishracinggreen}{rgb}{0.0, 0.26, 0.15}
\definecolor{cobalt}{rgb}{0.0, 0.28, 0.67}

\usepackage[
  a4paper,
  top=2.8cm,
  bottom=2.8cm,
  left=2.8cm,
  right=2.8cm,
  marginparwidth=60pt
]{geometry}

\usepackage{amsmath}
\usepackage{amsthm}
\usepackage[utopia]{mathdesign}
\usepackage{mathtools}

\usepackage{graphicx}
\usepackage{thmtools}
\usepackage{braket}
\usepackage{caption}
\usepackage{comment}
\usepackage{stmaryrd}
\usepackage{enumitem}
\usepackage{multirow}
\usepackage{booktabs}
\usepackage{relsize}
\usepackage[foot]{amsaddr}
\usepackage{float}
\usepackage[colorinlistoftodos]{todonotes}
\usepackage{tikz}
\usepackage{tikz-cd}
\usepackage{adjustbox}
\usepackage{rotating}

\usetikzlibrary{decorations.markings}
\usetikzlibrary{patterns,matrix,shapes,intersections,arrows,decorations.pathmorphing}

\usepackage[nopatch=showhyphens]{microtype}

\DeclareFontShape{T1}{mdput}{m}{up}{<->ssub*mdput/m/n}{}
\DeclareFontShape{T1}{mdput}{b}{up}{<->ssub*mdput/b/n}{}
\DeclareFontShape{T1}{mdput}{bx}{up}{<->ssub*mdput/b/n}{}

\SetSymbolFont{stmry}{bold}{U}{stmry}{m}{n}

\numberwithin{equation}{section}
\allowdisplaybreaks
\usepackage[
  colorlinks,
  bookmarks,
  hypertexnames=false
]{hyperref}

\hypersetup{
  colorlinks,
  citecolor=britishracinggreen,
  filecolor=black,
  linkcolor=cobalt,
  urlcolor=purple
}

\usepackage[capitalise]{cleveref}

\DeclareSymbolFont{cmarrows}{OMS}{cmsy}{m}{n}
\SetSymbolFont{cmarrows}{bold}{OMS}{cmsy}{b}{n}
\DeclareMathSymbol{\cmminus}{\mathbin}{cmarrows}{"00}

\DeclareMathSymbol{\leftrightarrow}{\mathrel}{cmarrows}{"24}
\DeclareMathSymbol{\leftarrow}{\mathrel}{cmarrows}{"20}
   
\DeclareMathSymbol{\rightarrow}{\mathrel}{cmarrows}{"21}
   \let\to=\rightarrow
\DeclareMathSymbol{\mapstochar}{\mathrel}{cmarrows}{"37}
   \def\mapsto{\mapstochar\rightarrow}    

\DeclareSymbolFont{usualmathcal}{OMS}{cmsy}{m}{n}
\DeclareSymbolFontAlphabet{\mathcal}{usualmathcal}
\DeclareMathAlphabet\BCal{OMS}{cmsy}{b}{n} 

\newcommand{\ext}{\operatorname{ext}}
\newcommand{\Ext}{\operatorname{Ext}}
\newcommand{\Hom}{\operatorname{Hom}}

\newcommand{\dimhom}{\operatorname{hom}}

\DeclareMathOperator{\lHom}{\mathscr{H}\kern-0.3em\mathit{om}}
\DeclareMathOperator{\lEnd}{\mathscr{E}\kern-0.1em\mathit{nd}\,}
\DeclareMathOperator{\RRlHom}{\mathbf{R}\kern-0.025em\mathscr{H}\kern-0.3em\mathit{om}}
\DeclareMathOperator{\RHom}{{\mathbf{R}\mathrm{Hom}}}
\newcommand{\lExt}{\mathop{\mathscr{E}\kern-0.2em\textit{xt}}\nolimits}
\DeclareMathOperator{\LTor}{{\mathscr{T}\kern-0.2em\mathit{or}}}

\newcommand{\Id}{\operatorname{Id}}
\newcommand{\id}{\operatorname{id}}

\newcommand{\Spec}{\operatorname{Spec}}
\newcommand{\Supp}{\operatorname{Supp}}

\newcommand{\surj}{\mathsf{surj}}

\newcommand{\dd}{\mathrm{d}}

\newcommand{\simto}{\,\widetilde{\to}\,}
\newcommand{\ch}{{\mathrm{ch}}}
\newcommand{\hd}{\operatorname{hd}}
\newcommand{\Pic}{\mathop{\rm Pic}\nolimits}
\newcommand{\PT}{\mathsf{PT}}

\newcommand{\RR}{\mathbf{R}}
\newcommand{\HH}{\mathrm{H}}
\newcommand{\OO}{\mathscr{O}}

\newcommand{\into}{\hookrightarrow}
\newcommand{\onto}{\mathrel{\mathrlap{\rightarrow}\mkern1mu\rightarrow}}
\newcommand{\derived}{\mathbf{D}}
\newcommand{\bounded}{\mathrm{b}}

\newcommand{\op}[1]{{\mathscr O}_{\mathbb{P}^{#1}}}

\makeatletter
  \newcommand{\mylabel}[2]{%
  #2%
  \phantomsection%
  \def\@currentlabel{#2}%
  \label{#1}%
  }
\makeatother

\newcommand*{\defeq}{\mathrel{\vcenter{\baselineskip0.5ex \lineskiplimit0pt
                     \hbox{\scriptsize.}\hbox{\scriptsize.}}}=}

\renewcommand{\leq}{\leqslant}
\renewcommand{\geq}{\geqslant}
\renewcommand{\le}{\leqslant}
\renewcommand{\ge}{\geqslant}

\newcommand{\boldit}[1]{\boldsymbol{#1}}

\DeclareMathOperator{\Kur}{\mathsf{Kur}}
\DeclareMathOperator{\Hilb}{\mathsf{Hilb}}
\DeclareMathOperator{\Quot}{\mathsf{Quot}}
\DeclareMathOperator{\Sym}{\mathsf{Sym}}

\DeclareMathOperator{\Gr}{\mathsf{Gr}}
\DeclareMathOperator{\red}{\mathsf{red}}
\DeclareMathOperator{\vir}{\mathsf{vir}}

\DeclareMathOperator{\tr}{\mathsf{tr}}
\newcommand{\num}{\mathsf{num}}

\newcommand{\main}{{\mathsf{main}}}
\newcommand{\Stab}{{\mathsf{Stab}}}

\newcommand{\st}{\mathsf{st}}
\newcommand{\pug}{\mathsf{pug}}

\DeclareMathOperator{\pr}{pr}

\DeclareMathOperator{\im}{im}
\DeclareMathOperator{\Sets}{Sets}
\DeclareMathOperator{\Sch}{Sch}

\DeclareMathOperator{\coker}{coker}

\DeclareMathOperator{\Perf}{Perf}
\DeclareMathOperator{\Coh}{Coh}

\DeclareMathOperator{\vd}{vd}

\DeclareMathOperator{\Art}{Art}
\DeclareMathOperator{\Sing}{Sing}

\DeclareMathOperator{\opp}{op}
\DeclareMathOperator{\CCoh}{\mathbf{Coh}}

\DeclareMathOperator{\Var}{Var}
\DeclareMathOperator{\Rep}{\mathsf{Rep}}
\DeclareMathOperator{\GL}{GL}
\DeclareMathOperator{\PGL}{PGL}
\DeclareMathOperator{\Mat}{Mat}

\DeclareMathOperator{\rk}{rk}

\DeclareMathOperator{\length}{length}

\DeclareMathOperator{\Aut}{Aut}

\DeclareFontFamily{OT1}{rsfs}{}
\DeclareFontShape{OT1}{rsfs}{n}{it}{<-> rsfs10}{}
\DeclareMathAlphabet{\curly}{OT1}{rsfs}{n}{it}

\makeatletter
  \newenvironment{proofof}[1]{\par
  \pushQED{\qed}%
  \normalfont \topsep6\p@\@plus6\p@\relax
  \trivlist
  \item[\hskip3\labelsep
        \itshape
    Proof of #1\@addpunct{.}]\ignorespaces
  }{%
  \popQED\endtrivlist\@endpefalse}
\makeatother

\tikzset{commutative diagrams/.cd,
mysymbol/.style={start anchor=center,end anchor=center,draw=none}}
\newcommand\MySymb[2][\square]{%
  \arrow[mysymbol]{#2}[description]{#1}}
\tikzset{
shift up/.style={
to path={([yshift=#1]\tikztostart.east) -- ([yshift=#1]\tikztotarget.west) \tikztonodes}
}
}

\theoremstyle{definition}

\newtheorem{definition}{Definition}[section]
\newtheorem{assumption}{Assumption}[section]

\newtheorem{conventions}[definition]{Conventions}

\newtheorem{situation}[definition]{Situation}

\newtheorem{notation}[definition]{Notation}
\newtheorem{remark}[definition]{Remark}

\newtheoremstyle{thm} 
        {3mm}
        {3mm}
        {\slshape}
        {0mm}
        {\bfseries}
        {.}
        {1mm}
        {}
\theoremstyle{thm}
\newtheorem{theorem}[definition]{Theorem}
\newtheorem{corollary}[definition]{Corollary}
\newtheorem{lemma}[definition]{Lemma}
\newtheorem{prop}[definition]{Proposition}
\newtheorem{thm}{Theorem}

\makeatletter

\@ifundefined{c@convention}{}{}
\@ifundefined{c@thm}{}{}
\@ifundefined{c@innercustomthm}{}{}
\makeatother

\crefname{definition}{Definition}{Definitions}
\Crefname{definition}{Definition}{Definitions}
\crefname{theorem}{Theorem}{Theorems}
\Crefname{theorem}{Theorem}{Theorems}
\crefname{lemma}{Lemma}{Lemmas}
\Crefname{lemma}{Lemma}{Lemmas}
\crefname{prop}{Proposition}{Propositions}
\Crefname{prop}{Proposition}{Propositions}
\crefname{corollary}{Corollary}{Corollaries}
\Crefname{corollary}{Corollary}{Corollaries}
\crefname{remark}{Remark}{Remarks}
\Crefname{remark}{Remark}{Remarks}
\crefname{example}{Example}{Examples}
\Crefname{example}{Example}{Examples}
\crefname{assumption}{Assumption}{Assumptions}
\Crefname{assumption}{Assumption}{Assumptions}
\crefname{problem}{Problem}{Problems}
\Crefname{problem}{Problem}{Problems}
\crefname{fact}{Fact}{Facts}
\Crefname{fact}{Fact}{Facts}
\crefname{conventions}{Conventions}{Conventions}
\Crefname{conventions}{Conventions}{Conventions}
\crefname{question}{Question}{Questions}
\Crefname{question}{Question}{Questions}
\crefname{situation}{Situation}{Situations}
\Crefname{situation}{Situation}{Situations}
\crefname{openproblem}{Open Problem}{Open Problems}
\Crefname{openproblem}{Open Problem}{Open Problems}
\crefname{notation}{Notation}{Notations}
\Crefname{notation}{Notation}{Notations}
\crefname{conjecture}{Conjecture}{Conjectures}
\Crefname{conjecture}{Conjecture}{Conjectures}
\crefname{claim}{Claim}{Claims}
\Crefname{claim}{Claim}{Claims}
\crefname{observation}{Observation}{Observations}
\Crefname{observation}{Observation}{Observations}
\crefname{exercise}{Exercise}{Exercises}
\Crefname{exercise}{Exercise}{Exercises}
\crefname{convention}{Convention}{Conventions}
\Crefname{convention}{Convention}{Conventions}
\crefname{construction}{Construction}{Constructions}
\Crefname{construction}{Construction}{Constructions}
\crefname{thm}{Theorem}{Theorems}
\Crefname{thm}{Theorem}{Theorems}
\crefname{innercustomthm}{Theorem}{Theorems}
\Crefname{innercustomthm}{Theorem}{Theorems}

\newcommand{\SE}{{\mathscr E}}
\newcommand{\SF}{{\mathscr F}}

\newcommand{\SI}{{\mathscr I}}

\newcommand{\SK}{{\mathscr K}}
\newcommand{\SL}{{\mathscr L}}

\newcommand{\SQ}{{\mathscr Q}}
\newcommand{\SR}{{\mathscr R}}
\newcommand{\ST}{{\mathscr T}}

\newcommand{\SV}{{\mathscr V}}

\newcommand{\BA}{{\mathbb{A}}}

\newcommand{\BC}{{\mathbb{C}}}

\newcommand{\BE}{{\mathbb{E}}}

\newcommand{\BG}{{\mathbb{G}}}

\newcommand{\BI}{{\mathbb{I}}}

\newcommand{\BL}{{\mathbb{L}}}
\newcommand{\BM}{{\mathbb{M}}}

\newcommand{\BP}{{\mathbb{P}}}
\newcommand{\BQ}{{\mathbb{Q}}}

\newcommand{\BZ}{{\mathbb{Z}}}

\newcommand{\CA}{{\mathcal A}}

\newcommand{\CD}{{\mathcal D}}

\newcommand{\CH}{{\mathcal H}}

\newcommand{\CJ}{{\mathcal J}}
\newcommand{\CK}{{\mathcal K}}
\newcommand{\CL}{{\mathcal L}}
\newcommand{\CM}{{\mathcal M}}
\newcommand{\CN}{{\mathcal N}}

\newcommand{\CS}{{\mathcal S}}
\newcommand{\CT}{{\mathcal T}}

\newcommand{\Fm}{{\mathfrak{m}}}

\newcommand{\FM}{{\mathfrak{M}}}

\title[Moduli spaces of rank 2 sheaves on $\mathbb P^3$ via Quot schemes]
{Moduli spaces of rank 2 sheaves on $\mathbb P^3$ via Quot schemes}

\author{Charles Almeida, Nicol\`o Bignami, Marcos Jardim, Andrea T. Ricolfi}

\begin{document}

\begin{abstract}
We study the Gieseker--Maruyama moduli space $M_{c_3}$ of stable rank 2 torsion-free sheaves on $\mathbb{P}^3$ with Chern classes $c_1=-1$, $c_2=1$, and odd $c_3 \le 1$. Every such sheaf $E$ has reflexive hull $E^{\ast\ast}$ parametrised by $M_1$, and the quotient $E^{\ast\ast}/E$ is a 0-dimensional sheaf of length $n = (1-c_3)/2$. We prove that $M_{c_3}$ is isomorphic to the relative Quot scheme $\Quot_{\BP^3 \times M_1/M_1}(\SR,n)$, where $\SR$ is the universal reflexive sheaf over $M_1$. Our technique combines (i) Bridgeland-type wall-crossing in Donaldson--Thomas theory, (ii) the theory of good moduli spaces, (iii) deformation theory in dgla and Ext-quiver language. We deduce that $M_{c_3}$ is connected, and we analyse its \emph{smoothable component} $M_{c_3}^{\main}$, proving it is rational and generically smooth of dimension $4n+3$. We further prove that $\dim \Sing M_{c_3} \geq 4n-2$ for all $n \geq 2$. Finally, we provide a geometric characterisation for the cases $n\leq 2$. We show that $M_{-1}$ is a smooth rational projective variety of dimension 7 and compute its motive in the Grothendieck ring of varieties, while $M_{-3}$ is shown to be an integral, normal variety of dimension 11 that is regular in codimension 4.
\end{abstract}

\maketitle

{\hypersetup{linkcolor=black}
\tableofcontents}

\section{Introduction}
\label{sec:Introduction}

\subsection{Overview}
\label{sec:overview}
Let $(Y,H)$ be a polarised smooth complex projective variety of dimension $d$, let $c = (r,c_1,\ldots,c_d)$ be a tuple in $\BZ_{\geq 1} \times \HH^{2\ast}(Y,\BZ)$, and denote by $M_Y(c)$ the Gieseker--Maruyama moduli space of $H$-semistable sheaves on $Y$ with rank $r$ and Chern classes $c_1,\ldots,c_d$. We denote by $M_Y^{\st}(c)\subset M_Y(c)$ the open subscheme parametrising $H$-stable sheaves. A classical result of Maruyama established that $M_Y(c)$ is a projective scheme \cite{Maruyama1977}. These spaces have been studied extensively when $d\leq 2$ --- we refer to \cite[Part II]{modulisheaves} for references to the surface case. However, when $\dim Y \geq 3$, the global geometry of such schemes remains largely unknown. In fact, Vakil's Murphy's Law shows that moduli spaces of sheaves can exhibit arbitrarily bad local singularities \cite{Murphy}. Fundamental questions regarding connectedness, smoothness, irreducibility, the number of irreducible components remain open, even when we restrict our attention to the 3-dimensional projective space.

For instance, when $Y = \mathbb{P}^3$, the moduli space $M_{\BP^3}^{\st}(r, c_1, c_2, c_3)$ may possess intricate geometric properties. For rank $r=1$ and $c_1=0$, it is isomorphic to the Hilbert scheme of closed subschemes $Z\subset \BP^3$ of dimension at most 1. The Hilbert scheme is always connected \cite{MR0213368} but also almost always singular, with an unknown number of components.

For rank $r \ge 2$, the situation is notoriously even more involved. Therefore, it is quite remarkable when one can provide a detailed description of $M_{\BP^3}^{\st}(r, c_1, c_2, c_3)$ for certain values of $(c_1, c_2, c_3)$. An account of the previously studied cases will be given in \Cref{sec:related-work}. Despite those isolated successes, a systematic global understanding of these moduli spaces remained out of reach. The main goal of this paper is to build towards filling this gap. We focus on an infinite family of fine moduli spaces of rank 2 torsion-free sheaves on $\BP^3$, namely the spaces 
\[
M_{c_3}\defeq M_{\BP^3}^{\st}(2,-1,1,c_3)
\]
of Gieseker stable rank 2 torsion-free sheaves with odd $c_3 \leq 1$ (for $c_3$ even or larger than $1$, these spaces are empty). Remarkably, the space $N\defeq M_1$ consists entirely of \emph{reflexive} sheaves and parametrises the reflexive hulls $E^{\ast\ast}$ of the sheaves $E$ parametrised by $M_{c_3}$ for all $c_3 \leq -1$. 

Our key observation is that every torsion-free sheaf $[E] \in M_{c_3}$ not only determines, but  is determined by, the canonical short exact sequence
\[
\begin{tikzcd}
0 \arrow{r} & E \arrow{r} & E^{\ast\ast} \arrow{r} & Q \arrow{r} & 0,
\end{tikzcd}
\]
where the reflexive hull $E^{\ast\ast}$ defines a point of $N$, and the quotient $Q = E^{\ast\ast}/E$ is a $0$-dimensional\footnote{The 0-dimensionality follows by $c_2=1$, which is the minimal value for which there exist stable torsion-free sheaves of rank 2 with $c_1 = -1$, cf.~\cite[Prop.~3.1, Ex.~3.7]{Gholampour2017a}.} coherent sheaf of length $n = \frac{1-c_3}{2}$. A considerable technical effort in this paper is devoted to establishing this correspondence \emph{in families}, cf.~\Cref{MAIN-THM-a}.

Our main techniques combine the Bridgeland type DT/PT wall-crossing machinery introduced in \cite{JLMM25} with deformation theory (in particular, differential graded Lie algebras, Ext-quivers and the local structure of good moduli spaces) to scheme-theoretically identify $M_{c_3}$ with the relative Quot scheme of points 
\[
\Quot_{\BP^3 \times N/N}(\SR,n),
\]
where $\SR\in \Coh(\BP^3 \times N)$ is the universal reflexive sheaf. We also prove that the structure morphism $M_{c_3} \to N$ is a Zariski locally trivial fibration, so that all Quot schemes of points $\Quot_{\BP^3}(R,n)$ for $[R] \in N$ are isomorphic to each other, and we confirm that $M_{c_3}$ is \emph{connected}.

\subsection{Main results}
\label{sec:main-results}
Our main results establish a comprehensive dictionary between the moduli spaces $M_{c_3}$ and relative Quot schemes, providing global structural properties, describing the smoothable components, and completely characterising the geometry for small values of $n$. We formalise these achievements through three main results, presented below as Theorems \ref{MAIN-THM-a}, \ref{MAIN-gen-smooth-cpt} and \ref{thm:MAIN-small-n}.

\subsubsection{\texorpdfstring{Structural properties of $M_{c_3}$}{}}

Our first main result is the following. 

\begin{thm}
\label{MAIN-THM-a}
Let $M_{c_3}$ be the fine moduli space of Gieseker stable rank 2 torsion-free sheaves on $\BP^3$ with Chern classes $(c_1,c_2) = (-1,1)$ and odd $c_3 \in \BZ_{\leq 1}$.
Let $N\subset M_1$ be the open subscheme of reflexive sheaves. The following statements hold true.
\begin{enumerate}
\item [\mylabel{MTHM1}{\normalfont{(1)}}] The inclusion $N\subset M_1$ is an equality, i.e.~$N=M_1$, and there is an isomorphism of schemes
\[
\begin{tikzcd}
N \arrow{r}{\sim} & \BP^3.
\end{tikzcd}
\]
\item [\mylabel{MTHM2}{\normalfont{(2)}}] Let $\SR \in \Coh(\BP^3 \times N)$ be the universal reflexive sheaf. Set $n = \frac{1-c_3}{2}$.
There is an isomorphism of schemes
\[
\begin{tikzcd}
\Quot_{\BP^3 \times N/N}(\SR,n) \arrow{r}{\sim} & M_{c_3}.
\end{tikzcd}
\] 
\item [\mylabel{MTHM31}{\normalfont{(3)}}] The induced natural map $\pi\colon M_{c_3} \to N$ is Zariski locally trivial. In particular, 
\begin{enumerate}
    \item [\mylabel{same-fibres-intro}{\normalfont{(i)}}] the morphism $\pi$ is flat, and all the fibres are isomorphic to each other, and
    \item [\mylabel{irr-iff-irr-intro}{\normalfont{(ii)}}] the moduli space $M_{c_3}$, where $c_3=1-2n$, is irreducible if and only if $\Quot_{\BP^3}(R,n)$ is irreducible, where $\Quot_{\BP^3}(R,n)$ is an arbitrary fibre of $\pi$.
\end{enumerate}
\item [\mylabel{MTHM3}{\normalfont{(4)}}] The moduli space $M_{c_3}$ is connected.
\end{enumerate}
\end{thm}

Part \ref{MTHM1} is proved in \Cref{prop:reflexive=P^3}. Part \ref{MTHM2} is proved in \Cref{thm:iso_relative_quot_M} (but relies on a technical result proved in \Cref{app:Kuranishi}). Part \ref{MTHM31} is proved in \Cref{prop:pi-is-ZLT}. Part \ref{MTHM3} is proven in \Cref{cor:connectedess-of-M}.

We briefly outline the proof of the isomorphism in \Cref{MAIN-THM-a}\ref{MTHM2}, as this is the technical core of the paper. 
\begin{enumerate}
\item [\mylabel{DTwc}{(a)}] The first idea is to exploit (part of) the DT/PT wall-crossing that was established, remarkably for the first time in the language of Bridgeland stability, in the paper \cite{JLMM25}. The DT side of the DT/PT wall-crossing yields a morphism 
\[
\begin{tikzcd}
    M_{c_3} \arrow{r}{\Phi} & \BM_\sigma, \qquad E \arrow[mapsto]{r} & {[}E^{\ast\ast}\oplus E^{\ast\ast}/E[-1]{]},
\end{tikzcd}
\]
where $\BM_\sigma$ is the good moduli space of Bridgeland's moduli stack $\FM_\sigma$ of $\sigma$-semistable objects in class $(2,-1,1,c_3)$, for a suitable stability condition $\sigma \in \Stab(\BP^3)$ lying on the DT/PT wall. The stack $\FM_\sigma$ is open in Lieblich's moduli stack $\FM^{\pug}(\BP^3)$ of semi-Schur objects. By \cite[Thm.~1.1]{Mihai-Tuomas}, we know that $\BM_\sigma$ is a projective \emph{scheme}.
\item [\mylabel{map-GMS}{(b)}] Denote by $\CN$ the \emph{stack} of sheaves parametrised by $N$. We construct in \Cref{thm:iso-with-sim} a 1-morphism of algebraic stacks $\overline u \colon \CN \times \CCoh^n(\BP^3) \to \FM_\sigma$ sending $(R,Q) \mapsto E=R\oplus Q[-1]$, so that taking good moduli spaces we obtain a morphism of projective schemes 
\[
\begin{tikzcd}
N \times \Sym^n(\BP^3) \arrow{r}{u} & \BM_\sigma, \qquad (R,[Q]) \mapsto [E].  
\end{tikzcd}
\]
Its bijectivity is settled by the general results of \cite[Sec.~6]{Mihai-Tuomas}. 
\item [\mylabel{GMS}{(c)}] We state, also in \Cref{thm:iso-with-sim}, that $u$ is an isomorphism. Being bijective, it is enough to prove it is \'etale. In fact, this is the main technical contribution of the paper, whose proof (confined in Appendix \ref{app:Kuranishi} to keep the flow going) is of independent interest. To prove \'etaleness of $u$ we exploit: (i) the \'etale local structure of good moduli spaces as in \cite{Alper-local-properties,Chen-Pertusi-Zhao}, to identify $\widehat{\OO}_{\BM_\sigma,[E]}$ with the $\Aut(E)$-invariant part of the miniversal formal Kuranishi algebra $\mathsf A_E$, (ii) deformation theory via dgla's \cite{Manetti-Book}, to select a convenient Kuranishi slice around $E$ yielding a convenient presentation of $\mathsf A_E$, (iii) the theory of Ext-quivers applied to $E$ and (iv) the Le Bruyn--Procesi theorem \cite{LeBruyn-Procesi,craw-yamagishi} to identify the invariant parts in $\mathsf A_\bullet$, with $\bullet \in \set{R,Q,E}$. Under these identifications, we deduce that the canonical map $\widehat u^\#_{(R,[Q])}$ is an isomorphism
\[
\begin{tikzcd}
\mathsf A_E^{\Aut(E)} \arrow{r}{\sim} & \mathsf A_R \widehat{\otimes} \mathsf A_Q^{\Aut(Q)},
\end{tikzcd}
\]
giving \'etaleness of $u$ at closed points, as we wanted.
\item [\mylabel{Double-dual}{(d)}] Once the isomorphism $N\times \Sym^n(\BP^3) \cong \BM_\sigma$ is established, the composition
\begin{equation}
\label{map-delta}
\begin{tikzcd}
    \delta\colon M_{c_3} \arrow{r}{\Phi} & \BM_\sigma \arrow{r}{\sim} & N \times \Sym^{n}(\BP^3)\arrow{r}{\pr_1} & N
\end{tikzcd}
\end{equation}
yields an algebraic \emph{morphism} which on points takes a sheaf to its double dual.\footnote{In general, as is well-known, taking the double dual of a flat family of sheaves does not yield a flat family of reflexive sheaves, cf.~\Cref{rmk:dd-issues}. However, remarkably, having the map $\delta$ turns out to be enough for our purposes.}
\item [\mylabel{map-quot}{(e)}] We construct in \Cref{sec:iso-quot-M} a bijective morphism  
\[
\begin{tikzcd}
\Quot_{\BP^3 \times N/N}(\mathscr R,n) \arrow{r}{f_{c_3}} & M_{c_3}
\end{tikzcd}
\]
of $N$-schemes. We then exploit the map $\delta$ to construct an explicit inverse (cf.~\Cref{thm:iso_relative_quot_M}). After this isomorphism is proved, by a slight abuse of notation we identify $\delta$ with the structure map $\pi\colon \Quot_{\BP^3 \times N/N}(\mathscr R,n) \to N$, cf.~\Cref{MAIN-THM-a}\ref{MTHM2}.
\end{enumerate}

\subsubsection{\texorpdfstring{The smoothable component of $M_{c_3}$}{}}
Parallel to the case of the Hilbert scheme of points $\Hilb^n(\BP^3)$, we observe that $M_{c_3}$ has a distinguished irreducible component $M_{c_3}^{\main} \subset M_{c_3}$, which, for the sake of analogy, we call the \emph{smoothable component}. It is defined as the closure of the (irreducible, open) locus
\[
M_{c_3}^{\circ} = \Set{[E] \in M_{c_3} |
\begin{array}{cc}
E^{\ast\ast} / E = \OO_Z \mbox{ where}   &  \\
 Z\subset \BP^3 \mbox{ is reduced}
\end{array}
}
\subset M_{c_3},
\]
and we endow it with its reduced induced structure.
The next result gives geometric information about this main component.

\begin{thm}
\label{MAIN-gen-smooth-cpt}
Fix an odd $c_3 \in \BZ_{\leq -1}$, and set $n = \frac{1-c_3}{2}$.
The smoothable component $M_{c_3}^{\main} \subset M_{c_3}$ has the following properties. 
\begin{itemize}
\item [\mylabel{sm-cpt-i}{\normalfont{(1)}}] It is generically smooth of dimension $5-2c_3 = 4n+3$. 
\item [\mylabel{sm-cpt-ii}{\normalfont{(2)}}] Assume $n\geq 2$. Then
\[
\dim \Sing M_{c_3} \geq \dim M_{c_3}^{\main} - 5 = 4n-2.
\]
\item [\mylabel{sm-cpt-iii}{\normalfont{(3)}}] It is a rational variety. 
\end{itemize}
\end{thm}

Part \ref{sm-cpt-i} is \Cref{mainComponent}. Part \ref{sm-cpt-ii} is \Cref{Sing-main-Component}. Part \ref{sm-cpt-iii} is proved in \Cref{sec:rationality-main-component}.

\subsubsection{\texorpdfstring{Description of $M_{c_3}$ for the largest values of $c_3$}{}}
For large $c_3$ (i.e.~low $n$), we can give an explicit description of the full moduli space $M_{c_3}$. The first nontrivial cases, beyond the description in \Cref{MAIN-THM-a}\ref{MTHM1}, are $c_3 \in \Set{-1,-3}$, corresponding to $n \in \Set{1,2}$. 

\begin{thm}
\label{thm:MAIN-small-n}
The following statements hold true.
\begin{enumerate}
\item [\mylabel{M-SMALL1}{\normalfont{(1)}}] The moduli space 
$M_{-1}$ is a smooth projective variety of dimension 7. Its motive in the Grothendieck ring of varieties $K_0(\Var_{\BC})$ is given by the formula
\[
[M_{-1}] = 1 + 3 \BL + 6 \BL^2 + 8 \BL^3 + 8 \BL^4 + 6 \BL^5 + 3 \BL^6 + \BL^7.
\]
\item [\mylabel{M-SMALL2}{\normalfont{(2)}}] The moduli space $M_{-3}$ is an integral, 11-dimensional projective variety, which is Cohen--Macaulay, regular in codimension 4, and normal.
\end{enumerate}
\end{thm}

Part \ref{M-SMALL1} is proved in \Cref{sec:M-1}, part \ref{M-SMALL2} is proved in \Cref{sec:M-3}, exploiting also \Cref{prop:int+CM+8}, where we prove that $\Quot_{\BP^3}(R,2)$ is integral and Cohen--Macaulay for every reflexive sheaf $R$ parametrised by $N$, using a general zero locus description of $\Quot_{\BP^3}(R,n)$, cf.~\Cref{lemma:zero-locus-section}.

\subsubsection{Side comment on \Cref{MAIN-THM-a}\ref{MTHM2}}
The isomorphism 
\[
\begin{tikzcd}
f\colon \Quot_{\BP^3 \times N/N}(\SR,n) \arrow{r}{\sim} & M_{c_3}
\end{tikzcd}
\] 
of \Cref{MAIN-THM-a}\ref{MTHM2} is obtained through highly nontrivial machinery, strengthening the already solid intuition that the gap in complexity between moduli spaces of rank 1 sheaves and rank 2 sheaves is substantial. The (simplest) parallel result for $r=1$ is the existence of an isomorphism 
\[
\begin{tikzcd}
\Hilb^n(Y) \arrow{r}{\sim} & M^{\st}_Y(1,\OO_Y,0,-n) 
\end{tikzcd}
\]
between the Hilbert scheme of points and the moduli space of ideal sheaves (torsion-free sheaves of rank 1 with trivial determinant), for any smooth projective 3-fold with $\HH^1(Y,\OO_Y) = \HH^2(Y,\OO_Y) = 0$. The map sends a quotient $[\OO_Y \onto \OO_Z]$ to its kernel $\mathscr I_Z$ and is bijective on points. This part is identical to how we obtain the \emph{map} $f$, but in the case of $\Hilb^n(Y) \to M^{\st}_Y(1,\OO_Y,0,-n)$, to conclude it is an isomorphism one simply exploits $\HH^1(Y,\OO_Y) = 0$ to deduce that the induced map on tangent spaces is an isomorphism, and $\HH^2(Y,\OO_Y) = 0$ to deduce that the induced map on obstruction spaces is injective. By \cite[Prop.~2.2]{Cazzaniga-Ricolfi-Framedsheaves}, these properties alone prove \'etaleness, hence the two spaces are isomorphic, see e.g.~\cite[Thm.~11.2.11]{Ricolfi-Book} for the full argument.\footnote{An alternative argument for this isomorphism can be found in \cite[Lemma B.5.6]{Kuznetsov2018}.} In our case, the deformation-theoretic technology one needs to exploit is much more sophisticated.

\subsection{Related work}
\label{sec:related-work}
Our perspective closely parallels \cite{Guimaraes-Jardim}, where the moduli space of \emph{quasi-trivial} sheaves $\mathcal{N}_Y(r,n)$, defined as sheaves whose reflexive hulls are completely trivial, $E^{\ast\ast} \cong \OO_Y^{\oplus r}$, is investigated via the Quot scheme of points $\Quot_Y(\OO_Y^{\oplus r}, n)$. In the setting of this paper, the reflexive hulls \emph{vary} over the fine moduli space $N = M_1$ of stable rank 2 reflexive sheaves.

Previous investigations established nonemptiness bounds for $M_{\BP^3}(2,c_1,c_2,c_3)$ when $c_1 \in \{0, -1\}$, proving irreducibility for $c_1=c_2=0$ in the range $-20 \le c_3 \le -4$, see \cite{Guimaraes-Jardim}, as well as smoothness, irreducibility, and rationality at the maximal boundaries $c_3 = c_2^2 - c_2 + 2$ ($c_1=0$) and $c_3 = c_2^2$ ($c_1=-1$), see \cite{Schmidt}. Furthermore, $M_{\BP^3}(2,-1,c_2,c_2^2-2)$ was shown to be irreducible of dimension $c_2^2+3c_2+5$ in \cite{AJO} when $c_2\geq 4$, while the connectedness and precise number of irreducible components (and their dimensions) was determined for $M_{\BP^3}(2,0,2,c_3)$ in \cite{JMT}, for $M_{\BP^3}(2,-1,2,c_3)$ in \cite{AJT2022}, and for $M_{\BP^3}(2,-1,3,7)$ in \cite{AJO}.

On the other hand, these moduli spaces are also extremely relevant in enumerative geometry: they are the key players in (rank 2) Donaldson--Thomas theory. In this context, the available techniques are powerful enough to define, and in some cases compute, sensible invariants attached to these moduli spaces, despite the lack of a general understanding of their global geometry. For instance, a close relationship between $M_Y(r,c_1,c_2,c_3)$ and Quot schemes of points on a toric 3-fold $Y$ was established by Gholampour--Kool--Young in \cite{Gholampour2017a}, where classical and virtual enumerative invariants are discussed, partially at a conjectural level, in the combinatorial language of torus localisation: a vertex formalism is proposed in loc.~cit.~and proved (sometimes up to a conjecture) to yield natural formulas displaying a rich structure. In \cite{Gholampour-Kool-stablereflexive}, Gholampour--Kool prove several formulas and structural results for generating functions of Euler characteristics of moduli spaces of stable rank 2 reflexive sheaves on $\BP^3$. In \cite{Gholampour2017} they obtain, partially building on earlier work on the Serre correspondence \cite{Gholampour-Kool-Serre-corr}, beautiful formulas --- functional equations --- for Euler characteristics of moduli spaces of torsion-free sheaves and Quot schemes of points on 3-folds, establishing with Hall algebra and wall-crossing techniques a link with Pandharipande--Thomas stable pairs; this latter point of view, originated in rank 1 in \cite{ST}, was later picked up again and upgraded to the virtual level in higher rank in \cite{Beentjes-Ricolfi}. The role of Quot schemes in the wall-crossing theory for 3-folds is also present in the works \cite{Lo2018,Ricolfi-DT/PT-smooth-curves,Davison-Ricolfi,Ricolfi-Local-contributions-DT,Cazzaniga-Ricolfi-crepant,Cazzaniga-Ralaivaosaona-Ricolfi}.

\subsection{Organisation of contents}
\label{sec:organisation-contents}
This paper is organised as follows.

In \Cref{sec:Background}, we collect foundational material on moduli spaces of sheaves. In \Cref{sec:quot-of-points}, we collect several results on Quot schemes of points, including a general connectedness result (\Cref{prop:Quot-connected}), to our knowledge not yet present in the literature. In \Cref{sec:sheaves}, we introduce the moduli space $M_{c_3}$ and prove some basic properties, such as the existence of a universal family (\Cref{lemma:M_is_fine}) and the isomorphism $N\defeq M_1\cong \BP^3$ (\Cref{prop:reflexive=P^3}), corresponding to \Cref{MAIN-THM-a}\ref{MTHM1}. In \Cref{sec:Cohomological} we prove some properties of global and local self-Ext's attached to points $[E] \in M_{c_3}$. For instance, in \Cref{lem:dimext2} we prove that $\ext^2(E,E)=4n$ whenever $E^{\ast\ast}/E$ is supported on $n = \frac{1-c_3}{2}$ distinct points.
In \Cref{sec:proofmainthma} we prove \Cref{MAIN-THM-a}\ref{MTHM2}. In \Cref{sec:moduli-wall} we analyse the DT side of the wall-crossing picture, which is one of the technical cores of the paper. Our key result, \Cref{thm:iso-with-sim}, establishes that the good moduli space $\BM_\sigma$ of Bridgeland's moduli stack $\FM_\sigma$ at the DT/PT wall is isomorphic to $N \times \Sym^n(\BP^3)$ --- details are given in \Cref{app:Kuranishi}. In \Cref{sec:proofmainthmb} we prove \Cref{MAIN-gen-smooth-cpt}. We prove in \Cref{sec:zerolocus} a vanishing locus description of $\Quot_{\BP^3}(R,n)$, and the generic smoothness of its smoothable component as part of \Cref{prop:Quot-generic-smoothness}. In \Cref{sec:maincomponent} we introduce the main component and prove Items \ref{sm-cpt-i} and \ref{sm-cpt-ii} of \Cref{MAIN-gen-smooth-cpt}, whereas Item \ref{sm-cpt-iii} is proved in \Cref{sec:rationality-main-component}.

In \Cref{sec:M-1} and \Cref{sec:M-3} we provide a geometric description of the geometry of the moduli spaces $M_{-1}$ and $M_{-3}$, which completes the proof of \Cref{thm:MAIN-small-n}. In \Cref{sec:BN locus}, we introduce and investigate the \emph{Brill--Noether locus} $B_{c_3} \subset M_{c_3},$ parametrising sheaves $E$ with $h^0(\BP^3, E(1)) > 0$. Using the modular Serre correspondence \cite{Jardim-Mu}, we prove that $B_{c_3}$ is always connected, and irreducible when $c_3 \in \{-1, -3, -5\}$. Additionally, we prove that $B_{c_3}=M_{c_3}$ for $c_3 \in \{-1,-3\}$.

\begin{conventions}
\label{conventions}
We work over the field $\BC$ of complex numbers throughout, unless specified otherwise.
We denote by $\Coh X$ the abelian category of coherent sheaves on a scheme $X$. A \textit{variety} is an integral, quasiprojective scheme over $\BC$.
If $f\colon X \to S$ is a morphism, $\mathscr E \in \Coh X$ is a coherent sheaf (resp.~a complex of coherent sheaves) and $s \in S$ is a point, we let $\SE_s$, or sometimes $\SE|_s$, denote the sheaf $\mathscr E|_{X_s}$ (resp.~the derived restriction of $\SE$ to the fibre $X_s$).
Dimensions of cohomology groups $\HH^i(Y,F)$ are denoted $h^i(Y,F)$, and similarly we set $\ext^i(F,G)\defeq \dim_{\BC}\Ext^i(F,G)$.
We let $\derived^\ast(X)$ denote the derived category of $\OO_X$-modules with coherent cohomology sheaves, possibly nonvanishing in the interval described by the superscript $\ast$, and similarly for the categories $\Perf^\ast(X)$ of perfect complexes. When dealing with a fine moduli space of sheaves $M$, we usually denote by $[E]$ the closed point corresponding to the sheaf $E$. When this may create a clash of notation with the S-equivalence class, also denoted $[E]$, we write $E \in M$ instead.
\end{conventions}

\subsection*{Acknowledgements}
A.R.~wishes to thank Martijn Kool for introducing him to the problem of understanding the moduli spaces studied in this paper and, in particular, for suggesting the statement of \Cref{MAIN-THM-a}\ref{MTHM2}, as well as for the countless generous explanations over the years around the subject of this paper. We wish to thank Mihai Pavel and Tuomas Tajakka for offering very useful clarifications on their work. We wish to thank Okke van Garderen for very helpful discussions around Kuranishi theory and Ian Selvaggi for valuable conversations on good moduli spaces.

M.J.~is partially supported by the CNPQ grant number 305601/2022-9 and the FAPESP-CEPID Project 2024/00923-6; he would also like to acknowledge support from the ICTP through the Associates Programme and from the Simons Foundation, whose grant (Record ID: SFI-MPS-T-Institutes-00012057, AD) also supported this work. 

\subsection*{AI declaration}
LLMs were used for (i) a preliminary proofreading of an old draft, (ii) mathematical assistance in the case-by-case analysis needed for the proof of \Cref{cor:irreducible-quot}, and (iii) to find the shape of the Kuranishi relations appearing in \Cref{app:Kuranishi}. In the latter case, we independently found a direct proof of these relations via Ext-quivers, not relying on the AI ansatz. All ideas present in this paper are the authors' own. The authors take full responsibility for the scientific content and conclusions of this work.

\section{Background material}
\label{sec:Background}

Let $(Y,H)$ be a polarised smooth complex projective variety of dimension $d$. We refer to \cite[Sec.~4]{modulisheaves} for the construction of the moduli space of (semi)stable torsion-free sheaves $M_Y^{\st}(v)\subset M_Y(v)$ with fixed class $v \in \HH^\ast(Y,\BQ)$, or fixed rank and Chern classes. For background on stable pair moduli spaces, which we only use in \Cref{sec:BN locus}, we refer to \cite{Moduli-stable-pairs}. As for the various notions of stability on a torsion-free sheaf on a smooth projective 3-fold, we have implications 
\begin{equation}
\label{stab-notions}
\begin{split}
\mu\mbox{-stable}&\,\,\Rightarrow\,\,\mbox{2-Gieseker stable}\,\,\Rightarrow\,\,\mbox{Gieseker stable}\,\,\Rightarrow\,\,\mbox{Gieseker semistable}\\
&\,\,\Rightarrow\,\,\mbox{2-Gieseker semistable}
\Rightarrow\,\,\mu\mbox{-semistable},
\end{split}
\end{equation}
and all these notions coincide as soon as $\gcd(\rk E,c_1(E)\cdot H^{2})=1$, cf.~\cite[Sec.~2]{Mihai-Tuomas}. This will be exactly our situation, so the distinctions do not matter. Moreover, for moduli spaces of sheaves on $Y=\BP^d$, fixing the Chern classes and the rank is equivalent to fixing the Chern character, which in turn is equivalent to fixing the Hilbert polynomial \cite{zbMATH02138866}. We will use this fact without mention when studying sheaves on $Y=\BP^3$, so we will write $M_Y(r,c_1,c_2,c_3)$ and $M_Y(v)$ interchangeably.

\begin{remark}
\label{rmk:moduli_of_0-cycles}
There is a notion of semistability for torsion sheaves \cite[Def.~1.6.9]{modulisheaves}, which is automatically satisfied by every 0-\emph{dimensional} sheaf; moreover, for every $n \in \BZ_{\geq 1}$, there is a moduli space $M_Y(0,\ldots,0,n)$ of 0-dimensional (semistable) sheaves $Q\in \Coh Y$ of length $n=\chi(Q)=h^0(Y,Q)$. As proven in \cite[Ex.~4.3.6]{modulisheaves}, it is isomorphic to the symmetric product $\Sym^n(Y) = Y^n/S_n$, which can be seen as the moduli space of effective $0$-cycles of degree $n$, see \cite{Rydh1}. It can also be seen as the \emph{good moduli space} of the algebraic stack $\CCoh^n(Y)$ of 0-dimensional length $n$ sheaves. The good moduli space morphism $\gamma_n\colon \CCoh^n(Y) \to \Sym^n(Y)$ is the ``Coh-to-Chow morphism'', see \cite{Rydh1,Fantechi-Ricolfi-structural,Fantechi-Ricolfi-motivic}.
\end{remark}

Recall that a torsion-free sheaf $E$ is called
\emph{reflexive} if the canonical map $\theta_E \colon E \to E^{\ast\ast}$ is an isomorphism. The sheaf $E^{\ast\ast}$ is called the \emph{reflexive hull} of $E$. In general, a reflexive sheaf $E$ is locally free away from a closed subset $\Sing E \subset Y$ of codimension at least $3$, called the \emph{singularity set} of $E$.
We refer the reader to \cite{Hart_StableReflexive} for the general theory of reflexive sheaves and their moduli.

Fix a vector $c = (r,c_1,\ldots,c_d)$, where $r \geq 1$ is a positive integer and $c_i \in \HH^{2i}(Y,\BZ)$ for $i=1,\ldots,d$. There is an open subscheme (possibly empty)
\begin{equation}
\label{def:M-reflexive}
\begin{tikzcd}
N_Y^{\mathsf{ref}}(c) \arrow[hook]{r} & M_Y^{\st}(c)
\end{tikzcd}
\end{equation}
parametrising stable reflexive sheaves in class $c$: it will be referred to as `the moduli space of reflexive sheaves'.

While a torsion-free sheaf is $\mu$-stable if and only if its double dual is, the double dual $\SE^{\ast\ast}$ of a flat family $\SE$ of stable sheaves is not a \emph{flat} family of reflexive sheaves in general. In other words, assume for simplicity that all stabilities in \eqref{stab-notions} coincide; then, the set-theoretic map
\[
\begin{tikzcd}
 M_Y^{\st}(c) \arrow{r}{\mathsf d} &  \displaystyle\coprod_{c_2',\ldots,c_d'}N_Y^{\mathsf{ref}}(r,c_1,c_2',\ldots,c_d')
\end{tikzcd}
\]
sending a $\mu$-stable sheaf to its double dual need not underlie an algebraic morphism in general. However, results of Koll\'{a}r \cite{Kollar_hulls} prove that $\mathsf d$ is constructible, i.e.~one can stratify $M_Y^{\st}(c)$ by locally closed subschemes $Z_i \into M_Y^{\st}(c)$ in such a way that $\mathsf d|_{Z_i}$ is an algebraic morphism. 

\begin{remark}
\label{rmk:dd-issues}
Even when there is a unique stratum $Z$, so that $Z_{\red} = M_Y^{\st}(c)_{\red}$, and the Chern classes of \emph{all} double duals are the same, it is not obvious a priori that the map $\mathsf d$ is a morphism, unless $M_Y^{\st}(c)$ is known to be reduced. This is the main issue that had to be technically overcome to obtain some of our main results.
\end{remark}

\begin{remark}
\label{rmk:refl-on-3folds}
Let $Y$ be a smooth projective variety. Recall that the \emph{homological dimension} $\hd(E)$ of a coherent sheaf $E \in \Coh Y$ is the maximal $p \in \BZ_{\geq 0}$ such that $\lExt^p(E,\OO_Y)\neq 0$. 
Let $E$ be a torsion-free sheaf on a smooth projective 3-fold $Y$ and let $R=E^{\ast\ast}$. Then $\dim R/E \leq 1$, and $\hd(E)\leq 1$ if and only if $R/E$ is either 0 or has pure dimension 1.

In fact, if $\dim Y = 3$, a torsion-free sheaf $R$ is reflexive if and only if $\dim \lExt^q(R,\OO_Y) \leq 1-q$ for all $q>0$, see e.g.~\cite[Prop.~1.1.10]{modulisheaves}. 
\end{remark}

We close this section with an elementary application of local duality, which we will need (only) for the proof of \Cref{lem:exts}\ref{ext-iii}.

\begin{lemma}
\label{lemma:supp-ext-TF}
Let $Y$ be a smooth variety of dimension $d$, and let $Q$ be a 0-dimensional coherent sheaf, $E \neq 0$ a torsion-free sheaf. Then
\[
\Supp \lExt^d(Q,E) = \Supp Q.
\]
\end{lemma}

\begin{proof}
If $y \in Y \setminus \Supp Q$, then $Q_y = 0$ and hence $\lExt^d(Q,E)_y = 0$. Conversely, assume $y \in \Supp Q$. We must prove that $\lExt^d(Q,E)_y \neq 0$. Set $R=\OO_{Y,y}$, $M = Q_y$ and $F=E_y$, so that our claim reduces to $\Ext_R^d(M,F) \neq 0$. Since $R$ is regular and $M$ has finite length, $M$ is a perfect
$R$-module of projective dimension $d$, and
\[
\Ext^i_R(M,R)=0,\qquad i\neq d.
\]
In particular, we have $\RHom_R(M,R) \cong \Ext_R^d(M,R)[-d]$. By perfectness of $M$ again, we have 
\[
\RHom_R(M,F) \cong \RHom_R(M,R) \otimes^{\mathbf L}_RF \cong \Ext_R^d(M,R)[-d] \otimes^{\mathbf L}_RF,
\]
and thus
\[
\Ext_R^d(M,F) = \HH^d(\RHom_R(M,F)) \cong \HH^0(\Ext_R^d(M,R) \otimes^{\mathbf L}_RF) \cong \Ext_R^d(M,R)\otimes_RF.
\]
We have to prove that this tensor product is nonzero. By local duality, the $R$-module $\Ext_R^d(M,R)$ is nonzero, has finite length and support equal to the support of $M$. Thus
\[
\Ext_R^d(M,R) \otimes_R R/\Fm_y \neq 0.
\]
Since $E$ is torsion free, $E_z \neq 0$ for all $z \in Y$, so in particular $F \neq 0$, and hence $F \otimes_R R/\Fm_y \neq 0$ by Nakayama's lemma. But now
\[
\Ext_R^d(M,F) \otimes_R R/\Fm_y \cong (\Ext_R^d(M,R) \otimes_R R/\Fm_y) \otimes_{R/\Fm_y} (F \otimes_R R/\Fm_y).
\]
Now we have a tensor product over a field, which vanishes if and only if one of the factors vanishes. It follows that $\Ext_R^d(M,F) \otimes_R R/\Fm_y \neq 0$, i.e.~$\Ext_R^d(M,F) \neq 0$, again by Nakayama's lemma.
\end{proof}

\section{Auxiliary results on Quot schemes of points}
\label{sec:quot-of-points}

Let $Y \to B$ be a projective morphism of noetherian schemes, $L$ a relatively ample line bundle on $Y$. Let also $F \in \Coh Y$ be a coherent sheaf, and $P \in \BQ[t]$ a polynomial. The \emph{Quot functor} attached to the tuple $(Y \to B,F,P,L)$ is the functor
\[
\begin{tikzcd}
\underline{\Quot}_{Y/B}(F,P) \colon \Sch_B^{\opp} \arrow{r} & \Sets
\end{tikzcd}
\]
sending a locally Noetherian $B$-scheme $V \to B$ to the set of equivalence classes\footnote{Two surjections are in the same equivalence class if they share the same kernel.} of surjections 
\[
\begin{tikzcd}
F_V \arrow[two heads]{r} & T
\end{tikzcd}
\]
in $\Coh (Y\times_BV)$, where $T$ is $V$-flat, $F_V$ is the pullback of $F$ along $Y\times_BV \to Y$, and for all $x \in V$ the Hilbert polynomial of $T|_x \in \Coh Y_x$ with respect to $L|_x$ is equal to $P$. Then a classical result of Grothendieck \cite[Thm.~3.1]{Grothendieck_Quot} states that $\underline{\Quot}_{Y/B}(F,P)$ is representable by a projective $B$-scheme, denoted $\Quot_{Y/B}(F,P)$, see also \cite{Nit}. We are interested in the special situation where $P$ is just an integer $n \geq 0$, in which case $\Quot_{Y/B}(F,n)$ is called the (relative) \emph{Quot scheme of points} --- and note that the polarisation $L$ becomes irrelevant in this case. Also, the Quot scheme of points still exists if $Y \to B$ is only quasiprojective \cite{Nit}. If $F$ is locally free of rank 1, then $\Quot_{Y/B}(F,n)$ recovers the classical (relative) \emph{Hilbert scheme of points} $\Hilb^n(Y/B)$, parametrising finite closed subschemes of length $n$ in the fibres of $Y \to B$. When $B=\Spec \BC$, we will omit it from the notation.

We now present, for future use, a series of results on the Quot scheme of points on a projective variety. Our first result, \Cref{prop:Quot-connected}, concerns the connectedness of the Quot scheme of points in pretty general situations, to the best of our knowledge not fully covered in the literature, cf.~\Cref{rmk:quot-connected-known}. It relies on the following elementary lemma, which is also used several other times in this paper and we therefore record here.

\begin{lemma}
\label{lemma:conn-irr}
Let $f\colon X \to Y$ be a surjective morphism of schemes.\footnote{The first two items are of course valid for arbitrary topological spaces.} 
\begin{enumerate}
    \item [\mylabel{CON}{\normalfont{(i)}}] If $Y$ is connected, each fibre $X_y$ is connected, and $f$ is either open or closed, then $X$ is connected.
    \item [\mylabel{IRR}{\normalfont{(ii)}}] If $Y$ is irreducible, all fibres $X_y$ are irreducible, and $f$ is open, then $X$ is irreducible.
    \item [\mylabel{IRR2}{\normalfont{(iii)}}] If $Y$ is irreducible, $f$ is a proper morphism of noetherian schemes, all fibres $X_y$ are irreducible of the same (finite) dimension, then $X$ is irreducible.
\end{enumerate}
\end{lemma}

\begin{proof}
Assertion \ref{CON} follows from \cite[\href{https://stacks.math.columbia.edu/tag/0378}{Tag 0378}]{stacks-project} if $f$ is open and from \cite[\href{https://stacks.math.columbia.edu/tag/0377}{Tag 0377}]{stacks-project} if $f$ is closed. Assertion \ref{IRR} follows from \cite[\href{https://stacks.math.columbia.edu/tag/037A}{Tag 037A}]{stacks-project}. Assertion \ref{IRR2} is an easy application of the upper-semicontinuity of fibre dimension \cite[\href{https://stacks.math.columbia.edu/tag/0D4I}{Tag 0D4I}]{stacks-project}.
\end{proof}

\begin{prop}
\label{prop:Quot-connected}
Let $Y$ be a projective variety\footnote{Recall from \Cref{conventions} that varieties are integral schemes, hence connected.} of positive dimension, and let $0\neq F\in\Coh Y$ be a torsion-free sheaf. Then, the Quot scheme of points $\Quot_Y(F,n)$ is connected for every $n \geq 0$.
\end{prop}

\begin{proof}
We proceed by induction on $n$. The case $n=0$ is trivial, as $\Quot_Y(F,0)=\Spec \BC$. Let us assume $\Quot_Y(F,n-1)$ is connected. Consider the nested Quot scheme $\Quot_Y(F,n,n-1)$, namely the closed subscheme 
\[
\begin{tikzcd}
\Quot_Y(F,n,n-1) \arrow[hook]{r} & \Quot_Y(F,n) \times \Quot_Y(F,n-1),
\end{tikzcd}
\]
parametrising (isomorphism classes of) length 2 flags
\[
\begin{tikzcd}
F \arrow[two heads]{r} & T_1 \arrow[two heads]{r} & T_0
\end{tikzcd}
\]
where $\chi(T_1)=n$ and $\chi(T_0)=n-1$. Consider the two projections
\[
\begin{tikzcd}[row sep=large,column sep=large]
\Quot_Y(F,n,n-1) \arrow{r}{\pr_1} \arrow[swap]{d}{\pr_2} & \Quot_Y(F,n) \\
\Quot_Y(F,n-1) & 
\end{tikzcd}
\]
and fix a closed point $q_0 = [F \onto T_0] \in \Quot_Y(F,n-1)$. Then $q_0$ is the datum of an exact sequence
\begin{equation}
\label{ses:K0}
\begin{tikzcd}
    0 \arrow{r} & K_0 \arrow{r} & F \arrow{r}{\rho} & T_0 \arrow{r} & 0
\end{tikzcd}
\end{equation}
and $\pr_2^{-1}(q_0)$ is naturally identified with the projectivisation of the sheaf $K_0$, denoted $\BP(K_0)$.
Set-theoretically, the identification works as follows: a point in $\pr_2^{-1}(q_0)$ is an equivalence class of a surjection $F \onto T_1$ fitting in a commutative diagram
\[
\begin{tikzcd}
F \arrow[two heads]{r}{\rho}\arrow[two heads]{d} & T_0 \arrow[equal]{d} \\
T_1 \arrow[two heads]{r} & T_0
\end{tikzcd}
\]
so that considering the kernel $J \defeq \ker (T_1 \onto T_0)$, a length 1 sheaf, and the sequence \eqref{ses:K0}, one gets an induced surjection $K_0 \onto J$. Conversely, given $K_0 \onto J$, define $T_1$ as the pushout of $K_0 \to F$ and $K_0 \onto J$. This yields an injection $J \into T_1$ with cokernel $T_0$. In families, a length one quotient $\CK_0 \onto \CJ$ gives
$\CT_1$ by pushout, and $\CT_1$ is flat over the test base because it is an extension of the two flat families $\CT_0$
and $\CJ$. Conversely the kernel of $\CT_1 \onto \CT_0$ is a flat family of length 1 sheaves. Now, since $F$ is torsion-free, so is its subsheaf $K_0$. Therefore, the structure morphism $\BP(K_0) \to Y$ is surjective (using $\dim Y > 0$). Since $Y$ is connected, and the fibres of the proper (hence closed) morphism $\BP(K_0) \to Y$ are connected, it follows from \Cref{lemma:conn-irr}\ref{CON} that $\BP(K_0)$, the fibre of $\pr_2$, is connected. But the target of $\pr_2$ is connected by the inductive hypothesis, and $\pr_2$ is surjective and proper (hence closed). Therefore $\Quot_Y(F,n,n-1)$ is connected, again by \Cref{lemma:conn-irr}\ref{CON}. Now, each 0-dimensional sheaf has a filtration with length 1 subquotients; in particular, given $[F \onto T_1] \in \Quot_Y(F,n)$, the sheaf $T_1$ admits a length 1 subsheaf, and hence a quotient $T_0$ of length $n-1$. This shows that $\pr_1$ is surjective. It follows that $\Quot_Y(F,n)$ is connected.
\end{proof}

\begin{remark}
\label{rmk:quot-connected-known}
\Cref{prop:Quot-connected} generalises a special case of \cite[Thm.~1.4]{Monavari-Ricolfi-motive-nested-quot}, where the statement (involving more general \emph{nested} Quot schemes) implied the connectedness of the Quot scheme of points attached to a \emph{locally free} sheaf. See also Cheah's work~\cite[Sec.~0.4]{MR1616606} proving that the \emph{nested Hilbert scheme} $\Hilb^{\boldit{n}}(Y)$ is always connected. 
\end{remark}

Next, we consider the following setup. Let $Y$ be a smooth $d$-dimensional variety, and let $\Sym^n(Y) = Y^n/S_n$ be its $n$-th symmetric product, for $n \geq 0$. Then $\Sym^n(Y)$ is a variety of dimension $nd$, singular for all $n \geq 2$ when $d>1$. Let $F$ be a coherent sheaf of rank $r\geq 1$ over $Y$. The Quot-to-Chow morphism \cite{Rydh1,Fantechi-Ricolfi-structural,Fantechi-Ricolfi-motivic}
\begin{equation}
\label{quot-to-chow}
\begin{tikzcd}
\Quot_{Y}(F,n) \arrow{r} & \Sym^n(Y)  
\end{tikzcd}
\end{equation}
is a proper morphism defined on closed points by sending a quotient $[F \onto Q]$ to the effective 0-cycle $[Q] = \sum_{y \in \Supp Q}m_yy \in \Sym^n(Y)$ determined by the points $y$ in the support of $Q$, each taken with its geometric multiplicity $m_y = \length_{\OO_{Y,y}}Q_y$. 

For the next definition, recall that the `big diagonal' in $\Sym^n(Y)$ is the closed subscheme
\[
\Delta^{\mathsf{big}} = \Set{\displaystyle\sum_{1\leq i \leq n}n_iy_i|n_i > 1\mbox{ for some }i}\subset \Sym^n(Y),
\]
and its open complement $\Sym^n(Y)\setminus \Delta^{\mathsf{big}}$ is a smooth variety (of dimension $nd$), as the quotient map $Y^n \to \Sym^n(Y)$ is a finite \'etale $S_n$-torsor when restricted to this locus.

\begin{definition}
\label{def:QUOT-circ}
Let $Y$ be a smooth variety. Fix $F \in \Coh Y$ and an integer $n \geq 0$. Define the open subscheme $\Quot^{\circ}_{Y}(F,n) \subset \Quot_Y(F,n)$ by means of the cartesian diagram
\[
\begin{tikzcd}[row sep=large,column sep=large]
    \Quot^{\circ}_{Y}(F,n)\MySymb{dr} \arrow[swap]{d}{h}\arrow[hook]{r} & \Quot_Y(F,n) \arrow{d}\\
\Sym^n(Y)\setminus \Delta^{\mathsf{big}}\arrow[hook]{r} & \Sym^n(Y)
\end{tikzcd}
\]
and define the \emph{smoothable locus} of $\Quot_Y(F,n)$ to be the closure
\[
\Quot_Y^{\main}(F,n) \defeq \overline{\Quot^{\circ}_{Y}(F,n)}  \subset \Quot_Y(F,n).
\]
\end{definition} 

Explicitly, we have the description
\[
\Quot^{\circ}_{Y}(F,n) =  \Set{[F \onto Q] \in \Quot_{Y}(F,n) | 
\begin{array}{cc}
 Q = \OO_Z \mbox{ with}   &  \\
 Z\subset Y \mbox{ reduced}
\end{array}
}.
\]

\begin{lemma}
\label{lemma:h-smooth}
Let $Y$ be a smooth variety, $F \in \Coh Y$ a locally free sheaf of rank $r\geq 1$, and let $n$ be a positive integer.
The morphism
\[
\begin{tikzcd}
\Quot^\circ_Y(F,n) \arrow{r}{h} & \Sym^n(Y) \setminus \Delta^{\mathsf{big}},
\end{tikzcd}
\]
is smooth of relative dimension $n(r-1)$.
\end{lemma}

\begin{proof}
The Quot-to-Chow morphism $\Quot_Y(F,n) \to \Sym^n(Y)$ factors through the Coh-to-Chow morphism $\gamma_n\colon \CCoh^n(Y) \to \Sym^n(Y)$. The morphism $\Quot_Y(F,n) \to \CCoh^n(Y)$, retaining the quotient sheaf, is smooth, see \cite[Lemma 1.24]{Fantechi-Ricolfi-structural} and \cite[Prop.~1.4]{Ricolfi-Savvas}. Therefore, so is its base change $\gamma_n^\circ \colon \Quot^\circ_Y(F,n) \to \CCoh^n(Y)^\circ$ along the open substack $\CCoh^n(Y)^\circ \subset \CCoh^n(Y)$ of sheaves supported on $n$ distinct points: the morphism $\gamma_n^\circ$ has fibre $(\BC^{r} \setminus 0)^n$. As a consequence of \cite[Thm.~1.38(3)]{Fantechi-Ricolfi-structural}, the morphism $\CCoh^n(Y)^\circ \to \Sym^n(Y) \setminus \Delta^{\mathsf{big}}$ is \'etale locally trivial with fibre $\mathrm{B}\BG_m^n$. It follows that $h$ is itself \'etale locally trivial, hence smooth. Since the fibre of $h$ over a reduced cycle $\sum_{1\leq i\leq n}x_i$ is the smooth variety $\prod_{1\leq i\leq n} \BP(F|_{x_i}) \cong (\BP^{r-1})^n$, we obtain the result.   
\end{proof}

\begin{prop}
\label{prop:smooth-open-quot}
Let $Y$ be a smooth $d$-dimensional variety, $F \in \Coh Y$ a locally free sheaf of rank $r\geq 1$, and let $n$ be a positive integer. Then $\Quot^\circ_Y(F,n)$ is a smooth variety of dimension $n(r-1+d)$.
\end{prop}

\begin{proof}
The morphism $h$ is open; its target is irreducible of dimension $nd$, and its fibres $(\BP^{r-1})^n$ are irreducible of the same dimension, so $\Quot^\circ_Y(F,n)$ is irreducible of dimension $n(r-1+d)$ by \Cref{lemma:conn-irr}\ref{IRR}. Its smoothness follows from the smoothness of $\Sym^n(Y) \setminus \Delta^{\mathsf{big}}$ and \Cref{lemma:h-smooth}.
\end{proof}

For $F$ locally free, in light of \Cref{prop:smooth-open-quot}, we refer to $\Quot_Y^{\main}(F,n)$ as the \emph{smoothable component} of the Quot scheme.

We close this section with a result on the Quot scheme of 2 points attached to a locally free sheaf of rank 3 on a smooth 3-fold. This will be used in \Cref{prop:int+CM+8}, a crucial step in the proof of \Cref{thm:MAIN-small-n}\ref{M-SMALL2}.

\begin{prop}
\label{prop:quot-red-CM}
Let $Y$ be a smooth quasiprojective 3-fold, $F$ a locally free sheaf of rank 3 over $Y$. Then $\Quot_Y(F,2)$ is reduced and Cohen--Macaulay.
\end{prop}

\begin{proof}
By the \'etale-local structure of Quot schemes \cite[Prop.~A.3]{Beentjes-Ricolfi}, and because the properties to be checked are of \'etale-local nature, it is enough to treat the case $Y=\BA^3$ and $F=\OO^{\oplus 3}$. In this case, by the ADHM construction \cite{Nakajima, Beentjes-Ricolfi,Joachim-Sivic,Cazzaniga-Ralaivaosaona-Ricolfi,Ricolfi-ncQuot-motivic}, there is a $\GL_2$-principal bundle 
\begin{equation}
\label{eqn:prin-bdl-representation}
\begin{tikzcd}
    U \arrow{r} & \Quot_{\BA^3}(\OO^{\oplus 3},2),
\end{tikzcd}
\end{equation}
where $U \subset \Mat_2(\BC)^{\oplus 3} \times (\BC^2)^{\oplus 3} \cong \BA^{18}$ is the locally closed subset consisting of pairs $(\boldit{X},\boldit{v})$ where $\boldit{X}$ is a triple of matrices $(X_1,X_2,X_3) \in \Mat_2(\BC)^{\oplus 3}$ and $\boldit{v}$ is a triple of vectors $(v_1,v_2,v_3) \in (\BC^2)^{\oplus 3}$ such that the matrices pairwise commute and the $\BC$-linear span of all monomials  $X_1^{a_1}X_2^{a_2}X_3^{a_3}$ applied to the vectors $v_1$, $v_2$ and $v_3$ is the whole $\BC^2$.

Now, by \cite{Joachim-Sivic}, the Quot scheme $\Quot_{\BA^3}(\OO^{\oplus 3},2)$ is irreducible of dimension 10, and by \cite[Eqn.~(6.2)]{Pielasa} its singular locus is the 5-dimensional smooth variety
\[
\Sing \Quot_{\BA^3}(\OO^{\oplus 3},2) \cong \BA^3 \times \Gr(2,3)\cong \BA^3 \times \BP^2.
\]
Now fix a singular point $q = [\OO^{\oplus 3} \onto T] \in \Quot_{\BA^3}(\OO^{\oplus 3},2)$. By \Cref{prop:smooth-open-quot}, $T$ is supported at a single point. Via the principal bundle \eqref{eqn:prin-bdl-representation}, $q$ can be represented as a tuple $[X_1,X_2,X_3,v_1,v_2,v_3]$, where $[-]$ means ``up to $\GL_2$''. By \cite[Thm.~6.6]{Pielasa}, $X_1$, $X_2$, and $X_3$ are scalar matrices, i.e.~$X_i = \lambda_i\cdot \Id_2$ for some $\lambda_i \in \BC$ --- note that the triple $(\lambda_1,\lambda_2,\lambda_3) \in \BA^3$ is exactly the support of the quotient sheaf corresponding to $q$. On the other hand, the spanning triple $[v_1,v_2,v_3]$ corresponds to a surjection $(\varphi\colon \BC^3 \onto \BC^2) \in \Gr(2,3)$. Now fix a Zariski open chart $W\subset \Gr(\BC^3,2)$ such that, under the stated correspondence, the vectors $v_1,v_2$ span the module $T$, so that $W \cong \Spec \BC[v_{31},v_{32}]$, with $v_{31}$ and $v_{32}$ the coordinates of the free vector $v_3$. In this chart, $q$ takes the form $(\lambda_1,\lambda_2,\lambda_3;v_{31},v_{32})$. Now, in a formal neighbourhood of $q$, a typical point is of the form 
\[
\widehat q = (\widehat X_1,\widehat X_2,\widehat X_3;\widehat v_{31},\widehat v_{32}),
\]
with $\widehat X_i = (\lambda_i+\tau_i)\cdot \Id_2 + A_i$, for a unique traceless matrix $A_i \in \mathfrak{sl_2}$ and infinitesimal deformation parameters $\tau_i$. Now, the condition $[\widehat X_i,\widehat X_j]=0$ for $1\leq i<j\leq 3$ means
\[
0 = [(\lambda_i+\tau_i)\cdot \Id_2 + A_i,(\lambda_j+\tau_j)\cdot \Id_2 + A_j] = [A_i,A_j],
\]
therefore the Quot scheme is, \'etale locally around $q$, of the form
\[
\BA^3 \times \BA^2 \times \mathsf C_{3,2},
\]
where $\BA^3$ has coordinates $\lambda_1,\lambda_2,\lambda_3$, the second factor has coordinates $v_{31}$, $v_{32}$, and 
\[
\mathsf C_{3,2} = \Set{(A_1,A_2,A_3)\,\in\,\mathfrak{sl}_2^3 | [A_i,A_j] = 0 ~\mathrm{for}~1\leq i<j\leq 3} \subset \mathfrak{sl}_2^3
\]
is the commuting scheme of triples of matrices in $\mathfrak{sl}_2$. We are then reduced to proving that $\mathsf C_{3,2}$ is reduced and Cohen--Macaulay. This will be achieved in the next lemma.
\end{proof}

\begin{lemma}
The commuting scheme
\[
\mathsf C_{3,2} = \Set{(A_1,A_2,A_3)\,\in\,\mathfrak{sl}_2^3 | [A_i,A_j] = 0 ~\mathrm{for}~1\leq i<j\leq 3} \subset \mathfrak{sl}_2^3
\]
is reduced and Cohen--Macaulay.
\end{lemma}

\begin{proof}
Fix the usual basis
\[
e 
= \begin{pmatrix}
0 & 1 \\
0 & 0
\end{pmatrix}, \qquad
f 
= \begin{pmatrix}
0 & 0 \\
1 & 0
\end{pmatrix}, \qquad
h 
= \begin{pmatrix}
1 & 0 \\
0 & -1
\end{pmatrix}
\]
of $\mathfrak{sl}_2$. For a triple $A_\bullet \defeq (A_1,A_2,A_3)\,\in\,\mathfrak{sl}_2^3$, write $A_i = \alpha_i e + \beta_i f + \gamma_i h$, i.e.
\[
A_i = 
\begin{pmatrix}
    \gamma_i & \alpha_i \\
    \beta_i & -\gamma_i
\end{pmatrix}.
\]
Then the condition $[A_i,A_j]=0$ is equivalent to the 3 equations
\begin{equation}
\label{eqn:commutativity}
\begin{split}
\alpha_i \beta_j-\alpha_j \beta_i &= 0, \\
\gamma_i\alpha_j - \gamma_j \alpha_i &= 0, \\
\gamma_j\beta_i-\gamma_i\beta_j &=0.
\end{split}
\end{equation}
Consider the $3 \times 3$ matrix
\[
M_{A_\bullet} = 
\begin{pmatrix}
    \alpha_1 & \beta_1 & \gamma_1 \\
    \alpha_2 & \beta_2 & \gamma_2 \\
    \alpha_3 & \beta_3 & \gamma_3 \\
\end{pmatrix}.
\]
For each pair $1\leq i<j\leq 3$, the three equations in \eqref{eqn:commutativity} are precisely, up to sign, the $2\times 2$ minors obtained from the two rows $i,j$. As
$(i,j)$ varies, these are exactly all $2\times 2$ minors of
$M_{A_\bullet}$. In other words, we have the description 
\[
\mathsf C_{3,2} \cong \Set{M \in \Mat_{3}(\BC) | \rk M \leq 1}.
\]
So $\mathsf C_{3,2}$ is a determinantal variety. Since $\dim \mathsf C_{3,2} = 5$, which confirms that $\mathsf C_{3,2}$ has the expected codimension inside $\Mat_{3}(\BC)\cong \BA^9$, namely $(3-1)(3-1) = 4$, we can conclude that $\mathsf C_{3,2}$ is Cohen--Macaulay by \cite[Lemma A.7.2]{Ful} (see also \cite[Thm.~14.4(c)]{Ful}). Since $\mathsf C_{3,2}$ is generically reduced by \cite[Thm.~B]{Joachim-Sivic}, and we work over $\BC$, the scheme $\mathsf C_{3,2}$ is generically smooth (i.e.~it has property $R_0$), so being Cohen--Macaulay it must be reduced by \cite[\href{https://stacks.math.columbia.edu/tag/033P}{Tag 033P}]{stacks-project}.
\end{proof}


\section{\texorpdfstring{Moduli of stable rank 2 sheaves with $c_1=-1$ and $c_2=1$}{}} \label{sec:sheaves}

From now on, we fix $Y = \BP^3$ along with the polarisation $H$ corresponding to the hyperplane class. We are interested in the moduli space
\[
M_{c_3}\defeq M_{\BP^3}^{\st}(2,-1,1,c_3),
\]
whose $\BC$-points correspond to stable torsion-free sheaves of rank $2$ with $c_1=-1$, $c_2=1$ and third Chern class $c_3$.

\subsection{The reflexive hull short exact sequence}
\label{sec:Chern_classes}
If $E$ is a coherent sheaf of rank 2 on $\BP^3$, Hirzebruch--Riemann--Roch yields the identity \cite[Thm.~2.3]{Hart_StableReflexive}
\begin{equation}
\label{eqn:Chern_character}
\chi(E(m)) = 2 + \binom{c_1(E(m))+3}{3} - 2c_2(E(m)) + \frac{1}{2}(c_3-c_1(E(m))c_2(E(m)))-1.
\end{equation}

\begin{remark}
\label{same-stability}
Let $E$ be $\mu$-semistable of rank 2 and Chern classes $c_i$ for $i=1,2,3$. By \Cref{eqn:Chern_character}, one has $c_3-c_1c_2\,\in\,2\BZ$. So if $c_1=-1$ and $c_2=1$, one has that $c_3$ must be odd; since $\rk E=2$ and $c_1\cdot H^2=-1$ are coprime, we have that $E$ is in fact $\mu$-stable. Thus, by \eqref{stab-notions}, all notions of stability agree and $M_{c_3}$ is a projective scheme.     
\end{remark}

\begin{lemma}
\label{lemma:M_is_fine}
The moduli space $M_{c_3}$, if nonempty, is a fine moduli space.
\end{lemma}

\begin{proof}
This is a simple application of \cite[Cor.~4.6.6]{modulisheaves}.
\end{proof}

\begin{lemma} 
\label{lemma:length-quotient}
If $E$ is a semistable rank 2 sheaf with Chern classes $(c_1,c_2)=(-1,1)$, then it sits in an exact sequence
\begin{equation} 
\label{eq:e1}
\begin{tikzcd}
0 \arrow{r} & E \arrow{r} & R \arrow{r} & Q \arrow{r} & 0,
\end{tikzcd}
\end{equation}
where $R$ is a stable rank 2 reflexive sheaf with Chern classes $(-1,1,1)$ and the quotient $Q$ is a 0-dimensional sheaf of length 
\[
n\defeq \frac{1-c_3(E)}{2}.
\]
In particular, $c_3(E)\leq 1$, and $c_3(E)=1$ if and only if $E$ is reflexive. 
\end{lemma}

\begin{proof}
The exact sequence is given by setting $R \defeq E^{\ast\ast}$ and letting $E \into R$ be the canonical map $\theta_E$, cf.~\Cref{sec:Background}. In principle, $\dim Q \leq 1$. The fact that $\dim Q = 0$ is proved in \cite[Prop.~3.1]{Gholampour2017a}. Since $\frac{1-c_3(E)}{2} = -\chi(E)$ by \Cref{eqn:Chern_character}, all we need to do is show $\chi(R)=0$. Since $R/E$ is $0$-dimensional, we have $c_i(R)=c_i$ for $i\leq 2$. But for a $\mu$-stable reflexive rank $2$ sheaf $R$ on $\BP^3$ such that $c_1(R)=-1$, by Hartshorne's result \cite[Thm.~8.2(d)]{Hart_StableReflexive} one has $0 \leq c_3(R) \leq c_2(R)^2$, which yields $c_3(R)=1$. Therefore applying \Cref{eqn:Chern_character} again we find $\chi(R) = 0$, as required. 
\end{proof}

Let us denote by 
\begin{equation}
\label{eqn:reflexive-inclusion}
N \defeq N_{\BP^3}^{\mathsf{ref}}(2,-1,1,1) \subset M_1
\end{equation}
the moduli space of reflexive sheaves, cf.~\Cref{def:M-reflexive}. By (the last sentence of) \Cref{lemma:length-quotient}, all points of $M_1$ correspond to reflexive sheaves, so the inclusion \eqref{eqn:reflexive-inclusion} is an equality, i.e.~$N = M_1$.

\begin{remark}
\label{rmk:numerology}
The Chern character of the sheaves parametrised by $M_{c_3}$ is
\begin{equation}
\label{eqn:chern-char}
v = \left(2,-1,-\frac{1}{2},\frac{2+3c_3}{6}\right) = \left(2,-1,-\frac{1}{2},\frac{5}{6}-n\right),
\end{equation}
with $n = \frac{1-c_3}{2}$ as in \Cref{lemma:length-quotient}. This datum is equivalent to the Chern class $(-1,1,1-2n)$.

As seen during the proof of \Cref{lemma:length-quotient}, by \cite[Thm.~8.2(d)]{Hart_StableReflexive}, the reflexive hulls $E^{\ast\ast}$ of the sheaves $E$ parametrised by $M_{c_3}$ have Chern class $(-1,1,1)$, and thus define points of $M_1$. Since by \cite[Prop.~2.6]{Hart_StableReflexive} for a reflexive rank 2 sheaf $R$ one has that $c_3(R) = h^0(\BP^3,\lExt^1(R,\OO_{\BP^3}))$, and moreover the support of $\lExt^1(R,\OO_{\BP^3})$ coincides with $\Sing R$, there is a unique point $p \in \BP^3$ such that $R$ is not free at $p$.
\end{remark}

We now give a proof of the following structural result, corresponding to \Cref{MAIN-THM-a}\ref{MTHM1}.

\begin{prop}
\label{prop:reflexive=P^3}
There is an isomorphism
\[
\begin{tikzcd}
N \arrow{r}{\sim} & \BP^3.
\end{tikzcd}
\]
\end{prop}

\begin{proof}
Let $\mathscr R \in \Coh(\BP^3 \times N)$ be the universal reflexive sheaf. Then $\lExt^1_{\BP^3 \times N}(\mathscr R,\OO_{\BP^3 \times N})$ is $N$-flat (by an application of \cite[Thm.~1.9(ii)]{PicScheme}) and satisfies base change (by an application of \cite[Thm.~1.9(i)]{PicScheme}), i.e.~for any point $p = [R] \in N$ the canonical morphism
\[
\begin{tikzcd}
\lExt^1_{\BP^3 \times N}(\mathscr R,\OO_{\BP^3 \times N}) \otimes_{\OO_{N}} \kappa(p) \arrow{r} & \lExt^1_{\BP^3}(R,\OO_{\BP^3})
\end{tikzcd}
\]
is an isomorphism. Then, since $\lExt^1_{\BP^3}(R,\OO_{\BP^3})$ is just a skyscraper sheaf on $\BP^3$, in particular a stable sheaf, the $N$-flat family $\lExt^1_{\BP^3 \times N}(\mathscr R,\OO_{\BP^3 \times N})$ defines a modular morphism
\[
\begin{tikzcd}
\sigma\colon N \arrow{r} & M_{\BP^3}^{\st}(0,0,0,1)
\end{tikzcd}
\]
to the moduli space of stable torsion sheaves with Chern character $(0,0,0,1)$. This moduli space is isomorphic to $\BP^3$ (by \cite[Ex.~4.3.6]{modulisheaves} recalled in \Cref{rmk:moduli_of_0-cycles}), and the map $\sigma$, sending a reflexive sheaf to its unique singular point, is bijective on points by Hartshorne's results \cite[Sec.~9]{Hart_StableReflexive}. In \cite[Thm.~9.2]{Hart_StableReflexive} it is proven that $N$ is an irreducible, smooth rational projective $3$-fold. Therefore, since both source and target of $\sigma$ are smooth of dimension $3$, the morphism $\sigma$ is an isomorphism by Zariski's Main Theorem.
\end{proof}

\subsection{Resolution for the reflexive hulls}
\label{sec:resolution}
The case $c_3=1$ is well understood: by \Cref{prop:reflexive=P^3}, every sheaf parametrised by $M_1$ is reflexive, and by \cite[Sec.~9]{Hart_StableReflexive} (see also \cite[Sec.~6.1.3]{Jardim-Mu}), there is a resolution
\begin{equation}
\label{eqn:resolution-F}
\begin{tikzcd}
0 \arrow{r} & \OO_{\BP^3}(-2) \arrow{r}{\boldit{s}} & \OO_{\BP^3}(-1)^{\oplus3} \arrow{r} & R \arrow{r} & 0,
\end{tikzcd}
\end{equation}
which also gives a geometric meaning to the unique point $p = \Sing R$ where $R$ is not locally free.
In a little more detail, as explained in the proof of \cite[Lemma 9.4]{Hart_StableReflexive}, the sheaf $R(1)$ is generated by global sections, $h^0(\BP^3,R(1))=3$ and $\HH^i(\BP^3,R(1))=0$ for $i>0$. Choose a basis of $\HH^0(\BP^3,R(1))$, and let $\OO_{\BP^3}^{\oplus 3} \onto R(1)$ be the induced surjection. Since $R(1)$ has homological dimension 1, the kernel $\SK$ is a line bundle, and since $c_1(\SK)=-c_1(R(1))=c_1(R)=-1$, we must have $\SK = \OO_{\BP^3}(-1)$. Twisting by $\OO_{\BP^3}(-1)$ gives the resolution \eqref{eqn:resolution-F}. The morphism $\boldit{s} = (s_1,s_2,s_3)^{\mathsf t}$ is given by three linearly independent\footnote{If they were linearly dependent, the cokernel of $\boldit{s}$ would split as a direct sum, contradicting the stability of $R$.} linear forms $s_i\in \HH^0(\BP^3,\op3(1))$, and dualising \eqref{eqn:resolution-F}, we obtain 
\begin{equation} \label{eqn:dual-resolution-F} \begin{tikzcd}
0 \arrow{r} & R^\ast \arrow{r} & \OO_{\BP^3}(1)^{\oplus 3}
\arrow{r}{\boldit{s}^{\rm t}} & \OO_{\BP^3}(2) \arrow{r} & \lExt^1(R,\OO_{\BP^3}) \arrow{r} & 0,
\end{tikzcd}\end{equation}
so that $\lExt^1(R,\OO_{\BP^3})=\OO_p$, where $p$ is the common zero locus of $s_1,s_2,s_3$, namely the locus where $\boldit{s}$ drops rank; in other words, $p = \Sing R$.

\begin{remark}
    As explained in \cite[Sec.~9]{Hart_StableReflexive}, every sheaf $R$ fitting in a sequence of the form of \Cref{eqn:resolution-F} belongs to $M_1$.
\end{remark}

\smallbreak
Choosing two of the sections, say $(s_1,s_2)$, we obtain the commutative diagram
\begin{equation} \label{diag-res}
\begin{tikzcd}
&& 0\arrow{d} & 0\arrow{d} & \\
&& \OO_{\BP^3}(-1)\arrow{d}\arrow[equal]{r} & \OO_{\BP^3}(-1)\arrow{d}{\sigma} & \\
0\arrow{r} & \OO_{\BP^3}(-2)\arrow{r}{\boldit{s}}\arrow[equal]{d} & \OO_{\BP^3}(-1)^{\oplus 3}\arrow{d}\arrow{r} & R\arrow{d}\arrow{r} & 0 \\
0\arrow{r} & \OO_{\BP^3}(-2)\arrow{r}{\left(\begin{smallmatrix}
    s_1 \\ s_2
\end{smallmatrix}\right)} & \OO_{\BP^3}(-1)^{\oplus 2}\arrow{d}\arrow{r} & \coker(\sigma)\arrow{r}\arrow{d} & 0 \\
&& 0 & 0 & 
\end{tikzcd}
\end{equation}
where $\sigma\in \HH^0(\BP^3,R(1))$ is an induced section. From the exact bottom row, we conclude that $\coker(\sigma) \cong \SI_\ell$, where $\ell\subset\BP^3$ is the line defined by $s_1=s_2=0$. In particular, $\ell$ contains the singular point $p = \Sing R$.

\smallbreak
The resolution in \eqref{eqn:resolution-F} will play a central role throughout this paper.

\subsection{Construction of stable rank 2 sheaves}
\label{sec:construction}
One can easily construct examples of $\mu$-stable rank 2 sheaves on $\BP^3$ with Chern classes $(c_1,c_2,c_3)=(-1,1,1-2n)$ for every $n\geq 1$. Take a $\mu$-stable reflexive sheaf $R$ with Chern classes $(c_1,c_2,c_3)=(-1,1,1)$ as in \Cref{sec:resolution}, and let $Z\subset\BP^3$ be a reduced, 0-dimensional scheme of length $n$ not containing the singular point of $R$, i.e.
\begin{equation}
\label{eqn:no-intersection}
\Supp \lExt^1(R,\OO_{\BP^3}) \cap Z = \emptyset.
\end{equation}
Let $q_1,\ldots,q_n \in \BP^3$ be the points in the support of $Z$, which we assume to be distinct, so that  $\OO_Z=\bigoplus_{i=1}^n\OO_{q_i}$ with $q_i\ne q_j$ for all $i \neq j$. We have
\[
\Hom(R,\OO_Z) = \bigoplus_{i=1}^n\Hom(R,\OO_{q_i}) = \bigoplus_{i=1}^n \kappa(q_i)^2 \cong \BC^{2n},
\]
since $R$ is locally free around $Z$ by \eqref{eqn:no-intersection}. Choosing nontrivial elements at each summand yields an epimorphism $\varphi\colon R\onto\OO_Z$ and thus a short exact sequence
\begin{equation} 
\label{eq:e0}
\begin{tikzcd}
0 \arrow{r} & E \arrow{r} & R \arrow{r}{\varphi} & \OO_Z \arrow{r} & 0.
\end{tikzcd}
\end{equation}
Dualising yields an isomorphism $R^\ast \simto E^\ast$, and by \cite[Prop.~1.10]{Hart_StableReflexive} we have $R^\ast \cong R(1)$, so that in particular $E$ is $\mu$-stable. In addition, $\ch(E) = \ch(R)-\ch(\OO_Z) = (2,-1,-\frac{1}{2},\frac{5}{6}-n)$, so by \Cref{rmk:numerology} the Chern classes of $E$ are precisely $(-1,1,1-2n)$, as desired.

\subsection{Construction of quotients}
\label{sec:construct-quotients}

We saw in the previous section that it is straightforward to construct points of $M_{c_3}$ out of points of $M_1$. In this section, we consider instead the space of homomorphisms $\Hom(R,Q)$ for $[R]$ a point in $N=M_1$ and $Q$ a 0-dimensional sheaf. Even though one always has $2n \leq \hom(R,Q) \leq 3n$ for $n \geq 1$, the existence of a \emph{surjection} $R \onto Q$ is not granted in general. We can nevertheless prove the following result.

\begin{lemma}
\label{lemma:homR-O}
Let $R$ be a sheaf corresponding to a point of $N=M_1$, and let $p = \Sing R$. Let $Q$ be a 0-dimensional sheaf of length $n$. Then
\begin{itemize}
    \item [\mylabel{bound}{\normalfont{(i)}}]$\dimhom(R,Q) \leq 3n$,
    \item [\mylabel{hom(R,Q)}{\normalfont{(ii)}}] If $p\notin\Supp Q$, then $\dimhom(R,Q) = 2n$,
    \item [\mylabel{red-pt}{\normalfont{(iii)}}] $\dimhom(R,\OO_p) = 3$.
    \item [\mylabel{surj}{\normalfont{(iv)}}] If $Q = \OO_Z$ where $Z\subset \BP^3$ is a reduced 0-dimensional scheme, then the locus of surjective maps $\Hom^{\surj}(R,\OO_Z)\subset \Hom(R,\OO_Z)$ is a nonempty open subset of $\Hom(R,\OO_Z)$.  
\end{itemize}
\end{lemma}

\begin{proof}
Applying $\Hom(-,Q)$ to \eqref{eqn:resolution-F}, we get an inclusion $\Hom(R,Q) \subset \Hom(\op3(-1)^{\oplus3},Q) = \BC^{3n}$, which proves \ref{bound}. To prove \ref{hom(R,Q)}, notice that $\Hom(R,Q) \cong\HH^0(\BP^3,Q)^{\oplus 2}\cong\BC^{2n}$, since $p\notin\Supp Q$. As for \ref{red-pt}, the map
\[
\begin{tikzcd}
\Hom(\op3(-1)^{\oplus 3},\OO_p) \arrow{r}{-\,\circ \boldit{s}} &
\Hom(\op3(-2),\OO_p)
\end{tikzcd}
\]
obtained applying $\Hom(-,\OO_p)$ to \eqref{eqn:resolution-F} vanishes, confirming $\dimhom(R,\OO_p) = \dimhom(\op3(-1)^{\oplus 3},\OO_p) = 3$. Finally, to prove \ref{surj}, let $Z \subset \BP^3$ be a reduced 0-dimensional subscheme. Then 
\begin{equation}
\label{eqn:hom(R,OZ)}
\Hom(R,\OO_Z) = \bigoplus_{q \in Z} \Hom(R,\OO_q).
\end{equation}
The surjective locus $\Hom^{\surj}(R,\OO_Z)$ is Zariski open and nonempty, since any tuple of nonzero homomorphisms $R \to \OO_q$ assembles to a global surjective homomorphism $R \onto \OO_Z$.
\end{proof}

\begin{remark}
\label{rmk:twohoms}
In fact, thanks to Items \ref{hom(R,Q)} and \ref{red-pt} of \Cref{lemma:homR-O}, in Equation \eqref{eqn:hom(R,OZ)} one can say more, namely
\[
\Hom(R,\OO_Z) = \bigoplus_{q \in Z} \Hom(R,\OO_q) = 
\begin{cases}
    \BC^{2n} & \mbox{ if }Z\cap \Sing R = \emptyset \\ \\
    \BC^{2n+1} & \mbox{ if }Z\cap \Sing R \neq \emptyset.
\end{cases}
\]
\end{remark}

\section{Cohomological calculations}
\label{sec:Cohomological}
We collect in this section a series of calculations of Ext groups and Ext sheaves attached to sheaves $E$ parametrised by the moduli space $M_{c_3}$. Recall from \Cref{lemma:length-quotient} that, for a point $[E] \in M_{c_3}$, the length $n = \chi(Q)$ of $Q = E^{**}/E$ and the third Chern class $c_3 = c_3(E)$ are related by the formula $c_3 = 1-2n$.

The next result is completely general, and applies to an arbitrary point $[E]$ of $M_{c_3}$.

\begin{lemma} 
\label{lem:cohomology}
Let $E$ be a rank 2 stable sheaf with Chern classes $(c_1,c_2,c_3)=(-1,1,1-2n)$ for some $n \geq 1$. Set $R=E^{\ast\ast}$ and $Q=R/E$. Then
\begin{enumerate}
\item [\mylabel{coh-1}{\normalfont{(i)}}] $h^0(\BP^3,E(1))\leq 2$.
\item [\mylabel{coh-2}{\normalfont{(ii)}}] $\lExt^1(E,\op3)\cong\OO_p$ where $p=\Sing R$.
\item [\mylabel{coh-3}{\normalfont{(iii)}}] $\lExt^2(E,\op3) \cong \lExt^3(Q,\op3)$. In particular, $\lExt^2(E,\op3)$ is a 0-dimensional sheaf of length $n$, with $\Supp\lExt^2(E,\op3)=\Supp Q$.
\item [\mylabel{coh-4}{\normalfont{(iv)}}] E satisfies $\hd(E)=2$ and
\begin{equation}
\label{eqn:Sing(E)}
\Sing E=\Sing R\cup\Supp Q.   \end{equation}
\item [\mylabel{coh-5}{\normalfont{(v)}}] For every $i>0$, the sheaf $\lExt^i(E,E)$ is 0-dimensional.
\end{enumerate}
\end{lemma}
\begin{proof}
Let us start proving \ref{coh-1}. If $h^0(\BP^3,E(1))>0$, then we obtain, putting together the sequences in \eqref{diag-res} and \eqref{eq:e1}, the following commutative diagram
\begin{equation} 
\label{eq:diag1}
\begin{tikzcd}
& 0\arrow{d} & 0\arrow{d} & & \\
& \op3(-1)\arrow[equal]{r}\arrow{d} & \op3(-1)\arrow{d} & & \\
0\arrow{r} & E\arrow{r}\arrow{d} & R\arrow{d}\arrow{r} & Q\arrow[equal]{d}\arrow{r} & 0 \\
0\arrow{r} & \SI_L\arrow{r}\arrow{d} & \SI_\ell\arrow{r}\arrow{d} & Q\arrow{r} & 0 \\
& 0 & 0 & &
\end{tikzcd}
\end{equation}
where $\ell$ is a line in $\BP^3$. The bottom sequence implies that $L$, a ``line plus $n$ points'', defines a point of the Hilbert scheme $\Hilb^{t+1+n}(\BP^3)$. It follows that $h^0(\BP^3,E(1))=1$ when $L$ is not (schematically) planar, and $h^0(\BP^3,E(1))=2$ when $L$ is planar. This proves \ref{coh-1}.

As for Items \ref{coh-2} and \ref{coh-3}, simply dualise the sequence 
\[
\begin{tikzcd}
0 \arrow{r} & E \arrow{r} & R \arrow{r} & Q \arrow{r} & 0
\end{tikzcd}
\]
to obtain
\[
\begin{tikzcd}
\OO_p = \lExt^1(R,\op3) \arrow{r}{\sim} & \lExt^1(E,\op3) 
\end{tikzcd}
\]
thanks to the vanishings $\lExt^{1}(Q,\OO_{\BP^3})=\lExt^{2}(Q,\OO_{\BP^3})=0$, and
\[
\begin{tikzcd}
\lExt^2(E,\op3) \arrow{r}{\sim} &\lExt^3(Q,\op3)    
\end{tikzcd}
\]
thanks to the vanishing $\lExt^{>1}(R,\OO_{\BP^3})=0$ due to $R$ being of homological dimension 1. 

To prove \ref{coh-4}, note that the exact piece 
\[
\begin{tikzcd}
0 = \lExt^3(R,\OO_{\BP^3}) \arrow{r} & \lExt^3(E,\OO_{\BP^3}) \arrow{r} & \lExt^4(Q,\op3) = 0
\end{tikzcd}
\]
forces $\lExt^3(E,\OO_{\BP^3}) = 0$, which by Item \ref{coh-3} confirms $\hd(E)=2$. The identity \eqref{eqn:Sing(E)} follows from the general formula $\Sing E=\Supp\lExt^1(E,\op3)\cup\Supp\lExt^2(E,\op3)$ for a 3-fold, see e.g.~\cite[p.~73]{OSS}, combined with Items \ref{coh-2} and \ref{coh-3}.

Now Item \ref{coh-5} follows from the fact that if $E$ is locally free at $x \in \BP^3$ then $\lExt^i(E,E)_x = \Ext^i(E_x,E_x) = 0$ for all $i>0$, thus $\Supp \lExt^i(E,E) \subset \Sing E$, which is 0-dimensional by Item \ref{coh-4}.
\end{proof}

\begin{lemma}
\label{lem:exts}
Let $E$ be a rank 2 semistable sheaf with Chern classes $(c_1,c_2,c_3)=(-1,1,1-2n)$ for some $n \geq 1$. Set $R=E^{\ast\ast}$ and $Q=R/E$. Then
\begin{enumerate}
\item [\mylabel{ext-i}{\normalfont{(i)}}] $\Ext^1(E,E)\cong  \HH^1(\BP^3,\lHom(E,E)) \oplus \HH^0(\BP^3,\lExt^1(E,E))$,
\item [\mylabel{ext-ii}{\normalfont{(ii)}}] $\Ext^2(E,E)\cong  \HH^0(\BP^3,\lExt^2(E,E)) \cong \Hom(E,Q)^\ast$,
\item [\mylabel{ext-iii}{\normalfont{(iii)}}] $\Supp \lExt^2(E,E)=\Supp Q$.
\end{enumerate}
\end{lemma}

\begin{proof}
By \Cref{lem:cohomology}\ref{coh-5} we know that $\dim\lExt^{>0}(E,E)=0$, so $\HH^p(\lExt^q(E,E))=0$ for $p>0$ and $q>0$ and the only potentially nontrivial differentials in the local-to-global spectral sequence $\mathrm{E}^{pq}_2=\HH^p(\BP^3,\lExt^q(E,E))\Rightarrow \Ext^{p+q}(E,E)$ are 
\[
\begin{tikzcd}[row sep=tiny]
\mathrm d_2^{01} \colon \HH^0(\BP^3,\lExt^1(E,E)) \arrow{r} & \HH^2(\BP^3,\lHom(E,E)), \\
\mathrm d_3^{02} \colon \HH^0(\BP^3,\lExt^2(E,E)) \arrow{r} & \HH^3(\BP^3,\lHom(E,E)).
\end{tikzcd}
\]
They are the connecting maps in the associated low-degree terms exact sequence 
\[
\begin{tikzcd}
0 \arrow[r]
    & \HH^1(\BP^3,\lHom(E,E)) \arrow[r]
    & \Ext^1(E,E) \arrow[r]
        \arrow[d, phantom, ""{coordinate, name=Z}]
    & \HH^0(\BP^3,\lExt^1(E,E)) \arrow[dll,rounded corners,
to path={ -- ([xshift=2ex]\tikztostart.east)
|- (Z) [near end]\tikztonodes
-| ([xshift=-2ex]\tikztotarget.west) -- (\tikztotarget)}] \\
& \HH^2(\BP^3,\lHom(E,E)) \arrow[r]\arrow[d, phantom, ""{coordinate, name=U}]
    & \Ext^2(E,E) \arrow[r]
& \HH^0(\BP^3,\lExt^2(E,E)) 
\arrow[dll,rounded corners,
to path={ -- ([xshift=2ex]\tikztostart.east)
|- (U) [near end]\tikztonodes
-| ([xshift=-2ex]\tikztotarget.west) -- (\tikztotarget)}] \\
&\HH^3(\BP^3,\lHom(E,E))\arrow[r] & 0. &
\end{tikzcd}
\]
We used $\Ext^3(E,E)=0$, which follows from $\BP^3$ being Fano and $E$ being stable. Both Item \ref{ext-i} and the first isomorphism of Item \ref{ext-ii} would follow from the vanishings $\HH^{\geq 2}(\BP^3,\lHom(E,E))=0$, which we now prove. Applying the functor $\lHom(-,E)$ to the reflexive hull sequence 
\begin{equation}
\label{eq:e113}
\begin{tikzcd}
0 \arrow{r} & E \arrow{r} & R \arrow{r} & Q \arrow{r} & 0,
\end{tikzcd}
\end{equation}
we get, noting that $\lHom(Q,E)=0$,
\[
\begin{tikzcd}
0 \arrow{r} & \lHom(R,E) \arrow{r}{b} & \lHom(E,E) \arrow{r} & \lExt^1(Q,E) \arrow{r} & \cdots
\end{tikzcd}
\]
But $\dim\lExt^1(Q,E)=0$, so the cokernel of the map $b$ is also 0-dimensional, and we conclude that $\HH^j(\BP^3,\lHom(E,E))\cong \HH^j(\BP^3,\lHom(R,E))$ for $j\geq 2$. Now, the canonical morphism
$R^\ast\otimes E \to \lHom(R,E)$ is an isomorphism away from the point $\Sing R$, so it has 0-dimensional kernel and cokernel. Factoring it through its image and taking cohomology then yields isomorphisms 
\[
\begin{tikzcd}
\HH^j(\BP^3,R^\ast\otimes E) \arrow{r}{\sim} & \HH^j(\BP^3,\lHom(R,E)), \qquad j\geq 2.
\end{tikzcd}
\]
We need to show the source vanishes. But $R^\ast=R(1)$, and tensoring \eqref{eqn:resolution-F} with $E(1)$ yields\footnote{Note that $\LTor_1(R,E)=0$ by \cite[Lemma A.1(ii)]{CJM22}.}
\begin{equation}
\label{eqn:res-again}
\begin{tikzcd}
0 \arrow{r} & E(-1) \arrow{r} & E^{\oplus3} \arrow{r} & R(1)\otimes E \arrow{r} & 0,
\end{tikzcd}
\end{equation}
so that, passing to cohomology, we get an exact piece
\[
\begin{tikzcd}
\HH^j(\BP^3,E)^{\oplus3} \arrow{r} & \HH^j(\BP^3,R(1)\otimes E) \arrow{r} & \HH^{j+1}(\BP^3,E(-1)). \end{tikzcd}
\]
Note that $\HH^j(\BP^3,E(k))=\HH^j(\BP^3,R(k))$ for $j\geq 2$ and all $k$ (since $Q=R/E$ is 0-dimensional), and $\HH^2(\BP^3,R)=\HH^3(\BP^3,R(-1))=0$ by \eqref{eqn:resolution-F}. Thus $ \HH^j(\BP^3,R(1)\otimes E)=0$ for $j\geq 2$, and it follows that $\HH^j(\BP^3,\lHom(E,E))=0$ in the same range, as desired.

To prove the second isomorphism of Item \ref{ext-ii}, apply again the functor $\lHom(-,E)$ to the exact sequence \eqref{eq:e113} to get
\[
\begin{tikzcd}
\lExt^2(R,E)=0 \arrow{r} & \lExt^2(E,E) \arrow{r}{\sim} & \lExt^3(Q,E) \arrow{r} & 0=\lExt^3(R,E), 
\end{tikzcd}
\]
with the vanishings because $R$ has homological dimension 1. Moreover, since $\dim\lExt^i(Q,E) \leq 0$ for every $i \geq 0$, the local-to-global spectral sequence implies $\Ext^3(Q,E)\cong \HH^0(\BP^3,\lExt^3(Q,E))$, while Serre duality together with 0-dimensionality of $Q$ yield $\Ext^3(Q,E)\cong \Hom(E,Q)^\ast$. We thus obtain
\[
\HH^0(\BP^3,\lExt^2(E,E)) \cong \HH^0(\BP^3,\lExt^3(Q,E)) \cong \Hom(E,Q)^\ast,
\]
as desired. Finally, the isomorphism $\lExt^2(E,E) \simto \lExt^3(Q,E)$ combined with \Cref{lemma:supp-ext-TF} gives $\Supp\lExt^2(E,E) = \Supp Q$, proving \ref{ext-iii}.
\end{proof}

\begin{prop}
\label{lem:dimext2}
Let $E$ be a rank 2 semistable sheaf with Chern classes $(c_1,c_2,c_3)=(-1,1,1-2n)$ for some $n \geq 1$. Assume that $Q \defeq E^{\ast\ast}/E = \OO_Z$ for some reduced 0-dimensional subscheme $Z \subset \BP^3$. Then 
\[
\ext^2(E,E)=4n.
\]
\end{prop}

\begin{proof}
Since by \Cref{lem:exts}\ref{ext-iii} the sheaf $\lExt^2(E,E)$ is 0-dimensional and has support equal to $Z$, we have
\[
\lExt^2(E,E) = \bigoplus_{x\in Z} \lExt^2(E,E)_x
\]
and thus, by \Cref{lem:exts}\ref{ext-ii},
\begin{equation}
\label{eqn:ext2-decomposition}
\Ext^2(E,E) = \HH^0(\lExt^2(E,E)) = \bigoplus_{x\in Z} \Ext^2_{\OO_{\BP^3,x}}(E_x,E_x).
\end{equation}

Let $p\in\BP^3$ denote the unique singular point of the reflexive sheaf $R = E^{**}$, and assume first that $p\notin Z$. Let $x \in Z$ be a point. Localising $E \into R \onto \OO_Z$ at $x$ yields
\[
\begin{tikzcd}
0 \arrow{r} & E_x \arrow{r} & R_x \arrow{r} & \kappa(x) \arrow{r} & 0,
\end{tikzcd}
\]
where $\kappa(x) \cong \BC$  is the residue field at $x$ and $E_x \cong \mathfrak m_x \oplus \OO_{\BP^3,x}$. Here, $\mathfrak m_x \subset \OO_{\BP^3,x}$ is the maximal ideal, which coincides with the localisation of the ideal sheaf $\SI_x \subset \OO_{\BP^3}$ at $x$. Since $\Ext^{>0}(\OO_{\BP^3,x},-) = 0$,
\[
\Ext^2_{\OO_{\BP^3,x}}(E_x,E_x) = \Ext^2_{\OO_{\BP^3,x}}(\mathfrak m_x,\mathfrak m_x)\oplus \Ext^2_{\OO_{\BP^3,x}}(\mathfrak m_x,\OO_{\BP^3,x}).
\]
Since $\lExt^{>0}(\SI_x,-)|_{\BP^3 \setminus \set{x}} = 0$ and $\HH^{>0}(\BP^3,\lHom(\SI_x,\SI_x))=0$, we have, again using the spectral sequence for local-to-global exts, 
\[
{\Ext^2_{\OO_{\BP^3,x}}}{(\mathfrak m_x,\mathfrak m_x)}
\cong \lExt^2(\mathscr I_x,\mathscr I_x)_x \cong \HH^0(\BP^3,\lExt^2(\SI_x,\SI_x))\cong \Ext^2(\SI_x,\SI_x).
\]
On the other hand, $\Ext^2(\SI_x,\SI_x)$ is dual to $\Ext^1(\SI_x,\SI_x)$, see e.g.~\cite[Sec.~2.1.1]{Ricolfi-Virtual-Classes-Quot} or \cite[Sec.~11.2.2]{Ricolfi-Book}. Since $\ext^1(\SI_x,\SI_x) = 3$, this summand contributes 3 dimensions. Similarly, for the second summand, we have  
\[
\Ext^2_{\OO_{\BP^3,x}}(\mathfrak m_x,\OO_{\BP^3,x}) \cong \Ext^2(\SI_x,\OO_{\BP^3}),
\]
which is Serre dual to $\HH^1(\BP^3,\SI_x(-4)) \cong \HH^0(\BP^3,\kappa(x)) = \BC$. In conclusion, each $x \in Z$ contributes 4 dimensions, whence $\ext^2(E,E) = 4n$ in the case $Z \cap \Sing R = \emptyset$.

Now let $Z=p \amalg Z'$ with $Z'$ reduced and not containing $p$. We consider the commutative diagram
\[
\begin{tikzcd}
0\arrow{r} & E\arrow{r}\arrow{d} & R\arrow{r}\arrow[equal]{d} & \OO_{p}\oplus\OO_{Z'}\arrow[two heads]{d}\arrow{r} & 0 \\
0\arrow{r} & E'\arrow{r} & R\arrow{r} & \OO_{Z'}\arrow{r} & 0
\end{tikzcd}
\]
where $E'\defeq \ker(R\onto\OO_{p}\oplus\OO_{Z'}\onto\OO_{Z'})$. By the snake lemma, $E'$ fits into the exact sequence
\[
\begin{tikzcd}
0 \arrow{r} & E \arrow{r} & E' \arrow{r} & \OO_{p} \arrow{r} & 0.
\end{tikzcd}
\]
It follows that $E_{x}\cong  E_{x}'$ for all $x \in Z'$. We then get, using \Cref{eqn:ext2-decomposition}, that
\[
\Ext^2(E,E) =  \Ext^2_{\OO_{\BP^3,p}}(E_p,E_p) \oplus \bigoplus_{x\in Z'}\Ext^2_{\OO_{\BP^3,x}}(E_x',E_x').
\]
The second summands have been computed above: each has dimension 4. We now confirm that the first summand $\Ext^2_{\OO_{\BP^3,p}}(E_p,E_p)$ also has dimension 4. By \Cref{lem:exts}\ref{ext-ii}, it is enough to prove that for a quotient $R \onto \OO_p$ with kernel $E$, one has $\hom(E,\OO_p)=4$, where $R=E^{\ast\ast}$ is singular at $p$. Consider the exact sequence
\[
\begin{tikzcd}[column sep=small]
    0 \arrow{r} & \Hom(\OO_p,\OO_p) \arrow{r} & \Hom(R,\OO_p) \arrow{r} & \Hom(E,\OO_p)  \arrow{r} & \Ext^1(\OO_p,\OO_p) \arrow{r}{w} & \Ext^1(R,\OO_p).
\end{tikzcd}
\]
We have that $\hom(R,\OO_p)=3$ by \Cref{lemma:homR-O}\ref{red-pt}. Since $\hom(R,\OO_p)-\ext^1(R,\OO_p)=\chi(R,\OO_p)=2\cdot \chi(\OO_p) = 2$, we deduce $\ext^1(R,\OO_p)=1$. Moreover, $\hom(\OO_p,\OO_p)=1$ and $\ext^1(\OO_p,\OO_p)=3$. Thus
\[
\hom(E,\OO_p) = 5-\dim \im w.
\]
Since $\ext^1(R,\OO_p)=1$, we have to prove that $w$ is nonzero. The map $R \onto \OO_p$ is equivalent to a nonzero linear functional $\kappa(p)^{\oplus 3} \onto \kappa(p)$, determined by a nonzero tuple $\boldit a = (a_1,a_2,a_3) \in \BC^{3}$. Set $A = \widehat{\OO}_{\BP^3,p} = \BC\llbracket x,y,z \rrbracket$. The Koszul resolution of $\OO_p$, given by 
\[
\begin{tikzcd}[column sep=large]
   A^{\oplus 3} \arrow{r}{(x~y~z)} & A \arrow{r} & \OO_p \arrow{r} & 0,
\end{tikzcd}
\]
is compatible with the resolution of $R$, i.e.~the diagram 
\[
\begin{tikzcd}[column sep=large]
  0\arrow{r} & A \arrow{r}{(x~y~z)^{\mathsf t}}\arrow{d}{\boldit a^{\mathsf t}} & A^{\oplus 3} \arrow{r}\arrow{d}{\boldit a} & R \arrow{r}\arrow[two heads]{d} & 0 \\
   A^{\oplus 3} \arrow{r} & A^{\oplus 3} \arrow{r}{(x~y~z)} & A \arrow{r} & \OO_p \arrow{r} & 0 
\end{tikzcd}
\]
commutes. It follows that the map $w \colon \Ext^1(\OO_p,\OO_p) \to \Ext^1(R,\OO_p)$, computed through these compatible resolutions, is identified with the $\kappa(p)$-linear map $(b_1,b_2,b_3)\mapsto a_1b_1+a_2b_2+a_3b_3$, and is therefore nonzero.
\end{proof} 

\begin{remark}
The ideas in the proofs of Lemmas \ref{lem:exts} and \ref{lem:dimext2} were already present in \cite[Prop.~7]{AJT2022} and \cite[Prop.~6]{JMT} under different hypotheses.
\end{remark}

We conclude this section by computing the virtual dimension of the moduli space $M_{c_3}$.

\begin{lemma}
\label{thm:tangent_space}
Let $[E] \in M_{c_3}$ be a closed point. Then, $\ext^1(E,E) - \ext^2(E,E) = 3$.
\end{lemma}

\begin{proof}
Since $\BP^3$ is Fano and $E$ is stable, we have $\Ext^3(E,E)=0$ and $\Hom(E,E)=\BC$. Set $R=E^{\ast\ast}$. An easy calculation with Hirzebruch--Riemann--Roch\footnote{Alternatively, one could use Hartshorne's results \cite{Hart_StableReflexive}, in particular the vanishings $\Ext^{>1}(R,R)=0$ to confirm that $\chi(E,E) = \chi(R,R) = -2$.} shows that $\chi(E,E)=-2$, as claimed.
\end{proof}

\begin{remark}
By the seminal work of Thomas \cite{ThomasThesis}, $M_{c_3}$ admits a perfect obstruction theory in the sense of \cite{BFinc}. The \emph{virtual dimension} of the obstruction theory is
\[
\vd = \ext^1(E,E)-\ext^2(E,E) = 3,
\]
by \Cref{thm:tangent_space}. Since $\dim N = 3 = \vd$, \Cref{MAIN-THM-a}\ref{MTHM2} (proved in \Cref{thm:iso_relative_quot_M} below) combined with Donaldson--Thomas theory for $\BP^3$ suggests the existence of a 0-dimensional perfect obstruction theory on $\Quot_{\BP^3}(R,n)$, for every $[R] \in N$, yielding a virtual fundamental class
\begin{equation}
\label{VFC?}
[\Quot_{\BP^3}(R,n)]^{\vir} \,\in\,A_0(\Quot_{\BP^3}(R,n)).
\end{equation}
We also note that Mizyuk \cite{Mizyuk-APOT} constructed a 0-dimensional almost perfect obstruction theory on the Quot scheme of finite quotients of any locally free sheaf on any smooth 3-fold. In particular, her results yield a virtual fundamental class
\begin{equation}
\label{VFC-mizyuk}
[\Quot_{\BP^3}(\OO_{\BP^3}(-1)^{\oplus 3},n)]^{\vir} \,\in\,A_0(\Quot_{\BP^3}(\OO_{\BP^3}(-1)^{\oplus 3},n)).
\end{equation}
In light of the closed embedding $\Quot_{\BP^3}(R,n) \into \Quot_{\BP^3}(\OO_{\BP^3}(-1)^{\oplus 3},n)$, cf.~\Cref{sec:zerolocus}, we believe it is an interesting problem to construct the class \eqref{VFC?} and understand its relation with \eqref{VFC-mizyuk}. We leave these problems for future research.
\end{remark}


\section{\texorpdfstring{Proof of \Cref{MAIN-THM-a}}{}}\label{sec:proofmainthma}

This section is devoted to settling the proof of \Cref{MAIN-THM-a}. We already established \Cref{MAIN-THM-a}\ref{MTHM1} in \Cref{prop:reflexive=P^3}. 
The isomorphism 
\[
\begin{tikzcd}
\Quot_{\BP^3 \times N/N}(\SR,n) \arrow{r}{\sim} & M_{c_3},
\end{tikzcd}
\]
announced in \Cref{MAIN-THM-a}\ref{MTHM2}, will be confirmed in \Cref{thm:iso_relative_quot_M} (the key technical input here is \Cref{thm:iso-with-sim}, for which, in turn, the key step is proved in Appendix \ref{app:Kuranishi}). The local triviality of the map $M_{c_3} \to N$, namely \Cref{MAIN-THM-a}\ref{MTHM31}, will be proven in \Cref{prop:pi-is-ZLT}. Finally, a proof of the connectedness of $M_{c_3}$, as stated in \Cref{MAIN-THM-a}\ref{MTHM3}, will be given, as a consequence of \Cref{MAIN-THM-a}\ref{MTHM2}, in \Cref{cor:connectedess-of-M}.

\subsection{DT/PT wall-crossing and the moduli space at the wall}
\label{sec:moduli-wall}

We start by recalling some results from \cite{JLMM25,Mihai-Tuomas}.
Let $(Y,H)$ be a smooth complex projective 3-fold of Picard rank 1 admitting Bridgeland stability conditions obtained via the double tilt procedure of \cite{BMT14}. Let $v \in \mathrm K_{\num}(Y)$ be a numerical class such that $\rk(v)>0$ and is coprime to $H^2\cdot \ch_1(v)$, so that $\mu$-stability coincides with Gieseker stability. It is shown in \cite{JLMM25} that both the moduli space $M_Y^{\st}(v)$ of $\mu$-stable sheaves in class $v$ and the moduli space $M_Y^{\PT}(v)$ of PT-stable objects can be realised as Bridgeland moduli spaces. Moreover, there is a single wall separating the two chambers, giving rise to a diagram
\begin{equation}
\label{WC-diagram}
\begin{tikzcd}
M_Y^{\st}(v) \arrow[swap]{dr}{\Phi} & & M_Y^{\PT}(v)
\arrow{dl}{\Psi} \\
& \BM_{\sigma}(v)
\end{tikzcd}
\end{equation}
where $\sigma \in \Stab(Y)$ is a suitable stability condition on the wall (chosen accordingly to \cite[Sec.~6.2]{Mihai-Tuomas}) and $\BM_{\sigma}(v)$ is the good moduli space\footnote{Such a good moduli space exists by \cite[Ex.~7.29]{zbMATH07770154}.} of Bridgeland's moduli stack $\FM_{\sigma}(v)$ of $\sigma$-semistable objects in class $v$.

The maps $\Phi$ and $\Psi$ come from open embeddings $\overline \Phi$ and $\overline \Psi$ at stack level into the stack $\FM_\sigma(v)$, and can be explicitly described on $\BC$-valued points.\footnote{We will not use the PT side for our arguments, only the existence of the morphism $\Phi$. However, for completeness, the morphism $\Psi$ is given by $\Psi (J) = [h^0(J) \oplus h^1(J)[-1]]$.} Let $E \in M_Y^{\st}(v)$ be a point corresponding to a $\mu$-stable sheaf $E$. Consider the diagram
\[
\begin{tikzcd}
E^{[2]} \arrow[hook]{dr}&& Z_E\arrow[hook]{d} \\
 E \arrow[hook]{r}\arrow[hook]{u} & E^{\ast\ast} \arrow[two heads]{r}\arrow[two heads]{dr} & Q_E \arrow[two heads]{d}  \\
&& T_E 
\end{tikzcd}
\]
where $Z_E \subset Q_E$ is the maximal 0-dimensional subsheaf of $Q_E$, so that $T_E$ is pure of dimension 1 (or trivial). The subsheaf $E^{[2]} \into E^{\ast\ast}$ is the kernel of the composite surjection $E^{\ast\ast} \onto T_E$. One has
\begin{align*}
\Phi (E) &= [E^{[2]}\oplus Z_E[-1]]\,\in\,\BM_\sigma(v),
\end{align*}
where $[-]$ denotes the S-equivalence class.
When $\dim Q_E = 0$ we have $T_E = 0$ and hence $E^{[2]} = E^{\ast\ast}$, so that 
\[
\Phi(E) = [E^{\ast\ast} \oplus Q_E[-1]]\,\in\,\BM_\sigma(v).
\]

\begin{theorem}
[{\cite[Thm.~1.1]{Mihai-Tuomas}}]
\label{thm:proj-scheme}
Both $M_Y^{\PT}(v)$ and $\BM_{\sigma}(v)$ are projective schemes.
\end{theorem}

The results of \cite[Sec.~6]{Mihai-Tuomas} prove that Bridgeland's moduli stack  $\FM_{\sigma}(v)$ is equivalent to the stack $\FM_Y(v)$ defined as follows. Consider the torsion pair $(\Coh_{\leq 1}(Y),\Coh_{\geq 2}(Y))$ on $\Coh Y$, where 
\begin{align*}
\Coh_{\leq 1}(Y) &= \Set{T \in \Coh Y | \dim T \leq 1}, \\
\Coh_{\geq 2}(Y) &= \Set{F\in \Coh Y|\dim F' \geq 2\mbox{ for all }0\neq F'\subset F}
\end{align*} 
are full subcategories of $\Coh Y$. There is an induced abelian heart
\[ \CA^p(Y) = \braket{\Coh_{\geq 2}(Y),\Coh_{\leq 1}(Y)[-1]} \subset \derived^{\bounded}(Y) \]
whose stack of objects $\CD^p$ is a well-defined open substack of Lieblich's moduli stack $\FM^{\pug}(Y)$ of semi-Schur objects \cite{Lieblich1}. In particular, $\CD^p$ is algebraic and locally of finite type over $\BC$. Now, for a scheme $S$, the groupoid of $S$-objects of $\FM_Y(v)$ is defined as the groupoid of $S$-perfect complexes $\BE \in \derived^{\bounded}(Y\times S)$ such that for every geometric point $s\in S$ the derived restriction $\BE_s$ to the point $s$ is an object of $\CD^p$ in class $v$, such that 
\begin{enumerate}
\item $h^0(\BE_s)$ is 2-Gieseker semistable, and
\item $h^1(\BE_s)$ is 0-dimensional.
\end{enumerate}

\begin{situation}
\label{situation-moduli-identified}
We apply the general theory described so far to $Y=\BP^3$ and $v$ the class corresponding to $c=(2,-1,1,c_3)$. In this case, the Gieseker moduli scheme is $M_{c_3}$ and 2-Gieseker stability of $h^0(\BE_s)$ is equivalent to Gieseker stability. We shall denote
\[
\FM \defeq \FM_{\BP^3}(v),\qquad \FM_\sigma \defeq \FM_\sigma(v),\qquad \BM \defeq \BM_{\BP^3}(v),\qquad \BM_\sigma \defeq \BM_\sigma(v),
\]
and exploit the equivalence $\FM_\sigma \cong \FM$ to identify, canonically, $\BM_\sigma$ with the good moduli space $\BM$ of $\FM_{\BP^3}(v)$.
\end{situation}

\begin{remark}
\label{rmk:u-bijective}
Let us put ourselves in the context of \Cref{situation-moduli-identified}. The stack $\FM_{\sigma}$ parametrises objects $\BE$ which fit in an exact triangle
\begin{equation}
\label{triangle-h01}
\begin{tikzcd}
    h^0(\BE) \arrow{r} & \BE \arrow{r} & h^1(\BE)[-1]
\end{tikzcd}
\end{equation}
where $h^0(\BE)$ is stable and $h^1(\BE)$ is 0-dimensional.
In fact, when passing to the corresponding good moduli space $\BM_{\sigma}$, more can be said: by \cite[Prop.~6.5]{Mihai-Tuomas}, every point is S-equivalent to an object of the form $R\oplus Q[-1]$, where $R$ is of homological dimension at most 1 and $Q$ is 0-dimensional. In fact, for our purposes, we need to make the following easy but crucial observation: the sheaves $R$ are parametrised by the moduli space $N$ of reflexive sheaves.
By \cite[Cor.~6.3]{Mihai-Tuomas}, a closed polystable representative is a direct sum of 2-Gieseker stable sheaves of positive rank, with $\hd \leq 1$, and shifted skyscraper sheaves. Coprimality of rank 2 and degree $-1$ forces a single such positive rank summand $R$. It has $c_1=-1$ and $c_2=1$. \Cref{lemma:length-quotient} gives that $R^{**}/R$ is 0-dimensional, whereas \Cref{rmk:refl-on-3folds} and $\hd \leq 1$ force this quotient to vanish. Hence $R$ is reflexive and has $c_3=1$. Additivity of the Chern character gives $\chi(Q)=n$.
\end{remark}

We are ready to provide the following description of the good moduli space $\BM \cong \BM_{\sigma}$ of the stack $\FM \cong \FM_{\sigma}$. 

\begin{theorem}
\label{thm:iso-with-sim}
There is an isomorphism of projective schemes
\[
\begin{tikzcd}
u\colon N \times \Sym^{n}(\BP^3) \arrow{r}{\sim} & \BM.
\end{tikzcd}
\]
\end{theorem}

\begin{proof}
Consider the diagram
\[
\begin{tikzcd}
\SE_{\CN}\arrow[dash]{d} &\BP^3\times \CN \times \CCoh^{n}(\BP^3)\arrow[swap]{dl}{\pr_{12}}\arrow{dr}{\pr_{13}} & \SE_{n} \arrow[dash]{d} \\
    \BP^3\times \CN && \BP^3\times\CCoh^{n}(\BP^3) &
\end{tikzcd}
\]
where $\CCoh^{n}(\BP^3)$ is the moduli stack of 0-dimensional sheaves on $\BP^3$ of length $n$ (with universal sheaf $\SE_{n}$) and $\CN$ is the moduli stack of stable sheaves in class $c = (2,-1,1,1)$ (with universal reflexive sheaf $\SE_{\CN}$). The complex
\[
\BE = \SE_{\CN} \boxplus \SE_{n}[-1] = \pr_{12}^\ast \SE_{\CN} \oplus \pr_{13}^\ast \SE_{n}[-1]
\]
defines a modular morphism
\[
\begin{tikzcd}
\CN \times  \CCoh^{n}(\BP^3)
\arrow{r}{\overline u} & \FM,
\end{tikzcd}
\]
sending $(R,Q) \mapsto E=R\oplus Q[-1]$. Passing to the good moduli spaces, we obtain a commutative diagram 
\[
\begin{tikzcd}
\CN \times  \CCoh^{n}(\BP^3)
\arrow{r}{\overline u}\arrow[swap]{d} & \FM \arrow{d} \\
N \times  \Sym^{n}(\BP^3)
\arrow{r}{u} & \BM
\end{tikzcd}
\]
where $u$ is the good moduli space morphism attached to $\overline u$ and the vertical morphisms are the good moduli space maps. The map $u$ sends $(R,[Q]) \mapsto [R\oplus Q[-1]]$. This map is bijective on $\BC$-valued points by the results of \cite[Sec.~6]{Mihai-Tuomas} and \Cref{rmk:u-bijective}. Since both the source and the target of $u$ are finite type $\BC$-schemes (the target thanks to \Cref{thm:proj-scheme}), to conclude that it is an isomorphism, it is enough to prove that it is \'etale. This is the content of Appendix \ref{app:Kuranishi}.
\end{proof}

Now, exploiting the equivalence $\BM \cong \BM_{\sigma}$, we can consider the composition
\begin{equation}
\label{eqn:key-maps}
\begin{tikzcd}
\delta\colon M_{c_3} \arrow{r}{\Phi} & \BM_\sigma \arrow{r}{\sim} & N \times \Sym^{n}(\BP^3)\arrow{r}{\pr_1} & N.
\end{tikzcd}
\end{equation}
This morphism is given, on closed points, by taking the double dual.


\subsection{The isomorphism with the relative Quot scheme}
\label{sec:iso-quot-M}

Fix an odd $c_3\in\BZ_{\leq -1}$ and set, as ever, $n=\frac{1-c_3}{2}$. We use the universal family of reflexive sheaves
\[
\begin{tikzcd}
\mathscr R \arrow[dash]{d} & \\
\BP^3 \times N \arrow{r}{p_N} & N
\end{tikzcd}
\]
to form the relative Quot scheme of points
\[
\begin{tikzcd}
\Quot_{p_N}(\mathscr R,n) \arrow{r}{\pi} & N.
\end{tikzcd}
\]
The fibre of $\pi$ over the closed point of $N$ corresponding to the sheaf $R$ is the Quot scheme $\Quot_{\BP^3}(R,n)$. 
We shall denote by 
\begin{equation}
\label{eqn:universal-sequence}
\begin{tikzcd}
0 \arrow{r} &
\mathscr F_{c_3} \arrow{r}{\eta} &
(\id_{\BP^3} \times \pi)^\ast \mathscr R \arrow{r} &
\SQ
\arrow{r} & 0
\end{tikzcd}
\end{equation}
the universal short exact sequence living over
\[
(\BP^3 \times N)\times_{N} \Quot_{p_N}(\mathscr R,n) = \BP^3 \times \Quot_{p_N}(\mathscr R,n).
\]
Now, $M_{c_3}$ is itself a fine moduli space by \Cref{lemma:M_is_fine}, so to give a flat family of stable sheaves parametrised by a scheme $B$ (up to a line bundle pulled back from $B$) is the same as to give a morphism $B \to M_{c_3}$. The universal kernel $\mathscr F_{c_3}$ defines a flat family of stable sheaves on $\BP^3$ parametrised by the relative Quot scheme (reason for flatness: $\SQ$ is flat over $\Quot_{p_N}(\mathscr R,n)$ by definition of the Quot functor, and $(\id_{\BP^3} \times \pi)^\ast \mathscr R$ is also flat, since $\mathscr R$ is flat over $N$), thus we obtain a morphism
\begin{equation} 
\label{eq:fc3}
\begin{tikzcd}
f_{c_3}\colon \Quot_{p_N}(\mathscr R,n) \arrow{r} & M_{c_3},
\end{tikzcd}
\end{equation}
defined on closed points by sending $(R,[\phi\colon R \onto Q]) \mapsto \ker(\phi)$.

In fact, we can already prove the connectedness of $M_{c_3}$.

\begin{theorem}
\label{cor:connectedess-of-M}
Fix an odd $c_3 \in \BZ_{\leq 1}$. Then $M_{c_3}$ is connected.
\end{theorem}

\begin{proof}
The morphism $f_{c_3}$ is bijective. Injectivity follows from the definition of $f_{c_3}$ and the universal property of the Quot functor. As for surjectivity, a point $[E] \in M_{c_3}$ is the image of the point $(E^{\ast\ast},[E^{\ast\ast} \onto E^{\ast\ast}/E])$.
It is then enough to prove that $\Quot_{p_N}(\SR,n)$ is connected, where $n = \frac{1-c_3}{2}$. The structure morphism $\pi\colon \Quot_{p_N}(\SR,n) \to N$ is proper and surjective, and its fibres are connected by \Cref{prop:Quot-connected}. Thus the result follows from \Cref{lemma:conn-irr}\ref{CON}.
\end{proof}

The proof of \Cref{MAIN-THM-a}\ref{MTHM2}, namely that $f_{c_3}$ is an isomorphism, will be given as a combination of the constructions and results in the next two subsections: in \Cref{sec:map-to-RM} we construct an \emph{injective} map $i\colon \SE \to \SR_M \defeq (\id_{\BP^3}\times \delta)^\ast \SR$ of sheaves on $\BP^3 \times M_{c_3}$, whose cokernel will yield a morphism
\[
\begin{tikzcd}
g_{c_3}\colon M_{c_3} \arrow{r} & \Quot_{p_N}(\mathscr R,n).
\end{tikzcd}
\]
In \Cref{sec:check-inverse}, we will confirm that $g_{c_3}$ and $f_{c_3}$ are inverse to each other.

\subsubsection{\texorpdfstring{Constructing a morphism $g_{c_3}\colon M_{c_3}\to \Quot_{p_N}(\mathscr R,n)$}{}}
\label{sec:map-to-RM}
Throughout, we keep the notation introduced in \Cref{situation-moduli-identified}. We furthermore set $M \defeq M_{c_3}$ and we call $\CM$, resp.~$\CN$, the stack of stable sheaves corresponding to $M$, resp.~to $N$. We denote by $\SE_\bullet$ their universal families. Throughout, to simplify the notation, we omit pullbacks along projection maps, we omit deriving functors, and for a map $h\colon S \to T$ we still denote by $h$ the product map $\id_{\BP^3} \times h \colon \BP^3 \times S \to \BP^3 \times T$.

The constructions of \Cref{sec:moduli-wall} allow us to consider the diagram
\[
\begin{tikzcd}[column sep=large]
\SE_{\CM}\arrow[dash]{d} & \SE_{\FM_{\sigma}}\arrow[dash]{d} & \\ 
\BP^3\times\mathcal{M}\arrow{r}{\overline \Phi}\arrow{d}{\alpha} & \BP^3 \times \FM_\sigma\arrow{d}{\beta} & \BP^3 \times \CN \times \CCoh^{n}(\BP^3)\arrow[swap]{l}{\overline u}\arrow{d}{\gamma} \\
\BP^3\times M\arrow{r}{\Phi} & \BP^3 \times \BM_\sigma & \BP^3 \times N \times \Sym^{n}(\BP^3)\arrow[l,swap,"u","\sim"']
\end{tikzcd}
\]
where the vertical arrows are good moduli space morphisms. 
By definition of $\overline{\Phi}$ and $\overline u$ we have, respectively,
\[
\overline{\Phi}^*\SE_{\FM_{\sigma}}\cong  \SE_{\mathcal{M}},\qquad \overline{u}^*\SE_{\FM_{\sigma}}\cong  \SE_{\mathcal{N}}\boxplus \SE_{n}[-1].
\]
By \cite[Cor.~4.6.6]{modulisheaves}, there exists a line bundle $\SL_{\mathcal{M}}$ on $\mathcal{M}$ such that stabilisers at closed points act on the fibres with weight $-1$, so that $\SL_{\mathcal{M}}\otimes\SE_{\mathcal{M}}$ descends to a universal family $\SE $ on $\BP^3\times M$, i.e. $\alpha^*\SE \cong \SL_{\mathcal{M}}\otimes\SE_{\mathcal{M}}$. Similarly, there exists a line bundle $\SL_{\mathcal{N}}\in\Pic(\mathcal{N})$ such that $\SL_{\mathcal{N}}\otimes\SE_{\mathcal{N}}$ descends to a universal family $\mathscr{R}$ on $\BP^3\times N$, i.e. $\gamma^*\mathscr{R}\cong  \SL_{\mathcal{N}}\otimes\SE_{\mathcal{N}}$.

Let 
\[
V\defeq \Set{(x,R,[Q])\in \BP^3\times N\times\Sym^{n}(\BP^3)|x\notin\Sing R\cup \Supp Q},
\]
and define $U\defeq u(V)$. Since $u$ is an isomorphism by \Cref{thm:iso-with-sim}, $U$ is open. We can pull back the previous diagram via the open immersion $U\into\BP^3 \times\BM_{\sigma}$, and we obtain
\begin{equation}
\label{eqn:triple-comm-diag}
\begin{tikzcd}[column sep=large]
	\mathcal{W}\arrow{d}{\alpha}\arrow{r}{\overline \Phi} & \mathcal{U}\arrow{d}{\beta} & \mathcal{V}\arrow{d}{\gamma}\arrow[swap]{l}{\overline u} \\
	W\arrow{r}{\Phi} & U & V \arrow[l,swap,"u","\sim"']
\end{tikzcd}
\end{equation}
By restriction, we obtain relations
\begin{align}
\overline{\Phi}^*\SE_{\FM_{\sigma}}|_\mathcal{U} &\cong  \SE_{\mathcal{M}}|_\mathcal{W},
\label{eq:relpullbackrestricted1}
\\
\overline{u}^*\SE_{\FM_{\sigma}}|_\mathcal{U} &\cong  \left(\SE_{\mathcal{N}}\boxplus \SE_{n}[-1]\right)|_\mathcal{V}=\SE_{\mathcal{N}}|_{\mathcal{V}}.\label{eq:relpullbackrestricted2}
\end{align}
We further notice that $\SE_{\mathcal{N}}|_{\mathcal{V}}$ is locally free.

Our aim is to descend \Cref{eq:relpullbackrestricted1} and \Cref{eq:relpullbackrestricted2} at the level of good moduli spaces, in order to obtain a relation of the form
\begin{equation}
\SE |_{W}\otimes \SK|_W\cong  (u^{-1}\circ \Phi)^*\mathscr{R}|_V, 
\end{equation}
for some $\SK\in\Pic(M)$. To achieve this, we tensor $\SE_{\FM_{\sigma}}|_\mathcal{U}$ with a suitable line bundle in such a way that it descends to $U$. We follow ideas from \cite{Tajakka,Mihai-Tuomas}.

As explained in \cite[Sec.~4]{Mihai-Tuomas}, $\FM_\sigma$ carries a natural class of line bundles, called \textit{determinantal line bundles}. Indeed, for a vector bundle $F$ on $\mathbb{P}^3$, we can define a line bundle
\[
\lambda(F) \defeq \mathrm{det}\left(\RR\pr_{\FM_{\sigma}*}\left(\SE_{\FM_{\sigma}}\otimes \pr_{\BP^3}^\ast F\right)\right).
\]
For our purposes, we consider
\begin{equation}
\label{eqn:our-line-bundle}
    \SL_\sigma\defeq \lambda(\OO_{\BP^3})^{-1}\otimes \lambda(\OO_{\BP^3}(-1)^{\oplus 2})\otimes\lambda(\OO_{\BP^3}(-2))^{-1}.
\end{equation}
Let $E\defeq R\oplus \bigoplus_{1\leq a\leq k} W_a\otimes\OO_{q_a}[-1]$ represent a closed point of $\FM_{\sigma}$. By \cite[Prop.~4.2]{Tajakka}, an element $(t,f_1,\dots,f_k)\in\Aut(E)=\mathbb{G}_m\times\prod_{1\leq a\leq k}\operatorname{GL}(W_a)$ acts on the fibre of $\lambda(F)$, for a vector bundle $F$, by multiplication by
\[
t^{\chi(R\otimes F)}\det(f_1)^{-\rk(F)}\cdots\det(f_k)^{-\rk(F)}.
\]
Consequently, since $\chi(R)=\chi(R(-1))=0$ and $\chi(R(-2))=1$, the characters associated to the three determinantal line bundles in \eqref{eqn:our-line-bundle} are given in the following table.
\begin{equation*}
\renewcommand{\arraystretch}{1.1}
\begin{array}{@{}l@{\qquad}l@{}}
\text{\bfseries Sheaf} & \text{\bfseries Character} \\
\midrule[1.2pt]
\lambda(\OO_{\BP^3})
& \det(f_1)^{-1}\cdots\det(f_k)^{-1} \\
\midrule[0.4pt]
\lambda(\OO_{\BP^3}(-1)^{\oplus 2})
& \det(f_1)^{-2}\cdots\det(f_k)^{-2} \\
\midrule[0.4pt]
\lambda(\OO_{\BP^3}(-2))
& t\det(f_1)^{-1}\cdots\det(f_k)^{-1} \\
\bottomrule[0.6pt]
\end{array}
\end{equation*}
Thus the action on the fibre of $\SL_{\sigma}$ is multiplication by $t^{-1}$.

Notice that $\SE_{\FM_{\sigma}}|_{\mathcal{U}}$ is a vector bundle, since it is $\FM_{\sigma}$-perfect and its geometric fibres are locally free sheaves concentrated in degree zero. Consider the vector bundle
\begin{equation}
\left(\SE_{\FM_{\sigma}}\otimes\SL_\sigma\right)\big|_{\mathcal{U}}.
\end{equation}
By \cite[Thm.~10.3]{Alper-gms-artin-stacks}, it descends to a vector bundle on $U$ if and only if it has trivial stabiliser action at closed points. Let $y=(x,R\oplus \bigoplus_{1\leq a\leq k} W_a\otimes\OO_{q_a}[-1])\in\mathcal{U} \subset \BP^3 \times \FM_\sigma$ be a closed point. Since $x\notin\Sing R\cup\{q_1,\dots,q_k\}$, we have that
\[
\SE_{\FM_{\sigma}}|_y\cong  R|_x.
\]
Thus $(t,f_1,\dots f_k)\in\Aut(E)$ acts on this fibre by multiplication by $t$. The calculations above show that $(t,f_1,\dots f_k)\in\Aut(E)$ acts on the fibre $\SL_\sigma|_y$ by multiplication by $t^{-1}$.
Hence, the tensor product $\SE_{\FM_{\sigma}}|_{\mathcal{U}}\otimes \SL_\sigma|_{\mathcal{U}}$ has trivial stabiliser action at closed points and it descends to a vector bundle $\mathcal{G}$ on $U$, i.e.~$\beta^*\mathcal{G}\cong  \SE_{\FM_{\sigma}}|_{\mathcal{U}}\otimes \SL_\sigma|_{\mathcal{U}}$.

By \cite[Thm.~10.3]{Alper-gms-artin-stacks} again, the line bundle $\overline{\Phi}^*\SL_\sigma\otimes\SL^{-1}_{\mathcal{M}}\in\Pic(\mathcal{M})$ descends to a line bundle $\SK_{1}\in\Pic(M)$, and similarly $\overline{u}^*\SL_\sigma\otimes \SL^{-1}_{\mathcal{N}}\in\Pic(\mathcal{N}\times\CCoh^n(\BP^3))$ descends to a line bundle $\SK_{2}\in\Pic(N\times\Sym^n(\mathbb{P}^3))$. In particular, $\alpha^*\SK_1\cong \overline{\Phi}^*\SL_\sigma\otimes\SL^{-1}_{\mathcal{M}}$ and $\gamma^*\SK_2\cong \overline{u}^*\SL_\sigma\otimes \SL^{-1}_{\mathcal{N}}$.

Tensoring \Cref{eq:relpullbackrestricted1} with the line bundles defined above, we obtain 
\begin{equation}
\label{eqn:belle-cose}
    \overline{\Phi}^*\left(\SE_{\FM_{\sigma}}|_{\mathcal{U}}\otimes \SL_\sigma|_{\mathcal{U}}\right)\cong  \SE_{\mathcal{M}}|_\mathcal{W}\otimes\SL_{\mathcal{M}}|_\mathcal{W}\otimes\SL^{-1}_{\mathcal{M}}|_\mathcal{W}\otimes \overline{\Phi}^*\SL_\sigma|_\mathcal{W}.
\end{equation}
Using $\beta^*\mathcal{G}\cong  \SE_{\FM_{\sigma}}|_{\mathcal{U}}\otimes \SL_\sigma|_{\mathcal{U}}$ and the commutativity of the leftmost square in Diagram \eqref{eqn:triple-comm-diag}, the left hand side of \Cref{eqn:belle-cose} can be rewritten as
\[
\overline{\Phi}^*\beta^*\mathcal{G}   
 = \alpha^*\Phi^*\mathcal{G}.
\]
Using $\alpha^*\SE \cong \SL_{\mathcal{M}}\otimes\SE_{\mathcal{M}}$ and $\alpha^*\SK_1\cong \overline{\Phi}^*\SL_\sigma\otimes\SL^{-1}_{\mathcal{M}}$, the right hand side of \Cref{eqn:belle-cose} can be rewritten as
\[
\alpha^*(\SE |_{W})\otimes\alpha^*({\SK_1}|_{W}) = \alpha^*(\SE |_{W}\otimes{\SK_1}|_{W}).
\]
Pushing forward along $\alpha$ the relation just obtained, we get
\[
\alpha_{*}\alpha^*\Phi^*\mathcal{G}\cong  \alpha_*\alpha^*(\SE |_{W}\otimes{\SK_1}|_{W}).
\]
By \cite[Prop.~4.5]{Alper-gms-artin-stacks}, we obtain
\begin{equation}
\label{eq:needed}
\Phi^*\mathcal{G}\cong \SE |_{W}\otimes{\SK_1}|_{W}.
\end{equation}
With analogous reasoning, we obtain
\begin{align*}
    u^*\mathcal{G}&\cong  \mathscr{R}|_V\otimes \SK_{2}|_{V},\\
    \mathcal{G}&\cong (u^{-1})^*(\mathscr{R}|_V\otimes \SK_{2}|_{V}),
\end{align*}
which combined with \Cref{eq:needed}, yields
\begin{equation}
\label{almost-there}
\SE |_{W}\otimes (\SK_{1}\otimes(u^{-1}\circ\Phi)^*\SK^{-1}_{2})|_W\cong  (u^{-1}\circ\Phi)^*\left(\mathscr{R}|_V\right),
\end{equation}
as wanted.

We now construct a morphism $g_{c_3}\colon M\rightarrow\Quot_{p_N}(\mathscr R,n)$. Set $\mathscr{R}_{M}\defeq (\id_{\BP^3}\times \delta)^*\mathscr{R}$. We start replacing $\SE $ by $\SE \otimes \SK_{1}\otimes(u^{-1}\circ\Phi)^*\SK^{-1}_{2}$, so that the isomorphism \eqref{almost-there} reads
\begin{equation}
\label{eq:maprestrictedW}
\begin{tikzcd}
i_W\colon\SE |_{W}\arrow{r}{\sim} & \mathscr{R}_{M}|_{W}.
\end{tikzcd}
\end{equation}
The open immersion $j\colon W\hookrightarrow\BP^3\times M_{c_3}$ has complement of fibrewise codimension 3 over $M_{c_3}$ by \Cref{lem:cohomology}\ref{coh-4}. Therefore, by \cite[Prop.~3.5]{zbMATH02108016}, we have $j_*(\mathscr{R}_{M}|_{W})=\mathscr{R}_{M}$. As a consequence, by adjunction, the isomorphism in \Cref{eq:maprestrictedW} extends uniquely to a morphism
\[
\begin{tikzcd}
i\colon\SE \arrow{r} & \mathscr{R}_{M}.
\end{tikzcd}
\]
Its restriction to each fibre is generically an isomorphism (since $i$ extends $i_W$) and has a torsion-free source; hence, it is injective. We are then in a position to apply the following lemma.

\begin{lemma}
\label{Lemma:inj+flat}
Let $X \to B$ be a morphism of locally noetherian schemes, $i\colon \SE \to \SF$ a morphism of coherent sheaves over $X$ such that $\SF$ is $B$-flat  and the fibre $i_b\colon \SE_b \to \SF_b$ is a monomorphism in $\Coh X_b$ for all $b \in B$. Then $i$ is a monomorphism in $\Coh X$ and $\coker i = \SF/\SE$ is $B$-flat.     
\end{lemma}

\begin{proof}
Follows immediately from \cite[Prop.~11.3.7]{EGA-IV-3}. Alternatively, it can be obtained by combining the local criterion for flatness with \cite[\href{https://stacks.math.columbia.edu/tag/00ME}{Tag 00ME}]{stacks-project}.
\end{proof}

By \Cref{Lemma:inj+flat}, $i$ is injective and fits into an exact sequence
\begin{equation}
\label{eqn:ses-g}
\begin{tikzcd}
0\arrow{r} & \SE \arrow{r}{i} & \mathscr{R}_M\arrow{r} & \SR_M/\SE \arrow{r} & 0,
\end{tikzcd}
\end{equation}
such that the cokernel $\SR_M/\SE$ is flat over $M_{c_3}$. As a consequence, the exact sequence \eqref{eqn:ses-g} classifies a morphism
\[
\begin{tikzcd}
g_{c_3}\colon M_{c_3}\arrow{r} & \Quot_{p_N}(\mathscr R,n),
\end{tikzcd}
\]
which is our candidate inverse to $f_{c_3}$.

\subsubsection{Finishing the proof of \Cref{MAIN-THM-a}\ref{MTHM2}}
\label{sec:check-inverse}
To complete the proof of \Cref{MAIN-THM-a}\ref{MTHM2}, it remains to prove the following result.

\begin{theorem}
\label{thm:iso_relative_quot_M}
The morphism
\[
\begin{tikzcd}
f_{c_3}\colon \Quot_{p_N}(\mathscr R,n) \arrow{r} & M_{c_3},
\end{tikzcd}
\]
constructed in \eqref{eq:fc3} is an isomorphism, whose inverse is the morphism $g_{c_3}$.
\end{theorem}

We need the following lemma first. 

\begin{lemma}
\label{lemma:N-schemes}
The morphism $f_{c_3} \colon \Quot_{p_N}(\SR,n) \to M_{c_3}$ is a morphism of $N$-schemes, i.e.~we have $\delta\circ f_{c_3}=\pi$.
\end{lemma}

\begin{proof}
By the universal property of $M_{c_3}$, there exists a line bundle $\SL\in\Pic(\Quot_{p_N}(\mathscr R,n))$ and an isomorphism
\begin{equation*}
\begin{tikzcd}
h\colon\SF_{c_3}\arrow{r}{\sim} & (\id_{\BP^3} \times f_{c_3})^\ast\SE \otimes \SL.
\end{tikzcd}
\end{equation*}
Consider the open subset
\[
U'\defeq (\id_{\BP^3} \times f_{c_3})^{-1}(W)\subset\BP^3\times\Quot_{p_N}(\mathscr R,n),
\]
where $W \subset \BP^3 \times M_{c_3}$ is the open subset of points $(x,E)$ such that $x \notin \Sing E$, cf.~Diagram \eqref{eqn:triple-comm-diag}.
Exploiting \Cref{lem:cohomology}\ref{coh-4} and \cite[Prop.~3.5]{zbMATH02108016}, the composite isomorphism
\begin{equation*}
\begin{tikzcd}
\left((\id_{\BP^3} \times f_{c_3})^\ast i\otimes \id_{\SL}\right)|_{U'}\circ h|_{U'}\circ(\eta|_{U'})^{-1}\colon (\id_{\BP^3} \times \pi)^\ast \mathscr R|_{U'}\arrow{r}{\sim} & ((\id_{\BP^3} \times f_{c_3})^\ast \SR_M\otimes \SL)|_{U'}
\end{tikzcd}
\end{equation*}
extends uniquely to an isomorphism
\begin{equation*}
\begin{tikzcd}
\beta\colon (\id_{\BP^3} \times\pi)^\ast \mathscr R\arrow{r}{\sim} & (\id_{\BP^3} \times f_{c_3})^\ast \SR_M\otimes \SL.
\end{tikzcd}
\end{equation*}
By the universal property of $N$, this settles the equality $\delta\circ f_{c_3}=\pi$. 
\end{proof}

\begin{proofof}{\Cref{thm:iso_relative_quot_M}}
It is clear that $f_{c_3}\circ g_{c_3} = \id$ because the kernel of $\SR_M \onto \SR_M/\SE$ is the universal sheaf on $M_{c_3}$. It remains to prove that $g_{c_3}\circ f_{c_3} = \id$. The universal sequence on the Quot scheme is given in \eqref{eqn:universal-sequence}. On the other hand, the morphism $g_{c_3}\circ f_{c_3}$ is classified by the short exact sequence 
\[
\begin{tikzcd}[column sep=large]
0 \arrow{r} & (\id_{\BP^3} \times f_{c_3})^\ast\SE\arrow{r}{(\id_{\BP^3} \times f_{c_3})^\ast i} & (\id_{\BP^3} \times f_{c_3})^\ast \SR_M \arrow{r} &  (\id_{\BP^3} \times f_{c_3})^\ast (\SR_M/\SE)\arrow{r} & 0.
\end{tikzcd}
\]
Our goal is to show that $\im(\eta \colon \SF_{c_3} \to (\id_{\BP^3} \times \pi)^\ast \mathscr R)=\im  (\id_{\BP^3} \times f_{c_3})^\ast i$, after identifying $(\id_{\BP^3} \times f_{c_3})^\ast \SR_M$ with $(\id_{\BP^3} \times \pi)^\ast \mathscr R$, using the equality $\delta\circ f_{c_3}=\pi$ proved in \Cref{lemma:N-schemes}.
 
By uniqueness of the extension outside $U'$, we have the equality
\[
\beta\circ \eta=((\id_{\BP^3} \times f_{c_3})^\ast i\otimes\id_{\SL})\circ h.
\]
By \cite[Lemma 4.6.3]{modulisheaves}, the natural map
\[
\begin{tikzcd}
\OO_{\Quot_{p_N}(\mathscr R,n)}\arrow{r} & \pr_{\Quot\,\ast}\lEnd((\id_{\BP^3}\times\pi)^*\mathscr{R})
\end{tikzcd}
\]
is an isomorphism. In other words, for a Zariski open subset $\Spec A\subset\Quot_{p_N}(\mathscr R,n)$ where $\SL$ is trivial, the morphism $\beta|_{\BP^3\times \Spec A}$ is given by multiplication by a unit in $A$. Since $h$ is an isomorphism, this gives
$\im \eta=\im(\id_{\BP^3} \times f_{c_3})^\ast i$. This settles $g_{c_3}\circ f_{c_3} = \id$.
\end{proofof}

We have thus concluded the proof of the isomorphism in \Cref{MAIN-THM-a}\ref{MTHM2}.

\subsection{Proof of \Cref{MAIN-THM-a}\ref{MTHM31}}
We now prove the following property of the structure morphism of the relative Quot scheme.

\begin{prop}
\label{prop:pi-is-ZLT}
Fix an integer $n \geq 0$. The structure morphism
\[
\begin{tikzcd}
\Quot_{p_N}(\SR,n) \arrow{r}{\pi} & N
\end{tikzcd}
\]
is Zariski locally trivial. In particular, 
\begin{enumerate}
    \item [\mylabel{same-fibres}{\normalfont{(i)}}] the morphism $\pi$ is flat, and all the fibres are isomorphic to each other, and
    \item [\mylabel{irr-iff-irr}{\normalfont{(ii)}}] the moduli space $M_{c_3}$, where $c_3=1-2n$, is irreducible if and only if $\Quot_{\BP^3}(R,n)$ is irreducible, where $\Quot_{\BP^3}(R,n)$ is an arbitrary fibre of $\pi$.
\end{enumerate}
\end{prop}

\begin{proof}
Fix a closed point $p_0 = [R_0] \in N$ in the base. By the isomorphism $N \cong \BP^3$ of \Cref{prop:reflexive=P^3}, we have that $N$ is a homogeneous space of the form $G/P$, where $G=\PGL_4$ and $P$ is the parabolic subgroup $G_{p_0}\subset G$, namely the stabiliser of $p_0$. Let $\rho\colon G \to N$ be the quotient map. Denote by $\overline{\rho}$ the product map $\id_{\BP^3} \times \rho \colon \BP^3 \times G \to \BP^3 \times N$ and by $\SR_G$ the pullback of $\SR$ to $\BP^3 \times G$ along $\overline \rho$. By standard properties of relative Quot schemes, we have a cartesian diagram
\[
\begin{tikzcd}
\Quot_{p_G}(\SR_G,n) \MySymb{dr} \arrow{d}\arrow{r}{\pi_G} & G \arrow{d}{\rho} \\
\Quot_{p_N}(\SR,n)\arrow{r}{\pi} & N
\end{tikzcd}
\]
where $p_G\colon \BP^3 \times G \to G$ is the second projection. Now, the family $\SR_G$ is naturally homogeneous with respect to $G$: if we denote by $\alpha$ the $G$-isomorphism $\BP^3 \times G \simto \BP^3 \times G$ sending $(x,g) \mapsto (g^{-1}x,g)$, using that $G$ is reduced of finite type over $\BC$ and $N$ is separated, it is easy to see that the point-wise identification $\SR_G \otimes \kappa(g) = \SR \otimes \kappa(g\cdot p_0) = g_\ast R_0$ globalises to an isomorphism 
\begin{equation}
\label{equivariance}
\SR_G \cong \alpha^\ast \pr_{1}^\ast R_0 
\end{equation}
up to a line bundle pulled back from $G$, where $\pr_1\colon \BP^3 \times G \to \BP^3$ is the projection. Since quotients are insensitive to twisting by line bundles, we can use \eqref{equivariance} to obtain isomorphisms
\begin{equation}
\label{pullback-isos}
\begin{split}
\Quot_{p_N}(\SR,n) \times_NG 
&\cong \Quot_{p_G}(\SR_G,n) \\
&\cong \Quot_{p_G}(\alpha^\ast \pr_1^\ast R_0,n) \\
&\cong \Quot_{p_G}(\pr_1^\ast R_0,n) \\
&\cong G \times \Quot_{\BP^3}(R_0,n).
\end{split}
\end{equation}
Choose homogeneous coordinates in which the point corresponding to $p_0$ is $(1:0:0:0)$. Consider the standard affine opens $U_i \subset N$ given by $x_i \neq 0$, for $i=0,1,2,3$. For each $i$, there is a local section $\sigma_i\colon U_i \to G$ of $\rho$. For instance, writing a point of $U_0$ in the form $(1:a:b:c)$, the section $\sigma_0$ sends such a point to the class, in $G=\PGL_4$, of the matrix whose first column is $(1,a,b,c)^{\mathsf t}$ and whose next columns are the coordinate vectors $e_1$, $e_2$ and $e_3$. Consider the open subscheme $j_i\colon G_i \defeq G \times_NU_i \subset G$, so that $\sigma_i$ is a map $U_i \to G_i$. Denote by $\varphi$ the isomorphism of $G$-schemes
\[
\begin{tikzcd}
\Quot_{p_N}(\SR,n) \times_NG \arrow[r,"\varphi",,"\sim"'] & G \times \Quot_{\BP^3}(R_0,n)  
\end{tikzcd}
\]
obtained in \eqref{pullback-isos}. Pull it back along $j_i \circ \sigma_i\colon U_i \to G$ to obtain the isomorphism
\[
\begin{tikzcd}
(\Quot_{p_N}(\SR,n)\times_NG_i)\big|_{U_i} = \Quot_{p_N}(\SR,n) \times_NU_i \arrow[rr,"(j_i \circ \sigma_i)^\ast \varphi",,"\sim"'] & & U_i \times \Quot_{\BP^3}(R_0,n).
\end{tikzcd}
\]
Ranging over $i=0,1,2,3$, these isomorphisms provide a Zariski local trivialisation of $\pi$.

Finally, both \ref{same-fibres} and \ref{irr-iff-irr} are standard properties of Zariski fibrations, and \ref{irr-iff-irr} of course exploits the isomorphism $M_{c_3} \cong \Quot_{p_N}(\SR,n)$ proved in \Cref{thm:iso_relative_quot_M}.
\end{proof}

\section{\texorpdfstring{The smoothable component: Proof of \Cref{MAIN-gen-smooth-cpt}}{}}\label{sec:proofmainthmb}

In this section, we confirm that $M_{c_3}$ has a distinguished irreducible component 
\[
M_{c_3}^{\main} \subset M_{c_3},
\]
called the \emph{smoothable component}, which is always generically smooth, of dimension $5-2c_3$, cf.~\Cref{mainComponent}. This settles the proof of \Cref{MAIN-gen-smooth-cpt}\ref{sm-cpt-i}. Before that, in the next subsection, we provide a zero locus description of the Quot scheme of points $\Quot_{\BP^3}(R,n)$ for $[R] \in N = M_1$, yielding a lower bound for its dimension: this will be useful in determining the generic smoothness of the smoothable component. Part \ref{sm-cpt-ii}, resp.~\ref{sm-cpt-iii}, of \Cref{MAIN-gen-smooth-cpt} is proved in \Cref{Sing-main-Component}, resp.~in \Cref{sec:rationality-main-component}.

\subsection{Zero locus description for \texorpdfstring{$\Quot_{\BP^3}(R,n)$}{}}
\label{sec:zerolocus}

Let $R$ be a stable rank 2 reflexive sheaf on $\BP^3$ parametrised by the moduli space $N=M_1$. The resolution 
\begin{equation}
\label{res-r}
\begin{tikzcd}
0 \arrow{r} & \OO_{\BP^3}(-2) \arrow{r} & \OO_{\BP^3}(-1)^{\oplus3} \arrow{r} & R \arrow{r} & 0
\end{tikzcd}
\end{equation}
in \eqref{eqn:resolution-F} induces, for every $n \geq 0$, a closed immersion of Quot schemes of points
\[
\begin{tikzcd}
\Quot_{\BP^3}(R,n) \arrow[hook]{r}{\iota} & \Quot_{\BP^3}(L_0,n),
\end{tikzcd}
\]
where we have set $L_0 = \OO_{\BP^3}(-1)^{\oplus 3}$. Set also $L_1 = \OO_{\BP^3}(-2)$, and denote by $\SL_i$ the pullback of $L_i$ along the projection $\BP^3 \times \Quot_{\BP^3}(L_0,n) \to \BP^3$. Let $q_0\colon \SL_0 \onto \ST_0$ be the universal quotient living on $\BP^3 \times \Quot_{\BP^3}(L_0,n)$, and consider the composition
\[
\begin{tikzcd}
s_0\colon \SL_1\arrow[hook]{r} & \SL_0\arrow[two heads]{r}{q_0} & \ST_0.
\end{tikzcd}
\]
Let $\pr_2\colon \BP^3 \times \Quot_{\BP^3}(L_0,n) \to \Quot_{\BP^3}(L_0,n)$ be the second projection. The direct image
\[
\SV \defeq \pr_{2\ast} \lHom(\SL_1,\ST_0) = \pr_{2\ast} (\SL_1^{\ast}\otimes \ST_0)
\]
is a locally free sheaf of rank $n$ over $\Quot_{\BP^3}(L_0,n)$, and $s^{[n]}\defeq \pr_{2\ast} s_0$ is a global section of $\mathbf V \defeq \Spec \Sym \SV^\ast$.

The following lemma is a standard consequence of the universal property of cokernels and of the usual description of the functor of points of a zero locus of a section of a vector bundle.

\begin{lemma}
\label{lemma:zero-locus-section}
The Quot scheme $\Quot_{\BP^3}(R,n)$ is isomorphic to the scheme-theoretic vanishing locus of the section $s^{[n]} \in \Gamma(\Quot_{\BP^3}(L_0,n),\mathbf V)$.
\end{lemma}

Recall from \Cref{def:QUOT-circ} that we have introduced, for an arbitrary coherent sheaf $F$ on a variety $Y$, the open locus
\[ 
\Quot^{\circ}_{Y}(F,n) 
= \Set{ [F\onto Q] \in \Quot_{Y}(F,n) | 
\begin{array}{cc}
Q = \OO_Z \mbox{ with }    &  \\
 Z\subset Y \mbox{ reduced}  \end{array}
}\,\subset \,\Quot_Y(F,n). 
\]

\begin{prop}
\label{prop:Quot-generic-smoothness}
Let $R$ be a stable rank 2 reflexive sheaf on $\BP^3$ parametrised by the moduli space $N=M_1$. Then, the scheme $\Quot^{\circ}_{\BP^3}(R,n)$ is a smooth variety of dimension $4n$.
\end{prop}

\begin{proof}
The zero locus of the restricted section $s^\circ = s^{[n]}|_{\Quot^{\circ}_{\BP^3}(L_0,n)}$ fits in a cartesian diagram
\[
\begin{tikzcd}[row sep=large,column sep=large]
Z(s^\circ)\arrow[equal]{r} & \Quot_{\BP^3}^\circ(R,n) \MySymb{dr}\arrow[hook]{d}\arrow[hook]{r}{\iota^\circ} & \Quot_{\BP^3}^\circ(L_0,n)\arrow[hook]{d} \\
Z(s)\arrow[equal]{r} & \Quot_{\BP^3}(R,n) \arrow[hook]{r}{\iota} & \Quot_{\BP^3}(L_0,n)
\end{tikzcd}
\]
where the vertical, resp.~horizontal inclusions are open, resp.~closed. Since $\mathbf V$ has rank $n$ and $\Quot^{\circ}_{\BP^3}(L_0,n)$ is smooth and irreducible of dimension $5n$ by \Cref{prop:smooth-open-quot}, every irreducible component of the zero locus of $s^\circ$ has dimension at least $5n-n = 4n$. Hence at every closed point $[q]\in \Quot^{\circ}_{\BP^3}(R,n)$, representing a short exact sequence $E \into R \onto \OO_Z$, one has
\[ 
4n\leq \dim_{[q]}\Quot^{\circ}_{\BP^3}(R,n) \leq  \hom(E,\OO_Z)=\ext^2(E,E)=4n,
\]
where the two equalities come from \Cref{lem:exts}\ref{ext-ii} and \Cref{lem:dimext2} respectively. Thus $\Quot^{\circ}_{\BP^3}(R,n)$ is smooth of dimension $4n$.

To conclude the proof, we argue that $\Quot^{\circ}_{\BP^3}(R,n)$ is connected. The morphism 
\[ 
\begin{tikzcd}
\Quot^{\circ}_{\BP^3}(R,n)\arrow{r} & \Sym^n(\BP^3)\setminus\Delta^{\mathsf{big}} 
\end{tikzcd}
\]
is surjective and proper: indeed, it is surjective by \Cref{lemma:homR-O}\ref{surj}, and it is proper because it is a base change of the proper morphism $\Quot_{\BP^3}(R,n) \to \Sym^n(\BP^3)$. The fibre over a reduced cycle $Z$ is either $(\BP^1)^n$ (when $Z\cap\Sing R=\emptyset$) or $\BP^2\times(\BP^1)^{n-1}$ (when $\Sing R$ is in $Z$), cf.~also \Cref{rmk:twohoms}. In either case, the fibres are connected, thus $\Quot^{\circ}_{\BP^3}(R,n)$ is connected by \Cref{lemma:conn-irr}\ref{CON}.
\end{proof}

We now prove an irreducibility result that will turn out crucial in establishing the integrality of the moduli space $M_{-3}$.

\begin{theorem}
\label{cor:irreducible-quot}
For every reflexive sheaf $R$ parametrised by $N=M_1$, the schemes $\Quot_{\BP^3}(R,1)=\BP(R)$ and $\Quot_{\BP^3}(R,2)$ are irreducible, and moreover $\BP(R)$ is smooth of dimension 4. In particular, the moduli spaces $M_{-1}$ and $M_{-3}$ are irreducible. 
\end{theorem}

\begin{proof}
The final statement follows from the first thanks to \Cref{prop:pi-is-ZLT}\ref{irr-iff-irr}. The statements regarding $\BP(R)$ are proved as follows: assuming $R$ is singular at $(1:0:0:0)$, the induced resolution
\[
\begin{tikzcd}
0 \arrow{r} & \OO_{\BP^3}(-1) \arrow{r}{\left(\begin{smallmatrix}
    x_1 \\ x_2 \\ x_3\end{smallmatrix}\right)} & \OO_{\BP^3}^{\oplus 3} \arrow{r} & R(1) \arrow{r} & 0
\end{tikzcd}
\]
yields a closed immersion 
\[
\begin{tikzcd}
\BP(R) = \BP(R(1))\arrow[hook]{r} & \BP^3 \times \BP^2,
\end{tikzcd}
\]
whose scheme-theoretic image is the hypersurface cut out by $x_1y_1+x_2y_2+x_3y_3=0$, where $(y_1:y_2:y_3)$ are homogeneous coordinates on $\BP^2$. This equation defines a smooth irreducible projective 4-fold. Smoothness can be checked directly via the Jacobian criterion.

Let us move to the irreducibility of $\Quot_{\BP^3}(R,2)$.
The claim is equivalent to the assertion that every quotient $[R \onto T]$ is a limit of quotients belonging to $\Quot^\circ_{\BP^3}(R,2)$. Let $p = \Sing R$. If $p \notin \Supp T$, there is nothing to prove, as $R$ is locally free around $\Supp T$ and Quot schemes of $n \leq 10$ points are irreducible in the locally free case \cite{Henni_Quot, Sivic}. We are left with analysing the cases $T = \kappa(p)^{\oplus 2}$ and $T = \kappa(p)[y]/y^2$, as the case $T=\OO_p \oplus \OO_q$ with $q\neq p$ already belongs to the open locus $\Quot^\circ_{\BP^3}(R,2)$. In both cases, we can argue locally, replacing $\BP^3$ with $\BA^3$ and assuming the reflexive module $R$ to be singular at the origin $0 \in \BA^3$. By the resolution \eqref{eqn:resolution-F}, the module $M\defeq R|_{\BA^3}$ is of the form
\begin{equation}
\label{M-presentation}
M = \frac{Ae_1\oplus Ae_2 \oplus Ae_3}{x_1e_1+x_2e_2+x_3e_3},
\end{equation}
where $A = \BC[x_1,x_2,x_3]$ is the coordinate ring of the local chart $\BA^3 \subset \BP^3$. Set $B=\BC[t]$. Given a surjection
\begin{equation}
\label{central-quotient}
\begin{tikzcd}
    M \arrow[two heads]{r}{q_0} & T,
\end{tikzcd}
\end{equation}
where $T$ is a finite $A$-module of length 2 supported at the origin, we want to find a $B$-flat family of quotients
\begin{equation}
\label{family-quotient}
\begin{tikzcd}
    M_B = \displaystyle\frac{B[x_1,x_2,x_3]e_1\oplus B[x_1,x_2,x_3]e_2 \oplus B[x_1,x_2,x_3]e_3}{x_1e_1+x_2e_2+x_3e_3} \arrow[two heads]{r}{q_B} & T_B
\end{tikzcd}
\end{equation}
whose restriction to $\Spec \kappa(0) \to \BA^1 = \Spec B$ is $q_0$ and such that $q_t\defeq q_B \otimes \kappa(t)$ defines a point of $\Quot^\circ_{\BP^3}(R,2)$ for all $t \in \BG_m\subsetneq \BA^1$. It is time to distinguish two cases, according to the type of $T$. 

\textbf{Case 1: $T = \kappa(0)^{\oplus 2}$}. Assume given a quotient \eqref{central-quotient} where $T = \kappa(0)^{\oplus 2}$ has generators $z_1$ and $z_2$ and $q_0$ sends $e_i$ to $a_iz_1+b_iz_2 \in T$ for $i=1,2,3$. To ensure the surjectivity of $q_0$, we assume the vectors $\boldit{a} = (a_1,a_2,a_3)$ and $\boldit{b} = (b_1,b_2,b_3)$ are linearly independent. Note that since $\Fm = (x_1,x_2,x_3)$ kills $T$, the relation defining the quotient \eqref{M-presentation}, namely $x_1e_1+x_2e_2+x_3e_3 = 0$, is automatically satisfied in $T$. Now, fix two distinct vectors $\boldit{\alpha},\boldit{\beta} \in \BC^3$ such that $\boldit{\alpha}\cdot \boldit{a} = 0 = \boldit{\beta}\cdot \boldit{b}$. Consider the closed subschemes
\[
\begin{tikzcd}
    Z_{\boldit{\alpha}}\arrow[hook]{r} & \BA^3 \times \BA^1,& Z_{\boldit{\beta}}\arrow[hook]{r} & \BA^3 \times \BA^1 
\end{tikzcd}
\]
having coordinate rings
\[
T_{\boldit{\alpha}} = \frac{\BC[x_1,x_2,x_3,t]}{(x_1-t\alpha_1,x_2-t\alpha_2,x_3-t\alpha_3)},\qquad
T_{\boldit{\beta}} = \frac{\BC[x_1,x_2,x_3,t]}{(x_1-t\beta_1,x_2-t\beta_2,x_3-t\beta_3)}.
\]
By our assumption $\boldit{\alpha}\neq\boldit{\beta}$, for any fixed $t \in \BG_m$ the sum $T_{\boldit{\alpha}}|_t \oplus T_{\boldit{\beta}}|_t$ is the structure sheaf of a reduced subscheme of length $2$ in $\BA^3$, given by $\Spec \kappa(t\boldit{\alpha})\amalg \Spec \kappa(t\boldit{\beta})$. Define $T_B \defeq T_{\boldit{\alpha}} \oplus T_{\boldit{\beta}}$ and note that this $B$-module is flat, being a direct sum of $B$-flat modules, cf.~\cite[\href{https://stacks.math.columbia.edu/tag/05NG}{Tag 05NG}]{stacks-project}. Next, construct the map $q_B$ as in \eqref{family-quotient} by sending $e_i\mapsto (a_i,b_i)$ for $i=1,2,3$. This is well-defined since the expression $x_1e_1+x_2e_2+x_3e_3$ gets sent to \begin{align*}
    x_1(a_1,b_1) + x_2(a_2,b_2) + x_3(a_3,b_3) &= t(\alpha_1a_1,\beta_1b_1) + t(\alpha_2a_2,\beta_2b_2) + t(\alpha_3a_3,\beta_3b_3) \\
    &= t(\boldit{\alpha}\cdot \boldit{a},\boldit{\beta}\cdot \boldit{b}) \\
    &= 0.
\end{align*}  
Finally, the restriction $q_B \otimes \kappa(0)$ is precisely the surjection $q_0$ and thus the map $q_B$ is surjective in an open neighbourhood of $0 \in \BA^1=\Spec B$. Thus $q_0$ is a limit of points $q_t$ in the open locus $\Quot_{\BP^3}^{\circ}(R,2)$.

\smallbreak
\textbf{Case 2: $T = \kappa(0)[y]/y^2$}. Let $\boldit{w} = (w_1,w_2,w_3) \in \BC^3 \setminus 0$ be the tangent direction of the curvilinear scheme defined by $T$. Then the $A$-module structure of $T$ is prescribed by $x_i\cdot 1 = w_iy$, and $x_i\cdot y = 0$. Fix a surjection $q_0\colon M \onto T$ as in \eqref{central-quotient}. Let $\boldit{a} = (a_1,a_2,a_3)$ and $\boldit{b} = (b_1,b_2,b_3)$ be the vectors in $\BC^3$ defining $q_0$, in the sense that $q_0(e_i) = a_i + yb_i$ for $i=1,2,3$. To achieve surjectivity, one needs to have $\boldit{a} \neq \mathbf{0}$. To have $q_0$ well-defined, on the other hand, one has to ensure that $x_1e_1+x_2e_2+x_3e_3$ goes to 0. But, by the $A$-module structure of $T$ recalled above, we have
\begin{align*}
q_0(x_1e_1+x_2e_2+x_3e_3) 
&=x_1(a_1+b_1y) + x_2(a_2+b_2y) + x_3(a_3+b_3y) \\
&=a_1x_1\cdot 1+a_2x_2\cdot 1+a_3x_3\cdot 1 \\
&= (a_1w_1+a_2w_2+a_3w_3)y.
\end{align*}
In other words, we must impose $\boldit{a}\cdot \boldit{w} = 0$. Now we make the following choices: we fix
\begin{itemize}
    \item [$\circ$] a vector $\boldit{u} \in \BC^3 \setminus \BC\cdot \boldit{w}$ such that $\boldit{u}\cdot \boldit{a} = 0$, and
    \item [$\circ$] a vector $\boldit{c} \in \BC^3$ such that $\boldit{u}\cdot \boldit{c} = 0$, and $\boldit{w}\cdot \boldit{c} = (\boldit{u}+\boldit{w})\cdot \boldit{b}$.
\end{itemize}
Note that the condition $\boldit{u}\notin \BC\cdot \boldit{w}$ makes it possible to find $\boldit{c}$ and implies $\boldit{u}+\boldit{w} \neq \mathbf{0}$.

Next, we construct a $B[x_1,x_2,x_3]$-algebra $C$ as follows: we set $C\defeq B[s] / (s^2-ts)$ and we equip it with the surjective $B$-algebra homomorphism \[
\begin{tikzcd}[row sep=tiny]
B[x_1,x_2,x_3] \arrow[two heads]{rr} && C     \\
x_i\arrow[mapsto]{rr} && tu_i+sw_i.
\end{tikzcd}
\]
This defines a closed subscheme
\[
\begin{tikzcd}
Z_{\boldit{u},\boldit{w}} = \Spec C\arrow[hook]{r} & \BA^3 \times \BA^1,
\end{tikzcd}
\]
flat over $B$, with the property that for each $t \in \BG_m$ the subscheme $Z_{\boldit{u},\boldit{w}}|_t \subset \BA^3$ consists of the (distinct) points $\Spec \kappa (t\boldit{u})$, corresponding to $s=0$ in the splitting $0 = s^2-ts = s(s-t)$, and $\Spec \kappa (t(\boldit{u}+\boldit{w}))$, corresponding to $s=t$ in the same splitting. By our choices above, all 3 vectors $\boldit{u}$, $\boldit{w}$ and $\boldit{u}+\boldit{w}$ are nonzero and distinct, so $Z_{\boldit{u},\boldit{w}}|_t$ is indeed reduced of length 2. On the other hand, $Z_{\boldit{u},\boldit{w}}|_0$ is clearly equal to $T$ (up to $s \mapsto y$).

Finally, construct the $B$-module homomorphism $q_B\colon M_B \to T_B \defeq C$, by sending $e_i \mapsto a_i-tc_i+sb_i$. Using vector notation, to check that this is well-defined, we must confirm that $\boldit{x}\cdot (\boldit{a}-t\boldit{c}+s\boldit{b}) = 0$ in $C$. Spelling this out, we find
\begin{align*}
\boldit{x}\cdot (\boldit{a}-t\boldit{c}+s\boldit{b}) 
&= (t\boldit{u}+s\boldit{w}) \cdot (\boldit{a}-t\boldit{c}+s\boldit{b})\\    
&= t\boldit{u}\cdot (\boldit{a}-t\boldit{c}+s\boldit{b})+s\boldit{w}\cdot \boldit{a} + st \boldit{w}\cdot(\boldit{b}-\boldit{c}) \\
&= st \left(\boldit{u}\cdot \boldit{b} + \boldit{w}\cdot(\boldit{b}-\boldit{c})\right) \\
&= 0,
\end{align*}
where we have used the relations $\boldit{u}\cdot \boldit{a} = \boldit{w}\cdot \boldit{a} = \boldit{u}\cdot \boldit{c} = 0$ for the third identity and $\boldit{w}\cdot \boldit{c} = (\boldit{u}+\boldit{w})\cdot \boldit{b}$ for the last one. Once again $q_B \otimes \kappa(0)$ agrees with $q_0$ and therefore $q_B$ is surjective in an open neighbourhood of $0 \in \BA^1$. By definition of $T_B$, nearby quotients belong to $\Quot^\circ_{\BP^3}(R,2)$, as desired.
\end{proof}

Thanks to \Cref{cor:irreducible-quot}, we can obtain the following useful statement on $\Quot_{\BP^3}(R,2)$.

\begin{prop}
\label{prop:int+CM+8}
For every reflexive sheaf $R$ parametrised by $N=M_1$, the scheme $\Quot_{\BP^3}(R,2)$ is integral and Cohen--Macaulay of dimension 8.    
\end{prop}

\begin{proof}
We have that $\Quot_{\BP^3}(R,2)$ is irreducible by \Cref{cor:irreducible-quot}, and thus generically smooth of dimension 8 by \Cref{prop:Quot-generic-smoothness}. On the other hand, $\Quot_{\BP^3}(R,2)$ is embedded in the irreducible 10-fold $\Quot_{\BP^3}(L_0,2)$  as the zero locus of a section $s^{[2]}$ of a rank 2 bundle. This section is regular, since its zero locus has the correct codimension and $\Quot_{\BP^3}(L_0,2)$ is Cohen--Macaulay by \Cref{prop:quot-red-CM}. In particular, $\Quot_{\BP^3}(R,2)$ is also Cohen--Macaulay by \cite[Thm.~14.4(c)]{Ful}. Finally since $\Quot_{\BP^3}(R,2)$ is generically smooth and Cohen--Macaulay, it is in fact reduced (e.g.~because it has properties $R_0$ and $S_1$, cf.~\cite[\href{https://stacks.math.columbia.edu/tag/033P}{Tag 033P}]{stacks-project}), and thus integral. 
\end{proof}

\subsection{Main component}
\label{sec:maincomponent}
Fix $c_3 \in \BZ_{\leq -1}$ odd, and set $n = \frac{1-c_3}{2}$. In $M_{c_3}$, consider the open subscheme
\[
M_{c_3}^{\circ} \subset M_{c_3}
\]
defined as the preimage of $\Sym^n(\BP^3)\setminus\Delta^{\mathsf{big}}$ along the composition
\begin{equation} \label{eq:Phipr2}
\begin{tikzcd}
M_{c_3}\arrow{r}{\Phi} & N\times\Sym^n(\BP^3)\arrow{r}{\pr_2} & \Sym^n(\BP^3),
\end{tikzcd}  
\end{equation}
where $\Phi$ is the wall-crossing map (after the identification of \Cref{thm:iso-with-sim}). Explicitly, it is given by 
\[
M_{c_3}^{\circ} = \Set{[E] \in M_{c_3} | 
\begin{array}{cc}
E^{\ast\ast} / E = \OO_Z\mbox{ with}    &  \\
 Z\subset \BP^3 \mbox{ reduced}
\end{array}
}
\subset M_{c_3}.
\]
Consider, similarly, the open subscheme
\[
\Quot_{p_N}^{\circ}(\SR,n) \defeq f_{c_3}^{-1}(M_{c_3}^{\circ}) \subset \Quot_{p_N}(\SR,n).
\]
Clearly, $f_{c_3}$ gives an isomorphism of schemes
\begin{equation}
\label{iso-circ}
\begin{tikzcd}
\Quot_{p_N}^{\circ}(\SR,n) \arrow{r}{\sim} & M_{c_3}^{\circ}.
\end{tikzcd}
\end{equation}

\begin{lemma}
\label{lemma:pi-circ-smooth}
The restriction
\[
\begin{tikzcd}
\Quot_{p_N}^{\circ}(\SR,n) \arrow{r}{\pi^\circ} & N
\end{tikzcd}
\]
of the structure morphism $\pi\colon \Quot_{p_N}(\SR,n) \to N$ is smooth of relative dimension $4n$.
\end{lemma}

\begin{proof}
Combine \Cref{prop:pi-is-ZLT} and \Cref{prop:Quot-generic-smoothness} with one another.
\end{proof}

\begin{theorem}
\label{mainComponent}
Fix $c_3 \in \BZ_{\leq -1}$ odd, and set $n = \frac{1-c_3}{2}$. Then,
\begin{itemize}
\item [\mylabel{smooth-open}{\normalfont{(i)}}] $M_{c_3}^{\circ}$ is a smooth irreducible open subset of $M_{c_3}$ of dimension $3+4n=5-2c_3$, and 
\item [\mylabel{irr-comp}{\normalfont{(ii)}}] the closure $M_{c_3}^{\main}\defeq \overline{M_{c_3}^\circ}$, equipped with its reduced induced structure, is an irreducible component of $M_{c_3}$ of dimension $5-2c_3$.
\end{itemize}
\end{theorem}

\begin{proof}
By \Cref{lemma:pi-circ-smooth} and the isomorphism \eqref{iso-circ}, the scheme $M_{c_3}^\circ$ is smooth of dimension $3+4n=5-2c_3$. We now prove it is connected. By \Cref{lemma:homR-O}\ref{surj}, the restricted map 
\[
\begin{tikzcd}
\Phi\big|_{M_{c_3}^{\circ}} ~:~ M_{c_3}^{\circ}\arrow{r} &  N\times(\Sym^n(\BP^3)\setminus\Delta^{\mathsf{big}})
\end{tikzcd}
\]
is surjective; it is also proper, being a base change of the proper morphism $M_{c_3} \to N\times\Sym^n(\BP^3)$, and the fibre over a point $(R,Z)\in N\times(\Sym^n(\BP^3)\setminus\Delta^{\mathsf{big}})$ is either isomorphic to $(\BP^1)^n$ (when $\Sing R\cap Z = \emptyset$) or to $\BP^2 \times (\BP^1)^{n-1}$ (when $\Sing R\cap Z \neq \emptyset$), cf.~also \Cref{rmk:twohoms}. In both cases, all fibres are nonempty and connected. Since the target is connected, $M_{c_3}^\circ$ is connected by \Cref{lemma:conn-irr}\ref{CON}. Thus Item \ref{smooth-open} is proved. Item \ref{irr-comp} is a direct consequence.
\end{proof}

By analogy with the case of the Hilbert scheme of points, we call $M_{c_3}^{\main}$ the \emph{smoothable component}. We now give a lower bound for the dimension of the singular locus of $M_{c_3}$, exploiting the smoothable component.

\begin{prop}
\label{Sing-main-Component}
Fix $c_3 \in \BZ_{\leq -3}$ odd, and set $n = \frac{1-c_3}{2}$. Then
\[
\dim \Sing M_{c_3} \geq \dim M_{c_3}^{\main} - 5 = 4n-2.
\]
\end{prop}

\begin{proof}
We will show that there exists a closed subset $\Sigma \subset M_{c_3}^{\main} \cap \Sing M_{c_3}$ of codimension $5$ in $M_{c_3}^{\main}$.
Write $n = 2+k$ with $k \geq 0$.
Choose $[R] \in N$ and $k+1$ pairwise distinct points $x_1,y_1,\ldots,y_k$, none of which coincides with $\Sing R$. Construct the 0-dimensional sheaf
\[
Q = \OO_{x_1}^{\oplus 2}\oplus \OO_{y_1}\oplus \cdots\oplus \OO_{y_k}.
\]
Note that $\Aut(Q)$ has dimension $4+k$, whereas $\Hom(R,Q)$ has dimension $4+2k$, as $R$ is locally free around $\Supp Q$. Define $\Sigma\subset M_{c_3}$ to be the closure of the locus parametrising isomorphism classes of quotients $R \onto Q$ of this form. Using that the automorphism action is free on the surjective locus, taking into account the choice of $R$ (accounting for 3 parameters), the choice of $Q$ (accounting for $3(k+1)$ parameters), we see that $\Sigma$ has dimension $4k+6=4n-2$. Note that all such quotients lie in $M_{c_3}^{\main}$, as $\OO_{x_1}^{\oplus2}$ is smoothable by \Cref{cor:irreducible-quot}, while keeping the reduced summands $\OO_{y_j}$ fixed. So we have $\Sigma\subset M_{c_3}^{\main}$. The codimension of $\Sigma \subset M_{c_3}^{\main}$ is then easily computed to be precisely $5$, using $\dim M_{c_3}^{\main} = 4n+3$. It remains to show that $\Sigma \subset \Sing M_{c_3}$. Let $E\subset R$ denote the kernel of one such quotient. Then
\[
\Ext^2(E,E) = \Ext^2(\Fm_{x_1},\Fm_{x_1})^{\oplus 4}
\oplus \bigoplus_{j=1}^k \Ext^2(\Fm_{y_j}\oplus \OO_{\BP^3,y_j},\Fm_{y_j}\oplus \OO_{\BP^3,y_j}),
\]
which has dimension $12 + 4k = 4n+4$ by the same calculation of \Cref{lem:dimext2}. Thus $\ext^1(E,E) = 4n+7 > 4n+3 = \dim M_{c_3}^{\main}$, which implies $\Sigma \subset \Sing M_{c_3}$.
\end{proof}

\begin{remark}
\label{Tcomponents}
The smoothable component $M_{c_3}^{\main}$ coincides with the component $\mathrm{T}(-1,1,1,n)$, described in \cite[Main Theorem 1]{AJT2022}. Indeed, the sheaves $[E] \in M_{c_3}^{\circ}$ are precisely those obtained as kernels of epimorphisms $E^{**} \onto \OO_Z$, where $E^{**}$ is a reflexive sheaf with Chern classes $c_1=-1, c_2=c_3=1$, and where $Z$ is a $0$-dimensional reduced subscheme of $\mathbb{P}^3$, with length $n=\frac{1-c_3}{2}$. 
\end{remark}

\subsection{Rationality of the main component}
\label{sec:rationality-main-component}
The goal of this subsection is to prove that $M^{\main}_{c_3}$ is rational. The key idea is to use the fact that
\[
\begin{tikzcd}
\Quot_{p_N}(\SR,n) \arrow{r}{\pi} & N
\end{tikzcd}
\]
is a Zariski locally trivial fibration, cf.~\Cref{prop:pi-is-ZLT}, together with the rationality of the base and of the smoothable component of the fibre.

First of all, we recall some basic facts about the smoothable component 
\[
\Quot^{\main}_{\mathbb{A}^d}(\OO^{\oplus m},n) \subset \Quot_{\mathbb{A}^d}(\OO^{\oplus m},n)
\]
first introduced in \Cref{def:QUOT-circ}.

\begin{prop}
\label{smoothable}
The smoothable component $\Quot^{\main}_{\mathbb{A}^d}(\OO^{\oplus m},n)$ is rational.
\end{prop}
\begin{proof}
Consider the projection map
\[
\begin{tikzcd}
\rho\colon \Sym^n(\mathbb{A}^d\times \mathbb{P}^{m-1})\arrow{r} & \Sym^n(\mathbb{A}^d).
\end{tikzcd}
\]
The open subset 
\[
V_{d,m,n}\defeq \rho^{-1}(\Sym^n(\mathbb{A}^d)\setminus\Delta^{\mathsf{big}}) \subset \Sym^n(\mathbb{A}^d\times \mathbb{P}^{m-1})
\]
is isomorphic to the open subset 
\[
\Quot^\circ_{\mathbb{A}^d}(\OO^{\oplus m},n) \subset \Quot^{\main}_{\mathbb{A}^d}(\OO^{\oplus m},n)
\]
via the map sending a quotient
\[
\begin{tikzcd}
    {\bigg[}\OO^{\oplus m} \arrow[two heads]{r}{\varphi} & \displaystyle \bigoplus_{1\leq i\leq n} \kappa(p_i)   {\,\bigg]}\,\in\,\Quot^\circ_{\mathbb{A}^d}(\OO^{\oplus m},n)
\end{tikzcd}
\]
to the 0-cycle
\[
\sum_{1\leq i\leq n} (p_i,[\varphi_i])\,\in\,V_{d,m,n},
\]
where $[\varphi_i] \in \BP^{m-1}$ is the class of the map $\varphi_i \defeq \varphi|_{p_i} \colon \kappa(p_i)^{\oplus m} \onto \kappa(p_i)$. 

In particular, $\Quot^{\main}_{\mathbb{A}^d}(\OO^{\oplus m},n)$ is birational to $\Sym^{n}(\mathbb{A}^d\times\mathbb{P}^{m-1})$. Since $\mathbb{A}^d\times\mathbb{P}^{m-1}$ is birational to $\mathbb{A}^{d+m-1}$, it follows that $\Sym^{n}(\mathbb{A}^d\times\mathbb{P}^{m-1})$ is birational to $\Sym^{n}(\mathbb{A}^{d+m-1})$, which is well-known to be rational.
\end{proof}

We now turn to the proof that the smoothable component $M^{\main}_{c_3}$ is rational. In what follows, we implicitly use the natural identification $M_{c_3}=\Quot_{p_N}(\mathscr{R},n)$, for $n=\frac{1-c_3}{2}$.

\begin{prop}\label{fibration}
    The smoothable component $M^{\main}_{c_3}$ is birational to $\Quot_{\mathbb{A}^3}^{\main}(\OO^{\oplus 2},n)\times \mathbb{A}^3$, where $n=\frac{1-c_3}{2}$.
\end{prop}
\begin{proof}
    Let $R_0$ be a stable reflexive sheaf, corresponding to a point of $N$. By Proposition~\ref{prop:pi-is-ZLT}, the morphism $\pi\colon\Quot_{p_{N}}(\mathscr{R},n)\rightarrow N$ is Zariski locally trivial; in particular $R_0$ admits an open neighbourhood $U\subset N$ such that $\pi^{-1}(U)\cong U\times \Quot_{\mathbb{P}^3}(R_0,n)$. Denote by $p_0\in\mathbb{P}^3$ the singularity of $R_0$ and let $H\subset \mathbb{P}^3$ be a hyperplane containing $p_0$. Restricting to the complement of $H$, gives an open immersion
\[
\begin{tikzcd}
\Quot_{\mathbb{A}^{3}}(\OO^{\oplus 2},n)\arrow[hook]{r} & \Quot_{\BP^3}(R_0,n).
\end{tikzcd}
\]
This allows one to refine $\pi^{-1}(U)$, by taking $W\defeq U\times \Quot_{\mathbb{A}^{3}}(\OO^{\oplus 2},n)\subset U\times \Quot_{\mathbb{P}^3}(R_0,n)=\pi^{-1}(U)$. We can further refine $W$, by taking $W'\defeq U\times \Quot^{\circ}_{\mathbb{A}^{3}}(\OO^{\oplus 2},n)$.
    
    Note that $M^\circ_{c_3}$ and $W'$ are irreducible open subsets of $M_{c_3}=\Quot_{p_N}(\mathscr{R},n)$ with nonempty intersection: they both contain a point of the form $[R_0\onto\oplus_{i=1}^n\kappa(p_i)]$ where $p_1,\dots, p_n$ are pairwise distinct outside of $H$. This is enough to conclude.
\end{proof}

\begin{prop}
\label{prop:main-is-rational}
  The smoothable component $M^{\main}_{c_3}$ is rational. 
\end{prop}
\begin{proof}
Combine \Cref{fibration} and \Cref{smoothable}.
\end{proof}

\section{\texorpdfstring{Low values of $n$: proof of \Cref{thm:MAIN-small-n}}{}}

In this section we explore in detail the geometry of $M_{-1}$ and $M_{-3}$, proving \Cref{thm:MAIN-small-n}. We already know that both spaces are irreducible, by \Cref{cor:irreducible-quot}.

\subsection{\texorpdfstring{The geometry of $M_{-1}$}{}}
\label{sec:M-1}

The fact that $M_{-1}$ is a smooth variety is now easily obtained, as a consequence of the irreducibility of the moduli space, combined with the generic smoothness of the smoothable component and the constancy of the tangent space dimension.

\begin{theorem}
\label{thm:smooth7}
The moduli space $M_{-1}$ is a smooth, rational projective variety of dimension $7$.
\end{theorem}

\begin{proof}
The moduli space $M_{-1}$ is irreducible by \Cref{cor:irreducible-quot}. Therefore $M_{-1} = M_{-1}^{\main}$, and in particular it is rational (by \Cref{prop:main-is-rational}) and generically smooth of dimension 7 by \Cref{mainComponent}\ref{smooth-open}. On the other hand, combining \Cref{lem:dimext2} and \Cref{thm:tangent_space}, we see that the Zariski tangent space at all of its points is 7-dimensional.
\end{proof}

We conclude the proof of \Cref{thm:MAIN-small-n}\ref{M-SMALL1} by computing the motive of $M_{-1}$.

\begin{theorem}
\label{cor:motive-1}
Let $K_0(\Var_{\BC})$ be the Grothendieck ring of varieties. There is an identity
\[
[M_{-1}] = 1 + 3 \BL + 6 \BL^2 + 8 \BL^3 + 8 \BL^4 + 6 \BL^5 + 3 \BL^6 + \BL^7\,\in\,K_0(\Var_{\BC}).
\]
\end{theorem}

\begin{proof}
By \Cref{prop:pi-is-ZLT}, we know that $M_{-1} \to N$ is Zariski locally trivial, thus we have the identity $[M_{-1}] = [\BP^3][\Quot_{\BP^3}(R,1)]$, where $[R] \in N$ is any point in the base. The motive of the fibre can be computed as follows. Let $U\subset \Quot(R,1)$ be the open subscheme parametrising quotients $[R\onto Q]$ such that $\Supp Q\cap \Sing R=\emptyset$. Since $R$ is locally free of rank $2$ away from $\Sing R$ and $\Quot(R,1)=\BP(R)$, the subscheme $U$ is a Zariski locally trivial $\BP^1$-bundle on $\mathbb{P}^3\setminus \Sing R$. Hence
\begin{equation*}
    [U]=[\mathbb{P}^3\setminus \Sing R]\cdot[\mathbb{P}^1]=(\mathbb{L}^3+\mathbb{L}^2+\mathbb{L})(\mathbb{L}+1)=\mathbb{L}^4+2\mathbb{L}^3+2\mathbb{L}^2+\mathbb{L}.
\end{equation*}
Define $Z\defeq \Quot(R,1)\setminus U$. Then $Z$ is isomorphic to $\BP(R|_p)$, where $R|_p$ denotes the fibre of $R$ over $p$. As computed in~\Cref{rmk:twohoms}, we have $\dim_\mathbb{C}R|_p=3$, and therefore $[Z]=[\mathbb{P}^2]=\mathbb{L}^2+\mathbb{L}+1$. Consequently,
\begin{equation*}
    [\Quot(R,1)]=[U]+[Z]=\mathbb{L}^4+2\mathbb{L}^3+3\mathbb{L}^2+2\mathbb{L}+1.
\end{equation*}
Multiplying this by $[\BP^3]$ gives the formula.
\end{proof}

We have now completed the proof of \Cref{thm:MAIN-small-n}\ref{M-SMALL1}.

\begin{remark}
\label{rmk:motive-quot-1}
The conclusion of \Cref{cor:motive-1} could have been reached also as follows. One easily checks that the restriction of the DT wall-crossing map $\Phi\colon M_{-1} \to N \times \BP^3$ to the open subscheme $V \subset N\times \BP^3$, whose points are pairs $(R,q)$ such that $q \neq \Sing R$, defines a $\BP^1$-bundle $M_{-1}|_V \to V$. Under the identification $N \cong \BP^3$, we see that $V$ is just the complement of the diagonal. Thus 
\begin{align*}
[M_{-1}]
   &= [V][\BP^1] + [\BP^3] [\BP^2].
\end{align*}
Using that $V$ has class $[\BP^3]([\BP^3]-1)$, the formula in \Cref{cor:motive-1} follows immediately.
\end{remark}

\begin{remark}
\Cref{cor:motive-1} specialises to $e(M_{-1})=36$, in agreement with the formula
\begin{equation}
\label{eqn:chi-gf}
\begin{split}
\sum_{c_3}e(M_{c_3})q^{c_3} &= 4(q+q^{-1})\mathsf M(q^{-2})^8 \\
&= 4q + \frac{36}{q} + \frac{240}{q^3} + \frac{1296}{q^5} + \frac{6152}{q^7} + \frac{26440}{q^9} + \cdots
\end{split}
\end{equation}
which is a special case of \cite[Thm.~3.5]{Gholampour2017a}. Here, $\mathsf M(q) = \prod_{m\geq 1}(1-q^m)^{-m}$ is the MacMahon function, enumerating plane partitions of the integers, cf.~\cite[Ex.~3.7]{Gholampour2017a}.
\end{remark}

\subsection{\texorpdfstring{The geometry of $M_{-3}$}{}}
\label{sec:M-3}

In contrast with $M_{-1},$ the moduli space $M_{-3}$ is singular (cf.~\Cref{Sing-main-Component}), but nonetheless, we will show that it is integral, Cohen--Macaulay and normal (cf.~\Cref{thm:M-3-integral}), thus proving \Cref{thm:MAIN-small-n}\ref{M-SMALL2}. 

Consider a closed point $[E] \in M_{-3}$ and the associated short exact sequence
\begin{equation}
\label{seq:onceagain}
\begin{tikzcd}
0 \arrow{r} & E \arrow{r} & R \arrow{r} & Q \arrow{r} & 0,   
\end{tikzcd}
\end{equation}
where $R = E^{**}$ is parametrised by $N=M_1$, and $Q$ is a 0-dimensional sheaf of length 2 on $\mathbb{P}^3$. In particular, $Q$ is either a structure sheaf $\OO_Z$ for some length 2 subscheme $Z \subset \BP^3$, possibly nonreduced, or of the form $\OO_q \oplus \OO_q$, for some $q \in \BP^3$. Therefore, the possible quotients $Q=R/E$ fall in exactly one of the following cases as $[E]$ varies in $M_{-3}$:
\begin{enumerate}\label{casesQ}
\item [\mylabel{case1}{(1)}] $Q=\OO_{q_1}\oplus\OO_{q_2}$ with $q_1\ne q_2$ and neither coincides with the singularity of $R$;
\item [\mylabel{case2}{(2)}] $Q=\OO_{p}\oplus\OO_{q}$ where $p\ne q$ and $p$ is the singularity of $R$;
\item [\mylabel{case3}{(3)}] $Q=\OO_{Z}$ where $Z$ is supported at a single point away from the singularity of $R$;
\item [\mylabel{case4}{(4)}] $Q=\OO_{Z}$ where $Z$ is supported at the singularity of $R$;
\item [\mylabel{case5}{(5)}] $Q=\OO_{q}^{\oplus2}$ and $q$ does not coincide with the singularity of $R$;
\item [\mylabel{case6}{(6)}] $Q=\OO_{p}^{\oplus2}$ and $p$ is the singularity of $R$.
\end{enumerate}

\begin{prop}
\label{prop-extstrata}
Fix a closed point $[E] \in M_{-3}$ with $Q = E^{**}/E$, as in the cases described in {\normalfont{\ref{case1}}}--{\normalfont{\ref{case6}}}. Then the dimensions of $\mathrm{Ext}^i(E,E)$ for $i=1,2$ are given in \Cref{table:-3}. 

In particular, only the cases {\normalfont{\ref{case5}}} and {\normalfont{\ref{case6}}} yield singularities. 
\end{prop}
\begin{table}[ht]
\centering
\renewcommand{\arraystretch}{1.1}
\begin{tabular}{@{}c@{\qquad}c@{\qquad}c@{\qquad}c@{}}
\textbf{Case}
& $\boldsymbol{\ext^1(E,E)}$
& $\boldsymbol{\ext^2(E,E)}$
& \textbf{Smooth point?} \\
\midrule[1.2pt]
\ref{case1}, \ref{case2}, \ref{case3}, \ref{case4}
& 11 & 8 & yes \\
\midrule[0.4pt]
\ref{case5}, \ref{case6}
& 15 & 12 & no \\
\bottomrule[0.6pt]
\end{tabular}
\caption{Tangent and obstruction space dimensions for each $[E] \in M_{-3}$.}
\label{table:-3}
\end{table}

\begin{proof}
Cases \ref{case1} and \ref{case2} involve reduced supports, so they are covered by \Cref{lem:dimext2}, giving $\ext^2(E,E)=8$, and thus $\ext^1(E,E)=11$ by \Cref{thm:tangent_space}. As for case \ref{case5}, since the question is local and $R$ is free at $q$, the defining sequence \eqref{seq:onceagain} is equivalent to its localisation at $q$, namely
\[
\begin{tikzcd}
0 \arrow{r} & E_q \arrow{r} & \OO_{\BP^3,q}^{\oplus 2} \arrow{r} & \kappa(q)^{\oplus 2} \arrow{r} & 0,
\end{tikzcd}
\]
where $E_q \cong \Fm_q^{\oplus 2}$. By \Cref{lem:exts}\ref{ext-ii}, we find
\[
\ext^2(E,E) =  \hom(E,Q) = \hom(\Fm_q^{\oplus 2}, \kappa(q)^{\oplus 2}) = 4\cdot \hom(\Fm_q, \kappa(q)) = 12,
\]
where the last identity uses $\Hom(\Fm_q,\kappa(q))\cong \BC^{\oplus 3}$, this Hom space being the Zariski tangent space of $\BP^3$ at $q.$ Thus $\ext^1(E,E) = 15$ by \Cref{thm:tangent_space}. 

As for \ref{case3}, \ref{case4} and \ref{case6}, we use again $\ext^2(E,E) =  \hom(E,Q)$ and we compute the latter via the sequence
\[
\begin{tikzcd}
0 \arrow{r} & \Hom(Q,Q) \arrow{r} & \Hom(R,Q) \arrow{r} & \Hom(E,Q) \arrow{r} & \Ext^1(Q,Q) \arrow{r}{w} & \Ext^1(R,Q),
\end{tikzcd}
\]
yielding the relation
\begin{equation}
\label{alt-sum}
\ext^2(E,E) = \hom(R,Q)-\hom(Q,Q) + \ext^1(Q,Q) - \dim \im w.
\end{equation}
In both cases \ref{case3} and \ref{case4}, $Q = \OO_Z$ is a structure sheaf, so the sequence
\[
\begin{tikzcd}
0 \arrow{r} & \Hom(\OO_Z,\OO_Z) \arrow{r}{\sim} & \Hom(\OO_{\BP^3},\OO_Z) \arrow{r}{0} & \Hom(\SI_Z,\OO_Z) \arrow{r}{\sim} & \Ext^1(\OO_Z,\OO_Z) \arrow{r} &  0 
\end{tikzcd}
\]
shows that $\hom(\OO_Z,\OO_Z)=2$ and $\ext^1(\OO_Z,\OO_Z)=6$, since $\Hom(\SI_Z,\OO_Z)$ is the tangent space of the smooth 6-fold $\Hilb^2(\BP^3)$ at $[Z]$. Now, in case \ref{case3}, $p \notin Z$, so $\hom(R,\OO_Z)=4$ by \Cref{lemma:homR-O}\ref{hom(R,Q)}. This shows that $\ext^1(R,\OO_Z)=0$, therefore \Cref{alt-sum} yields $\ext^2(E,E)=4-2+6-0=8$. In case \ref{case4}, the middle map in the exact sequence
\[
\begin{tikzcd}[column sep=small]
0 \arrow{r} & \Hom(R,\OO_Z) \arrow{r} & \Hom(\OO_{\BP^3}(-1),\OO_Z)^{\oplus 3} \arrow{r} & \Hom(\OO_{\BP^3}(-2),\OO_Z) \arrow{r} & \Ext^1(R,\OO_Z) \arrow{r} & 0,
\end{tikzcd}
\]
which can be identified with a map $\HH^0(\BP^3,\OO_Z)^{\oplus 3} \to \HH^0(\BP^3,\OO_Z)$,
has rank 1, thus $\ext^1(R,\OO_Z)=1$ and then $\hom(R,\OO_Z) = 5$, so that \Cref{alt-sum} yields $\ext^2(E,E)=5-2+6-1=8$, since $w \neq 0$ as one can easily check adapting the proof of \Cref{lem:dimext2}.

Finally, in case \ref{case6}, by \Cref{lemma:homR-O}\ref{red-pt} we have $\hom(R,\OO_p^{\oplus 2})=6$, which forces $\ext^1(R,\OO_p^{\oplus 2})=2$. On the other hand, $\hom(\OO_p^{\oplus 2},\OO_p^{\oplus 2})=4$ and $\ext^1(\OO_p^{\oplus 2},\OO_p^{\oplus 2}) = 4\cdot \ext^1(\OO_p,\OO_p) = 12$.
The argument given in the final bit of the proof of \Cref{lem:dimext2}, again with minor modifications, shows that $w$ has maximal rank, namely 2, thus \Cref{alt-sum} yields $\ext^2(E,E)=6-4+12-2=12$. This concludes the ext calculations, confirming the first two columns of \Cref{table:-3}.

Finally, $M_{-3}$ is irreducible of dimension 11 thanks to \Cref{cor:irreducible-quot} and \Cref{mainComponent}\ref{smooth-open}. Therefore our $\ext^1$ calculations precisely detect the smoothness of all possible points $[E] \in M_{-3}$, confirming that only cases \ref{case5} and \ref{case6} give rise to singular points.
\end{proof}

\begin{theorem}
\label{thm:M-3-integral}
The moduli space $M_{-3}$ is an integral, 11-dimensional projective variety, which is moreover Cohen--Macaulay, regular in codimension 4, and normal.
\end{theorem}

\begin{proof}
We already know $M_{-3}$ is irreducible of dimension 11. To prove it is reduced, we argue as follows. By \Cref{prop:pi-is-ZLT} and \Cref{thm:iso_relative_quot_M}, the map $\pi\colon M_{-3} \to N$ is flat. Its target is reduced, being isomorphic to $\BP^3$. But by \Cref{prop:int+CM+8} the fibres of $\pi$ are integral (hence reduced) and Cohen--Macaulay of dimension 8. It follows that $M_{-3}$ is reduced (hence integral), and also Cohen--Macaulay: by \cite[\href{https://stacks.math.columbia.edu/tag/0C0W}{Tag 0C0W}]{stacks-project} applied with $Z=\Spec \BC$, it is enough to observe that $\pi\colon M_{-3} \to N$ is a Cohen--Macaulay morphism, i.e.~a flat morphism with Cohen--Macaulay fibres.

To say that $M_{-3}$ is regular in codimension 4 means that $\Sing M_{-3}$ has codimension at least $5$. By \Cref{prop-extstrata}, only the cases \ref{case5} and \ref{case6} contribute to the singular locus. The former gives the stratum corresponding to $(n,k) = (2,0)$ in \Cref{Sing-main-Component}, so has codimension 5. The latter, on the other hand, has codimension 6. Indeed, to give a quotient $R \onto \OO_{p}^{\oplus 2}$ with $p = \Sing R$ one needs 3 parameters for the choice of $[R]\in N$, 6 further parameters for the choice of $R \onto \OO_{p}^{\oplus 2}$ (because $R\otimes \kappa(p)$ is 3-dimensional). Subtracting $4 = \dim \Aut(\OO_p^{\oplus 2})$, we find $3+6-4 = 5$ moduli, whence codimension 6.

To prove $M_{-3}$ is normal, by \cite[\href{https://stacks.math.columbia.edu/tag/033P}{Tag 033P}]{stacks-project} we need to show it has properties $S_2$ and $R_1$. But we just showed $M_{-3}$ is regular in codimension 4, and it is Cohen--Macaulay (hence has $S_k$ for all $k$).
\end{proof}


\section{Brill--Noether locus}\label{sec:BN locus}
Fix an odd $c_3 \in \BZ_{\leq 1}$ and set $n = \frac{1-c_3}{2}$. We now consider the \emph{Brill--Noether locus}
\[
B_{c_3} \defeq \Set{[E] \in M_{c_3} \,|\, h^0(\BP^3,E(1))>0} \subset M_{c_3}.
\]
It is a closed subset, that we endow with the reduced induced subscheme structure. We will show its connectedness in \Cref{prop:bn-connected}. For $n=0,1,2$, the situation is clear, by the following remark.

\begin{remark}
\label{BN-small-n}
We have a set-theoretic identity $B_{c_3}=M_{c_3}$ when $c_3=1,-1,-3$ (i.e.~$n=0,1,2$). We already noticed in \Cref{sec:resolution} that $h^0(\BP^3,R(1))=3$ for every $[R]\in M_1$, so $B_1=M_1$. For $c_3=-1,-3$ fix $[E] \in M_{c_3}$ and set $R=E^{\ast\ast}$. Taking cohomology of the twisted short exact sequence
\[
\begin{tikzcd}
    0 \arrow{r} & E(1) \arrow{r} & R(1) \arrow{r} & Q \arrow{r} & 0
\end{tikzcd}
\]
gives
\[
\begin{tikzcd}
0\arrow{r} & \HH^0(\BP^3,E(1)) \arrow{r} & \HH^0(\BP^3,R(1)) = \BC^3 \arrow{r} & V \arrow{r} & 0,
\end{tikzcd}
\]
with $V\subseteq \HH^0(\BP^3,Q) = \BC^n$ for $1\leq n\leq 2$ and $V\neq 0$ since $R(1)$ is generated by global sections and surjects onto $Q$. This shows that $h^0(\BP^3,E(1)) \in \set{1,2}$, and in fact $h^0(\BP^3,E(1))=2$ when $n=1$.     
\end{remark}

For $c_3\leq-5$, $B_{c_3}$ is a \emph{proper and nonempty} closed subset of $M_{c_3}$, thanks to the following two remarks.

\begin{remark}
\label{BN-proper-locus}
For $c_3\leq -5$ (or equivalently $n\geq 3$), the complement of the Brill--Noether locus is nonempty. Indeed, pick homogeneous coordinates $x_0,x_1,x_2,x_3$ on $\BP^3$ and define $R(1)$ via the associated resolution
    \[
\begin{tikzcd}
0 \arrow{r} & \OO_{\BP^3}(-1) \arrow{r}{\left(\begin{smallmatrix}
    x_1 \\ x_2 \\ x_3\end{smallmatrix}\right)} & \OO_{\BP^3}^{\oplus 3} \arrow{r} & R(1) \arrow{r} & 0,
\end{tikzcd}
\]
so that $\Sing R = p \defeq (1:0:0:0)$. Let $p_1,p_2,p_3\in\BP^3$ be the other coordinate points, and choose distinct points $q_1,\dots,q_{n-3}$ in $\BP^3\setminus\set{p,p_1,p_2,p_3}$. Set
\[
Q = \OO_{p_1}\oplus \OO_{p_2}\oplus \OO_{p_3}\oplus \OO_{q_1}\oplus \cdots\oplus \OO_{q_{n-3}}.
\]
Since $R$ is locally free around $\Supp Q$, one can easily choose a surjective morphism $R(1)\onto Q$, such that the induced map
\begin{equation}
\label{eq:mapbetweenglobalsectioninjective}
\begin{tikzcd}
\HH^0(\BP^3,R(1))\arrow[hook]{r} & \HH^0(\BP^3,Q)    
\end{tikzcd}
\end{equation}
is injective. Consider $E(1)$ fitting in the exact sequence
\[
\begin{tikzcd}
    0 \arrow{r} & E(1) \arrow{r} & R(1) \arrow{r} & Q \arrow{r} & 0.
\end{tikzcd}
\]
It follows by the injectivity of the morphism in \Cref{eq:mapbetweenglobalsectioninjective} that $\HH^0(\BP^3,E(1))=0$.
\end{remark}

\begin{remark}
When $c_3\leq -5$, the locus $B_{c_3}\subsetneq M_{c_3}$ is nonempty. Indeed given a 1-dimensional subscheme $L\subset \BP^3$ defining a point $[L]\in\Hilb^{t+1+n}(\BP^3)$, with maximal purely 1-dimensional component $\ell\subset L$, one has
\begin{equation}
\label{eqn:Ext^1-line}
\Ext^1(\SI_L(1),\OO_{\BP^3})^\ast 
    \cong \HH^2(\BP^3,\SI_L(-3)) 
    \cong \HH^1(\BP^3,\OO_L(-3)) 
    \cong \HH^1(\BP^3,\OO_\ell(-3)) 
    \cong \HH^0(\BP^3,\OO_\ell(1))^\ast 
    \cong\BC^2.
\end{equation}
Now, it is easy to check that any rank 2 torsion-free sheaf $E$ such that $E(1)$ fits in a nontrivial extension of the form
\begin{equation}
\label{eq:oeil}
\begin{tikzcd}
0 \arrow{r} & \op3 \arrow{r}{\alpha} & E(1) \arrow{r} & \SI_L(1) \arrow{r} & 0
\end{tikzcd}
\end{equation}
is stable with $(c_1(E),c_2(E),c_3(E))=(-1,1,1-2n)$, i.e.~it defines a point of $M_{c_3}$. Moreover, one clearly has $h^0(\BP^3,E(1))>0$. In addition, we remark that if $L$ is reduced, then $[E]$ in fact belongs to $M_{c_3}^{\circ}\subset M_{c_3}$, showing that $B_{c_3}\cap M_{c_3}^{\circ}\ne\emptyset$. Indeed, if $L$ is reduced, then its ideal sheaf sits in a short exact sequence
\[
\begin{tikzcd}
0 \arrow{r} & \SI_L \arrow{r} & \SI_\ell \arrow{r} & \OO_Z \arrow{r} & 0 
\end{tikzcd}
\]
where $Z\subset \BP^3$ is a reduced 0-dimensional subscheme disjoint from the line $\ell$. Comparing with the middle row of the diagram \eqref{eq:diag1}, we conclude that $E^{**}/E=\OO_Z$, so that $[E] \in M_{c_3}^\circ$.
\end{remark}

The following result is based on the main conclusions of \cite{Jardim-Mu}. It follows from \cite[Prop.~2.19(1)]{Jardim-Mu}, applied with $k=1$, and \cite[Thm.~3.13]{Jardim-Mu}.

\begin{prop}
\label{prop:S-fine}
Fix an odd $c_3 \in \BZ_{\leq 1}$ and $0 < \delta < 1/2$. Let $P$ be the Hilbert polynomial corresponding to sheaves of the form $E(1)$, where $E$ is in class $(-1,1,c_3)$.\footnote{The polynomial $P$ is easily seen to be determined by the Chern character $\ch(E(1)) = (2,1,-1/2,(3c_3-2)/6)$.} Let $\CS_{c_3} = \CS_{\BP^3}(P,\delta)$ be the moduli space of $\delta$-semistable pairs $\alpha\colon \OO_{\BP^3} \to E(1)$ of type $P$. Then all such pairs are $\delta$-stable, so that in particular $\CS_{c_3}$ is a fine, projective moduli space.

Moreover, $\CS_{c_3}$ comes with two morphisms
\begin{equation} 
\label{eq:diagSMH}
\begin{tikzcd}[row sep=large,column sep=large]
\CS_{c_3}\arrow[swap]{d}{\Psi}\arrow{r}{\Gamma} & \Hilb^{t+1+n}(\BP^3)\\
M_{c_3} &
\end{tikzcd}
\end{equation}
defined on closed points by
\[
\Psi (E(1),\alpha) = E, \qquad
\Gamma (E(1),\alpha) = \coker(\alpha) \otimes \OO_{\BP^3}(-1).
\]
Moreover, $\Gamma$ is surjective and $\Psi$ restricts to a surjection $\Psi\colon \CS_{c_3,\red} \to B_{c_3}$. The respective fibres can be described as
\[
\Psi^{-1}(E) = \mathbb{P}\HH^0(\BP^3,E(1)), \qquad \Gamma^{-1}(L)=\mathbb{P}\Ext^1(\SI_L(1),\op3). 
\]
\end{prop}

\begin{remark}
When $c_3=1$, the diagram \eqref{eq:diagSMH} simplifies to the classical incidence correspondence
$\BP^3\leftarrow\mathbb{F}(0,1;3)\rightarrow\BG(1,3) = \Gr(2,4)$, where $\mathbb{F}(0,1;3)$ denotes the flag of points in lines in $\BP^3$, cf.~\cite[Sec.~6.1.3]{Jardim-Mu}.
\end{remark}

\begin{prop} 
\label{prop:bn-connected}
Fix an odd $c_3 \in \BZ_{\leq -1}$ and set $n = \frac{1-c_3}{2}$, so that $n \in \BZ_{\geq 1}$. Then the Brill--Noether locus $B_{c_3}$ is connected. In addition, if  $~\Hilb^{t+1+n}(\BP^3)$ is irreducible, then so is $B_{c_3}$.
\end{prop}
\begin{proof}
First, note that the moduli space of pairs $\CS_{c_3}$ is connected by \Cref{lemma:conn-irr}\ref{CON}, for it maps surjectively and properly, via $\Gamma$, to a connected scheme (see \cite{MR0213368} for connectedness of the Hilbert scheme of $\BP^3$), with connected fibres. Since the morphism $\Psi\colon \CS_{c_3,\red} \to B_{c_3}$ is surjective, its target must be connected as well.

Next, note again that $\Gamma$ is proper and surjective, and by \eqref{eqn:Ext^1-line} all of its fibres are isomorphic to $\BP^1$, in particular irreducible and of the same (finite) dimension. It follows from \Cref{lemma:conn-irr}\ref{IRR2} that $\CS_{c_3}$ is irreducible whenever $\Hilb^{t+1+n}(\BP^3)$ is, and thus the same occurs with $B_{c_3}$.
\end{proof}

The following corollary gives an alternative proof of the final assertion in \Cref{cor:irreducible-quot}.

\begin{corollary}
\label{irreducibleB}
The spaces $B_{-1}=M_{-1}$, $B_{-3}=M_{-3}$ and $B_{-5}$ are irreducible.
\end{corollary}

\begin{proof}
Chen and Nollet proved that $\Hilb^{t+1+n}(\BP^3)$ is irreducible when $n\leq 3$ \cite[Prop.~3.7]{Chen-Nollet}. Thus the result follows from \Cref{prop:bn-connected}.
\end{proof}

\begin{remark}
Unlike the case $n=1$, where $h^0(\BP^3,E(1))=2$ constantly on $[E] \in M_{-1}$, in the case $n=2$ both possible values for $h^0(\BP^3,E(1))$, namely 1 and 2, are attained. This can be checked via examples, but we offer an argument exploiting the DT/PT wall-crossing. By \Cref{eqn:chi-gf}, we know that $e(M_{-3}) = 240$. On the other hand, $e(\Hilb^{t+3}(\BP^3)) = 174$. To see this, one can use the unweighted rank 1 DT/PT correspondence
\begin{align*}
\sum_{n \geq 0}e(\Hilb^{t+1+n}(\BP^3))q^n 
&= \mathsf M(q)^{e(\BP^3)} \cdot \sum_{n\geq 0}e(\Sym^n(\CL/\Gr(2,4))) q^n \\
&= \mathsf M(q)^4 \cdot e(\Gr(2,4))\cdot \sum_{n\geq 0}\,(n+1)q^n \\
&= \mathsf M(q)^4\cdot \frac{6}{(1-q)^2},
\end{align*}
where $\CL\subset \BP^3 \times \Gr(2,4)$ is the universal line. Expanding the right-hand side and extracting the coefficient of $q^2$ gives 174. It follows that $e(\CS_{-3})=174\cdot 2 = 348$. Since this is not a multiple of $240$, the fibres of $\CS_{-3} \to M_{-3}$ are not isomorphic to each other.
\end{remark}

\appendix

\section{\texorpdfstring{Proof of \'etaleness of $u$}{}}
\label{app:Kuranishi}
In the proof of \Cref{thm:iso-with-sim}, we left behind a technical bit, namely the proof that the morphism $u$ arising from the commutative diagram
\begin{equation}
\label{diag:gmspaces}
\begin{tikzcd}
\mathcal{N} \times \CCoh^{n}(\BP^3)
\arrow{r}{\overline u}\arrow[swap]{d} & \FM\arrow{d} \\
N \times  \Sym^{n}(\BP^3)
\arrow{r}{u} & \BM
\end{tikzcd} \qquad
\begin{tikzcd}
(R,Q)\arrow[mapsto]{d}\arrow[mapsto]{r}{\overline u} & E=R\oplus Q{[}-1{]} \arrow[mapsto]{d} \\
(R,{[}Q{]})\arrow[mapsto]{r}{u} & {[}E{]}
\end{tikzcd}
\end{equation}
is \'etale. Our claim is that, for every closed point $s = (R,[Q]) \in N \times \Sym^n(\BP^3)$, the canonical map
\begin{equation}
\label{eqn:completed-local-map}
\begin{tikzcd}
\widehat{\OO}_{\BM,[E]} \arrow{r}{\widehat u_s^{\#}} & \widehat{\OO}_{N \times  \Sym^{n}(\BP^3),s}
\end{tikzcd}
\end{equation}
is an isomorphism. 

Throughout, fix a polystable representative of $s$, namely a closed point 
\begin{equation}
\label{eqn:polystable-rep}
(R,Q) \,\in\,\mathcal{N} \times \CCoh^{n}(\BP^3), \qquad Q = \bigoplus_{\alpha=1}^k \OO_{q_\alpha}\otimes W_\alpha, \qquad m_\alpha \defeq \dim_{\BC}W_\alpha,
\end{equation}
where $q_1,\ldots,q_k$ are the distinct points in the support of $Q$, so that $n = \sum_\alpha m_\alpha$.

The deformation theory of complexes \cite{Lieblich1}, together with the formal local structure of good moduli spaces \cite{Alper-local-properties}, provides explicit presentations of the completed local rings appearing in \eqref{eqn:completed-local-map}. Indeed, by \cite[Lemma 2.7]{Chen-Pertusi-Zhao}, the map $\widehat u_s^{\#}$ can be identified with a morphism
\begin{equation}
\label{eqn:completed-local-map2}
\begin{tikzcd}
 \left(\dfrac{\BC\llbracket \Ext^1(E,E)^\ast \rrbracket}{\BI_E}\right)^{\Aut(E)}\arrow{r} & \left(\dfrac{\BC\llbracket \Ext^1(R,R)^\ast \rrbracket}{\BI_R}\right)^{\Aut(R)}\widehat{\otimes}\left(\dfrac{\BC\llbracket \Ext^1(Q,Q)^\ast \rrbracket}{\BI_Q}\right)^{\Aut(Q)}.
\end{tikzcd}
\end{equation}
In \Cref{sec:kuranishiideal} we review some general theory: for each $F\in\{E,R,Q\}$, 
\[
\BI_F \subset \BC\llbracket \Ext^1(F,F)^\ast \rrbracket = \widehat{\Sym} \Ext^1(F,F)^\ast
\]
denotes the ideal of Kuranishi relations attached to an $\Aut(F)$-equivariant Kuranishi map; the displayed quotients are the induced presentations of the miniversal formal deformation algebras $\mathsf A_F$. The ideal $\BI_F$ is $\Aut(F)$-invariant and the natural action of $\Aut(F)$ on $\BC\llbracket \Ext^1(F,F)^\ast \rrbracket$ induces an action on $\mathsf A_F = \BC\llbracket \Ext^1(F,F)^\ast \rrbracket/\BI_F$, allowing one to take invariants. In order to prove that $\widehat u_s^{\#}$, i.e.~\eqref{eqn:completed-local-map2} is an isomorphism, our argument requires two main ingredients:
\begin{enumerate}
\item the computation of the \emph{invariant rings}, by identifying the $\Ext^1$-spaces with representation spaces of the associated Ext-quivers and applying the Le Bruyn--Procesi theorem, as explained in \Cref{sec:extquiver}, and
\item the computation of the \emph{Kuranishi relations} using the differential graded Lie algebras (dgla's for short) controlling the deformations of complexes, as explained in \Cref{sec:B-quadratic}.
\end{enumerate}

These ingredients are combined in \Cref{sec:B-conclusion}, where we show that $\widehat u_s^{\#}$ is an isomorphism, cf. \Cref{thm:B-completed-local-isomorphism}.

We believe the content of this appendix is of independent interest. Indeed, it generalises some of the Ext-quiver calculations performed in \cite{Ricolfi-Savvas}, where the summand $R$ was not present; the main difference there is the use of derived geometry, which we do not need here.

\subsection{The Kuranishi ideal attached to a dgla splitting}
\label{sec:kuranishiideal}

In this section we do a brief recap of Kuranishi theory in the language of dgla's, following \cite[Sec.~13.2]{Manetti-Book}.

\begin{definition}
\label{L-splitting}
Let $L$ be a dgla with finite dimensional cohomology. A \emph{splitting of $L$} is a choice of splittings
\[
Z^r = B^r\oplus \CH^r,\qquad L^r = Z^r\oplus \CK^r, \qquad r \in \BZ,
\]
where $Z^r = \ker (\dd_L^r)\subset L^r$ are the $r$-cocycles and $B^r = \im(\dd_L^{r-1})\subset L^r$ are the $r$-coboundaries.
\end{definition}

We denote by $q^r \colon Z^r \onto \HH^r(L)$ the natural surjection.
A splitting of $L$ induces, for each $r \in \BZ$,
\begin{enumerate}
\item [\mylabel{split-i}{(i)}] two isomorphisms 
\[
\begin{tikzcd}
q^r\big|_{\CH^r}\colon \CH^r \arrow{r}{\sim} & \HH^r(L),\qquad \dd^r_L\big|_{\CK^r}\colon \CK^r \arrow{r}{\sim} & B^{r+1},
\end{tikzcd}
\]
\item [\mylabel{split-ii}{(ii)}] a projection $H_r\colon L^r \onto \CH^r$ killing $B^r\oplus \CK^r$, and
\item [\mylabel{split-iii}{(iii)}] a linear map 
\[
\begin{tikzcd}
\delta_r \colon L^{r+1} \arrow[two heads]{r}{} & B^{r+1} \arrow{r}{\sim} & \CK^r \arrow[hook]{r} & L^r,
\end{tikzcd}
\]
where the middle isomorphism is the inverse of $\dd^r_L|_{\CK^r}$, and the first surjection is the projection killing $\CK^{r+1}\oplus \CH^{r+1}$.
\end{enumerate}
There is a natural transformation of formal neighbourhoods\footnote{See e.g.~\cite[Sec.~10.2]{Manetti-Book} for the definition of formal neighbourhood attached to a vector space.} $\phi \colon \widehat L^1 \to \widehat L^1$ defined on an object $(A,\Fm_A)$ of $\Art_{\BC}$ by the collection of compatible maps 
\[
\begin{tikzcd}[row sep=tiny]
L^1 \otimes \Fm_A \arrow{r}{\phi_A} & L^1 \otimes \Fm_A \\
x \arrow[mapsto]{r} & x+\frac{1}{2}\delta_1 [x,x].
\end{tikzcd}
\]
Since its linear part is the identity, there is a well-defined formal inverse of $\phi$, denoted by $\phi^{-1}$, see \cite[Lemma 13.2.2]{Manetti-Book}.
\begin{remark}
\label{rmk:expansion}
    For $A\in\Art_{\BC}$ and for $x\in L^1\otimes \mathfrak{m}_A$, one has the relation
\begin{equation*}
    \phi^{-1}(x)=x-\frac{1}{2}\delta_1[\phi^{-1}(x),\phi^{-1}(x)].
\end{equation*}
Solving recursively gives $\phi^{-1}(x)=\sum_{n>0}f_n(x)$ where
\begin{equation}
\label{eq:recursion}
    f_1(x)=x ,\quad
    f_n(x)=-\frac{1}{2}\sum_{i=1}^{n-1}\delta_1[f_i(x),f_{n-i}(x)] \,\in\,L^1\otimes \mathfrak{m}_A,
\end{equation}
for $n>1$. Since $\mathfrak{m}_A^{N_A}=0$ for $N_A\gg 0$, we have that $f_n(x)=0$ for $n\gg 0$, thus $\sum_{n>0}f_n(x)$ is a finite sum.
\end{remark}

The \emph{Kuranishi map} attached to a splitting as in \Cref{L-splitting} is the natural transformation $\kappa\colon \widehat \CH^1 \to \widehat \CH^2$ defined on $(A,\Fm_A)$ in $\Art_{\BC}$ by the map
\begin{equation}
\label{kuranishi-map}
\begin{tikzcd}[row sep=tiny]
\CH^1 \otimes \Fm_A \arrow{r}{\kappa_A} & \CH^2 \otimes \Fm_A \\
x \arrow[mapsto]{r} & \frac{1}{2}H_2([\phi^{-1}(x),\phi^{-1}(x)]).
\end{tikzcd}
\end{equation}
For $A\in\Art_{\mathbb{C}}$ and $x=\sum_{i=1}^d x_i\otimes a_i\in\mathcal{H}^1\otimes \mathfrak{m}_A$, we can expand $\kappa_A$ as follows:
\begin{equation}\label{eq:expand}
    \kappa_A\left(\sum_{i=1}^d x_i\otimes a_i\right)=\frac{1}{2}\left(\sum_{i,j}H_2[x_i,x_j]\otimes a_ia_j-\sum_{i,j,k}H_2[x_i,\delta_1[x_j,x_k]]\otimes a_ia_ja_k+\cdots\right).
\end{equation}
The \emph{Kuranishi functor} is the fibre of the Kuranishi map over the origin, defining a deformation functor 
\[
\begin{tikzcd}
\Art_{\BC}\arrow{r}{\Kur} & \Sets,
\qquad
(A,\mathfrak m_A)\arrow[mapsto]{r} &
\Set{
x\in \CH^1\otimes\mathfrak m_A
|
\kappa_A(x)=0}.    
\end{tikzcd}
\]
Let $e_1,\dots, e_d$ be a basis of $\mathcal{H}^1$, with dual basis $x_1,\dots, x_d$. Let $\epsilon_1,\dots,\epsilon_m$ be a basis of $\mathcal{H}^2$. For $A\in\Art_\BC$ and $t_1,\dots,t_d\in\mathfrak{m}_A$ one can write
\begin{equation}
\label{eqn:explicit-Kur}
\kappa_A\left(\sum_i e_i\otimes t_i\right)=\sum_j\epsilon_j\otimes g_j(t_1,\dots,t_d),
\end{equation}
for $g_j \in \BC\llbracket x_1,\dots,x_d\rrbracket$ formal power series independent of $A$. We call $g_1,\ldots,g_m$ the \emph{Kuranishi relations} and we denote by $\BI_{\mathrm{Kur}}$ the ideal they generate.
The Kuranishi functor is prorepresented by the miniversal formal Kuranishi algebra
\[
\mathsf A \defeq \mathbb{C}\llbracket x_1,\dots,x_d\rrbracket / \BI_{\mathrm{Kur}}.
\]

\subsection{\texorpdfstring{A locally free replacement for $E$}{A locally free replacement for E}}
\label{sec:loc-free-replace}
Consider our fixed complex $E = R\oplus Q[-1]$, where $Q=\bigoplus_{1\leq \alpha \leq k} \OO_{q_\alpha} \otimes W_\alpha$, as in \Cref{eqn:polystable-rep}. 
Let $p \in \BP^3$ be the singular point of $R$, and fix three linear forms $\ell_1,\ell_2,\ell_3$ cutting out $p$, giving rise to a resolution of $R$ (cf.~\Cref{sec:resolution}). This gives us the perfect complex
\[
\begin{tikzcd}
C_0^\bullet = {\big[}\OO_{\BP^3}(-2) \arrow{r}{\rho_0} & \OO_{\BP^3}(-1)^{\oplus 3} {\big]}\,\in\,\Perf^{[-1,0]}(\BP^3),\qquad \rho_0 = 
\left(
\begin{smallmatrix}
    \ell_1 \\
    \ell_2 \\
    \ell_3
\end{smallmatrix}
\right),
\end{tikzcd}
\]
resolving $R$. 

\begin{notation}
For a vector $\boldit t=(t_1,t_2,t_3)^{\mathsf t}$, write
\[
S(\boldit t)=\begin{pmatrix}
0&-t_3&t_2\\
t_3&0&-t_1\\
-t_2&t_1&0
\end{pmatrix}.
\]
\end{notation}

For $\alpha=1,\ldots,k$ fix a vector of three linear forms $f_\alpha = (f_{\alpha,1},f_{\alpha,2},f_{\alpha,3})^{\mathsf t}$ cutting out $q_\alpha \in \Supp Q$, and form the associated Koszul complex
\begin{equation}
\label{Koszul-alpha}
\begin{tikzcd}
C_\alpha^\bullet = 
{\big[}\OO_{\BP^3}(-3) \arrow{r}{f_\alpha} 
& \OO_{\BP^3}(-2)^{\oplus 3} \arrow{r}{S(f_\alpha)} 
& \OO_{\BP^3}(-1)^{\oplus 3} \arrow{r}{f_\alpha^{\mathsf t}} 
& \OO_{\BP^3}
{\big]}\,\in\,\Perf^{[-2,1]}(\BP^3)
\end{tikzcd}
\end{equation}
resolving $\OO_{q_\alpha}[-1]$, see e.g.~\cite[Sec.~3.3]{Ricolfi-Savvas}. Then, the bounded complex
\[
C^\bullet = C_0^\bullet \oplus \bigoplus_{1\leq \alpha \leq k} C_\alpha^\bullet \otimes W_\alpha\,\in\,\Perf (\BP^3)
\]
is a perfect complex representing $E = R\oplus Q[-1]$.

\begin{assumption}
From now on, we assume the singular point $p$ of $R$ belongs to the support of $Q$. Without loss of generality, we assume $q_1=p$. Moreover, we set $f_1 = \rho_0$, i.e.~$f_{1,i}=\ell_i$ for $i=1,2,3$. We shall deal with the (easier) case $p \notin \Supp Q$ in \Cref{rmk:p-not-in-Q}.
\end{assumption}

\subsection{\texorpdfstring{The dgla attached to $C^\bullet$}{}}

Define the dgla
\[
L \defeq \Gamma(\BP^3,\lEnd^\bullet(C^\bullet)).
\]

A degree $n$ element $\xi \in L^{n} \defeq \prod_{r\in \BZ}\Hom_{\OO_{\BP^3}}(C^{r},C^{r+n})$ is a collection of maps $\xi^{r}\colon C^{r} \to C^{r+n}$ for all $r \in \BZ$.\footnote{In fact, finitely many maps, as $C^\bullet$ is bounded.} Moreover, $L$ has a differential $\dd_{L}$, a product, and a dg Lie bracket. The differential $\dd_{L}\colon L^{n} \to L^{n+1}$ is given by $\dd_{L}\xi = \dd_{C^\bullet}\circ \xi - (-1)^{\lvert \xi\rvert}\xi\circ \dd_{C^\bullet}$, i.e.~sends a homogeneous element $\xi \in L^{n}$ to $\dd_{L}\xi \in L^{n+1}$, where $(\dd_{L}\xi)^{r} \colon C^{r} \to C^{r+n+1}$ is the map $\dd_{C^\bullet}^{r+n}\circ \xi^{r} - (-1)^{n}\xi^{r+1}\circ \dd_{C^\bullet}^{r}$. The product of two homogeneous elements $\xi \in L^{n}$ and $\eta \in L^{k}$ is the element $\xi\eta \in L^{n+k}$ determined by $(\xi\eta)^{r} = \xi^{r+k}\circ \eta^{r}\colon C^{r} \to C^{r+k} \to C^{r+k+n}$. Finally, the dg Lie bracket is given by the graded commutator $[\xi,\eta] = \xi \eta - (-1)^{nk}\eta\xi$, so that in particular $[\xi,\xi] = 2\xi^{2}$ for a degree 1 element $\xi \in L^{1}$.

\begin{remark}
Note that $L^r=0$ for $r<0$. All summands of $\lEnd^r(C^\bullet)$ are line bundles $\OO_{\BP^3}(t)$ with $-3\leq t\leq3$, hence have no higher cohomology. Thus the natural inclusion into the Dolbeault~dgla~(cf.~\cite[Sec. 2.2]{Chen-Pertusi-Zhao})
\[
\begin{tikzcd}
L\arrow[r, hook]& A^{0,\bullet}(\BP^3,\lEnd^\bullet(C^\bullet))
\end{tikzcd}
\]
is a quasi-isomorphism. Consequently, $L$ computes $\RHom(C^\bullet,C^\bullet)$ and controls the deformation theory of $E$, see \cite[Lemma~2.4]{Chen-Pertusi-Zhao}.
\end{remark}

Set $E_0=R$, $W_0=\BC$, and $E_i = \OO_{q_i}[-1]$ for $i=1,\ldots,k$. If for all $i,j \in \set{0,1,\ldots,k}$ we define
\begin{equation}
\label{def:Lij}
\begin{split}
L_{ij} &\defeq \Gamma(\BP^3,\lHom^\bullet(C_i^\bullet,C_j^\bullet)), \\
\pmb{L}_{ij} &\defeq \Gamma(\BP^3,\lHom^\bullet(C_i^\bullet \otimes W_i,C_j^\bullet\otimes W_j)),  
\end{split}
\end{equation}
we have canonical isomorphisms
\begin{align*}
\pmb{L}_{ij} &\cong L_{ij} \otimes \Hom(W_i,W_j),\\
\HH^r(L_{ij}) &\cong \Ext^r(E_i,E_j),
\end{align*}
and a global decomposition
\begin{equation}\label{eq:block-decomposition}
  L = \bigoplus_{i,j=0}^k \pmb{L}_{ij}
  =\bigoplus_{i,j=0}^k
    L_{ij}\otimes\Hom(W_i,W_j).
\end{equation}
Indeed,
\begin{equation}
\begin{split}
\Ext^r(E_i,E_j) 
    &= \HH^r(\RHom(E_i,E_j)) \\
    &= \HH^r(\mathbf{R}\Gamma(\lHom^\bullet(C_i^\bullet,C_j^\bullet))) \\
    &= \HH^r(\Gamma(\lHom^\bullet(C_i^\bullet,C_j^\bullet))) \\
    &= \HH^r(\Hom^\bullet(C_i^\bullet,C_j^\bullet)) \\
    &=\HH^r(L_{ij}),
\end{split}
\end{equation}
where the equality $\mathbf{R}\Gamma(\lHom^\bullet(C_i^\bullet,C_j^\bullet))=\Gamma(\lHom^\bullet(C_i^\bullet,C_j^\bullet))$, follows from the fact that $\lHom^\bullet(C_i^\bullet,C_j^\bullet)$ is a complex of direct sums of line bundle of the form $\OO_{\BP^3}(r)$ with $-3\leq r\leq 3$, which are acyclic with respect to the global section functor $\Gamma$.

\subsection{A splitting of the dgla}
\label{sec:B-splitting}
This section is devoted to the choice of a splitting of $L$ of the form
\begin{equation}\label{splitting-choice}
\CH^r=\bigoplus_{i,j}\CH^r_{ij}\otimes\Hom(W_i,W_j),
\qquad
\CK^r=\bigoplus_{i,j}\CK^r_{ij}\otimes\Hom(W_i,W_j),
\end{equation}
where we ask compatibility with the source-target block decomposition
\eqref{eq:block-decomposition}, in the sense that $\CH^r_{ij}\subset Z^r(L_{ij})$ and $\CK^r_{ij} \subset L_{ij}^r$ will themselves define a splitting of $L_{ij}$. Since $G_E\defeq \Aut(E)=\BG_m\times\prod_\alpha\GL(W_\alpha)$ acts
on the $\Hom$-factors, such a splitting is $G_E$-invariant. In particular, the induced Kuranishi map is $G_E$-equivariant and thus the Kuranishi ideal is $G_E$-invariant, see \cite[Ex.~13.2.5]{Manetti-Book}.

\begin{notation}
We introduce a further piece of notation:
\[
\pmb{\mathcal{H}}_{ij}^r \defeq \mathcal{H}_{ij}^r\otimes \Hom(W_i,W_j),\qquad\pmb{\mathcal{K}}_{ij}^r \defeq \mathcal{K}_{ij}^r\otimes \Hom(W_i,W_j).
\]
\end{notation}

We specify certain choices for the splitting, as in Table~\ref{eq:B-representatives}. The remaining parts may be chosen arbitrarily, block by block, as in \Cref{splitting-choice}.

Each $L_{ij}$ is a direct sum of $\Hom$-spaces, see \Cref{def:Lij}. In the following table, we record the decomposition of the blocks which are relevant for us.

{\fontsize{8pt}{7.5pt}\selectfont
\begin{equation}\label{eq:structureoftheblocks}
\renewcommand{\arraystretch}{1.1}
\begin{array}{@{}c@{\quad}l@{\quad}l@{}}
\text{\bfseries $(i,j)$}
& \text{\bfseries Description of }L^1_{ij}
& \text{\bfseries Description of }L^2_{ij} \\
\midrule[1.2pt]
(0,0)
& \Hom(\OO(-2),\OO(-1)^{\oplus 3})
& 0 \\
\midrule[0.4pt]
(0,1)
& \Hom(\OO(-2),\OO(-1)^{\oplus 3})
  \oplus\Hom(\OO(-1)^{\oplus 3},\OO)
& \Hom(\OO(-2),\OO) \\
\midrule[0.4pt]
(1,0)
& \Hom(\OO(-3),\OO(-2))
  \oplus\Hom(\OO(-2)^{\oplus 3},\OO(-1)^{\oplus 3})
& \Hom(\OO(-3),\OO(-1)^{\oplus 3}) \\
\midrule[0.4pt]
(\alpha,\alpha)
& \Hom(\OO(-3),\OO(-2)^{\oplus 3})
  \oplus\Hom(\OO(-2)^{\oplus 3},\OO(-1)^{\oplus 3})
  \oplus\Hom(\OO(-1)^{\oplus 3},\OO)
& \Hom(\OO(-3),\OO(-1)^{\oplus 3})
  \oplus\Hom(\OO(-2)^{\oplus 3},\OO) \\
\bottomrule[0.6pt]
\end{array}
\end{equation}
}

\medskip
In what follows, we identify $\Hom(\OO(a)^{\oplus m},\OO(b)^{\oplus n})$ with the space of $n\times m$ matrices whose coefficients are degree $b-a$ homogeneous polynomials.
In order to define a splitting, we need to specify subspaces of $L_{ij}$.  We will define explicit bases, which we also exploit to present the Kuranishi map.

Choose a linear form $h$ independent of $\ell_1,\ell_2,\ell_3$ and $h_\alpha$ independent of  $f_{\alpha,1},f_{\alpha,2},f_{\alpha,3}$ for $\alpha \neq 1$. Also set $h_1=h$. Let $\boldit e_1,\boldit{e}_2,\boldit{e}_3$ be the standard basis
of $\BC^3$. Choose the following cohomology representatives,
where $r=1,2,3$ and the decompositions refer to the direct summands in \eqref{eq:structureoftheblocks}:
\begin{equation}\label{eq:B-representatives}
\renewcommand{\arraystretch}{1.1}
\begin{array}{@{}c@{\qquad}l@{\qquad}l@{}}
\text{\bfseries $(i,j)$}
& \text{\bfseries Basis of }\CH^1_{ij}
& \text{\bfseries Basis of }\CH^2_{ij} \\
\midrule[1.2pt]
(0,0)
& \lambda_r\defeq h\boldit{e}_r
& - \\
\midrule[0.4pt]
(0,1)
& a_r\defeq(h\boldit{e}_r,-h\boldit{e}_r^{\mathsf t})
& h^2 \\
\midrule[0.4pt]
(1,0)
& b\defeq(h,-h\id_3)
& h^2\boldit{e}_r \\
\midrule[0.4pt]
(\alpha,\alpha)
& X_{\alpha,r}\defeq
  (h_\alpha\boldit{e}_r,h_\alpha S(\boldit{e}_r),
  h_\alpha\boldit{e}_r^{\mathsf t})
& (h_\alpha^2\boldit{e}_r,h_\alpha^2\boldit{e}_r^{\mathsf t}) \\
\bottomrule[0.6pt]
\end{array}
\end{equation}

\medskip
To ease notation, we set $X_r = X_{1,r}$ for $r=1,2,3$. These matrices will correspond to the loops based at $\mathsf v_1$ (the vertex corresponding to $\OO_p[-1]$) in the Ext-quiver $\mathbf Q_E$, see \Cref{fig:Ext-quiver}.
It is an easy linear algebra check that these are indeed linearly independent representatives of the cohomology for the displayed blocks, which match the dimensions displayed in Table~\ref{eq:dimensionblocks}.
As mentioned above, the remaining parts of the splitting may be chosen arbitrarily. We only record the dimensions of the various relevant blocks in the following table, mainly to point out the vanishings:
\begin{equation}\label{eq:dimensionblocks}
\renewcommand{\arraystretch}{1.1}
\begin{array}{@{}c@{\qquad}c@{\qquad}c@{\qquad}c@{\qquad}c@{\qquad}c@{\qquad}c@{}}
\text{\bfseries $(i,j)$}
& (0,0)
& (0,1)
& (1,0)
& (\alpha,\alpha)
& (0,\alpha),\ \alpha\neq 1
& (\alpha,0),\ \alpha\neq 1 \\
\midrule[1.2pt]
\dim\mathcal{H}^1_{ij}
& 3 & 3 & 1 & 3 & 2 & 0 \\
\midrule[0.4pt]
\dim\mathcal{H}^2_{ij}
& 0 & 1 & 3 & 3 & 0 & 2 \\
\bottomrule[0.6pt]
\end{array}
\end{equation}
In the Table above, $\alpha$ ranges in $\Set{1,\dots,k}$ and for all the choices of $(i,j)$ not appearing in the table, both $\CH^1_{ij}$ and $\CH^2_{ij}$ vanish.

\subsection{The Ext-quiver and invariants}
\label{sec:extquiver}

\begin{figure}[H]
\centering
\begin{tikzpicture}[>=stealth,shorten >=2pt,looseness=.5,auto]
  \matrix [matrix of math nodes,
           column sep={3cm,between origins},
           row sep={3cm,between origins},
           nodes={circle, draw, minimum size=5.5mm}]
{ 
|(A)| \mathsf v_0 & |(B)| \mathsf v_\alpha \\         
};
\tikzstyle{every node}=[font=\small\itshape]

\node [anchor=west,right] at (-1.78,1.8) {$\BC$};
\node [anchor=west,right] at (1.2,1.8) {$\BC^{m_a}$};
\node [anchor=west,right] at (-1,0.7) {$\Hom(R,\OO_{q_a})$};    
\node [anchor=west,right] at (-1,-0.8) {$\Ext^2(\OO_{q_a},R)$}; 
\node [anchor=west,right] at (-4.5,0) {$\Ext^1(R,R)$}; 
\node [anchor=west,right] at (2.8,0) {$\Ext^1(\OO_{q_a},\OO_{q_a})$}; 

\draw (-1.52,1.6) -- (-1.52,0.3) {};
\draw (1.48,1.6) -- (1.48,0.32) {};
\draw[->] (A) to [bend left=25,looseness=1] (B) node [midway] {};
\draw[->, dotted] (B) to [bend left=25,looseness=1] (A) node [midway,below] {};
\draw[->] (B) to [out=-30,in=30,looseness=9] (B) node {};
\draw[->] (A) to [out=150,in=210,looseness=9] (A) node {};
\end{tikzpicture}
\caption{The part of the Ext-quiver $\mathbf Q_E$ involving the vertices $\mathsf v_0$ (corresponding to the $R$-summand) and $\mathsf v_\alpha$ (corresponding to $\OO_{q_\alpha}[-1]$). The arrow $\mathsf v_\alpha \to \mathsf v_0$ is optional: it is only there if $q_\alpha = p$. In this case, there are exactly 3 arrows $\mathsf v_0 \to \mathsf v_\alpha$. If $q_\alpha \neq p$, there are two arrows $\mathsf v_0 \to \mathsf v_\alpha$ and no arrows backwards.}
\label{fig:Ext-quiver}
\end{figure}
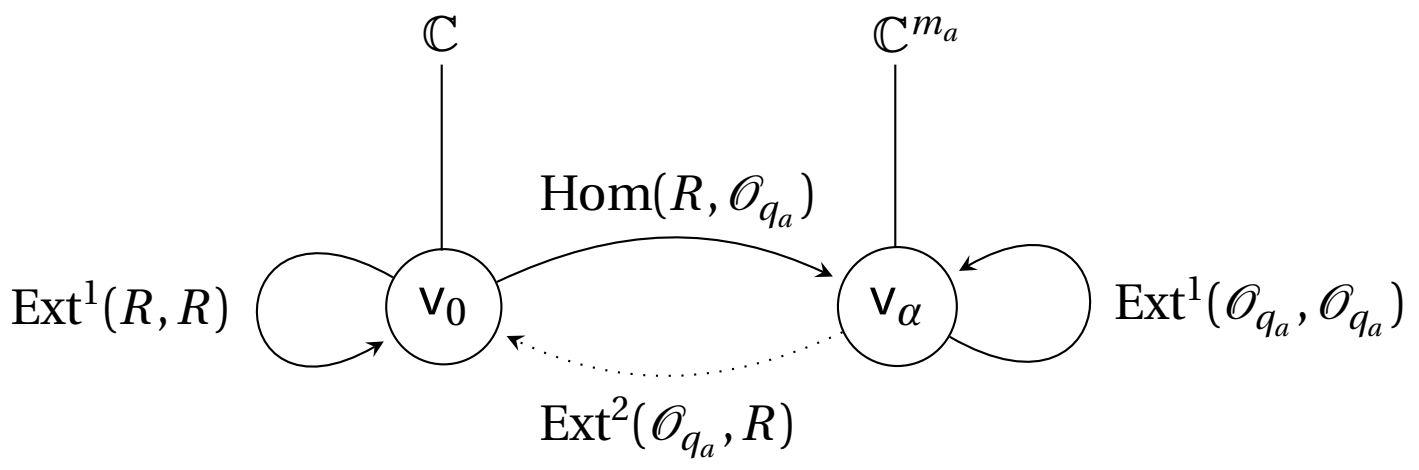

Set $m_\alpha \defeq \dim_{\BC} W_\alpha$, as in \Cref{eqn:polystable-rep} and form the dimension vector
$\boldit m=(1,m_1,\ldots,m_k)$.
The \emph{Ext-quiver} $\mathbf Q_E$ attached to $E=R\oplus Q[-1]$ (see e.g.~\cite{Toda:Moduli_Ext} for the general theory), represented in \Cref{fig:Ext-quiver}, has vertices
$\mathsf v_0,\mathsf v_1,\ldots,\mathsf v_k$ and one arrow $\mathsf v_i\to\mathsf v_j$
for each chosen basis vector of $\CH^1_{ij} \cong \Ext^1(E_i,E_j)$.

Introduce the polynomial $\BC$-algebra
\[
\mathsf S_E\defeq \BC[\Ext^1(E,E)^*] \defeq \Sym \Ext^1(E,E)^*\cong \Sym \mathcal{H}^{1,*}.
\]
Its completion at the maximal ideal of the origin 
\[
\widehat{\mathsf S}_E \defeq \BC\llbracket \Ext^1(E,E)^*\rrbracket \defeq \widehat{\Sym} \Ext^1(E,E)^*
\]
is exactly the formal power series ring appearing in \eqref{eqn:completed-local-map2}.

\begin{lemma}
\label{lem:B-traces}
The invariant ring $\mathsf S_E^{G_E}$ is generated by traces of
oriented cycles.
\end{lemma}

\begin{proof}
There is a $G_E$-equivariant identification
\cite{Toda:Moduli_Ext}
\[
\begin{tikzcd}
\Rep(\mathbf Q_E,\boldit m)\arrow{r}{\sim} & \CH^1
\end{tikzcd}
\]
where $G_E$ acts on $\Rep(\mathbf Q_E,\boldit m)$ by change of basis. It induces an equivariant isomorphism
\[
\BC[\Rep(\mathbf{Q}_E,\boldit m)]\cong \mathsf S_E
\]
between the corresponding coordinate rings. The result follows from the Le Bruyn--Procesi theorem, cf.~\cite{LeBruyn-Procesi,craw-yamagishi}.
\end{proof}

\subsection{The Kuranishi relations}
\label{sec:B-quadratic}

We now compute generators for the Kuranishi ideal $\BI_E \subset \widehat{\mathsf S}_E$. We first show that for our choice of splitting the Kuranishi map is quadratic.

\begin{prop}
\label{prop:B-quadratic}
The Kuranishi map attached to the splitting introduced in \Cref{sec:B-splitting} is quadratic, i.e. for $(A,\Fm_A)$ in $\Art_{\BC}$ and $x\in \mathcal{H}^1\otimes \mathfrak{m}_A$, we have
\[
    \kappa_A(x)=\frac{1}{2}H_2[x,x].
\]
\end{prop}
\begin{proof}
Let $p$ denote the singular point of $R$. We cover here the case where $p=q_1\in\Supp Q$, the other case is similar. Consider the subspace
\[
J\defeq \bigoplus_{2\leq \alpha \leq k}\pmb{\CK}^1_{0\alpha}\oplus \pmb{\CK}^1_{1\alpha}\subset L^1.
\]
Using the vanishings recorded in Table~\ref{eq:dimensionblocks}, we easily establish that
\begin{equation}\label{eq:relationsofJ}
    \delta_1[\mathcal{H}^1,J]\subset J,\qquad
    [J,J]=0,\qquad H_2[\mathcal{H}^1,J]=0.
\end{equation}

Indeed, $[J,J]=0$ because two elements of $J$ are not composable. Moreover, by Table~\ref{eq:dimensionblocks} the only nonzero brackets between $J$ and $\mathcal{H}^1$ are of the form $[\pmb{\mathcal{H}}^1_{00},J]\subset\bigoplus_{\alpha>1}\pmb{L}^2_{0\alpha}$, $[\pmb{\mathcal{H}}^1_{01},J]\subset\bigoplus_{\alpha>1}\pmb{L}^2_{0\alpha}$, $[\pmb{\mathcal{H}}^1_{10},J]\subset \bigoplus_{\alpha>1}\pmb{L}^2_{1\alpha}$,$[\pmb{\mathcal{H}}^1_{11},J]\subset\bigoplus_{\alpha>1}\pmb{L}^2_{1\alpha}$ and $[\bigoplus_{\alpha>1}\pmb{\mathcal{H}}^1_{\alpha\alpha},J]\subset\bigoplus_{\alpha>1}\pmb{L}^2_{0\alpha}\oplus\pmb{L}^2_{1\alpha}$ . Then, $\delta_1[\mathcal{H}^1,J]\subset J$ follows from the definition of the map $\delta_1$ and $H_2[\mathcal{H}^1,J]=0$ follows from $\pmb{\mathcal{H}}^2_{0\alpha}=\pmb{\mathcal{H}}^2_{1\alpha}=0$ for all $\alpha > 1$.

We can rewrite \Cref{kuranishi-map} as
\begin{equation}
\label{eq:formulaKfn}
    \kappa_A(x)=\frac{1}{2}\sum_{n,m}H_2[f_n(x),f_m(x)],
\end{equation}
where the terms $f_n(x)$ are defined in \Cref{rmk:expansion}.
We show by induction that $f_n(x)\in J\otimes\mathfrak{m}_A$ for all $n\geq 2$. First, notice that, $f_2(x)=-\frac{1}{2}\delta_1[x,x]\in J\otimes\Fm_A$. Indeed, the brackets between the blocks appearing in Table~\ref{eq:B-representatives} all land in $\mathcal{H}^2$, hence $\delta_1$ kills them. The remaining nonzero brackets land in
\[
\bigoplus_{2\leq \alpha \leq k}\pmb{L}^2_{0\alpha}\oplus \pmb{L}^2_{1\alpha},
\]
hence applying $\delta_1$, they are sent to $J$.

Combining the recursion to compute $f_m(x)$ given in \Cref{eq:recursion} with \Cref{eq:relationsofJ}, we conclude that $f_m(x)\in J\otimes\mathfrak{m}_A$, for all $m\geq 2$. Since $[J,J]=H_2[\mathcal{H}^1,J]=0$, we conclude that $H_2[f_n(x),f_m(x)]$ vanishes unless $n=m=1$.\footnote{Indeed, $f_1(x)=x\in\mathcal{H}^1\otimes\Fm_A$ and $f_m(x)\in J\otimes\Fm_A$, for all $m\geq 2$.} The claim then follows from \Cref{eq:formulaKfn}.
\end{proof}

After introducing the splitting and the bases for the relevant pieces of $\mathcal{H}^1$ and $\mathcal{H}^2$ (see \Cref{sec:B-splitting}), we are ready to compute the generators of the Kuranishi ideal $\BI_E \subset \widehat{\mathsf{S}}_E$, following the recipe of \Cref{sec:kuranishiideal}. In order to present the generators of $\BI_E$ in the most compact way, we exploit the Ext-quiver machinery of \Cref{sec:extquiver} and in particular the identification $\BC[\Rep(\mathbf Q_E,\boldit m)]\cong  \BC[\mathcal{H}^{1,*}]$.

Fix $i,j\in\set{0,1,\dots,k}$. Let $E_{i,j}$ be the set of arrows from $\mathsf v_i$ to $\mathsf v_j$ in the Ext-quiver $\mathbf Q_E$. A length 2 path $dc$ of the form
\[
\begin{tikzcd}
\mathsf v_i \arrow{r}{c} & \mathsf v_l \arrow{r}{d} & \mathsf v_j
\end{tikzcd}
\]
for $c \in \CH^1_{il}$ and $d \in \CH^1_{lj}$ defines a degree 2 element of $L^2$, that we still denote by $dc$.

Fix a basis $\epsilon_1,\dots, \epsilon_s$ of $\mathcal{H}^2_{ij}$ and let $\epsilon_1^\vee,\dots, \epsilon_s^\vee$ be its dual basis. Define, for $r=1,\ldots,s$, the scalars
\[
u_r(d,c)\defeq \epsilon^{\vee}_r(H_2(dc))\,\in\,\BC.
\]
Then, set
\begin{equation}
\label{eq:relationsmatrixform}
\mathsf M^{ij}_r=\sum_{l\in\{0,\dots,k\}}\sum_{c\in E_{i,l}}\sum_{d\in E_{l,j}}
u_r(d,c)dc.
\end{equation}
By definition, the sum is over the set of length-$2$ paths from $\mathsf v_i$ to $\mathsf v_j$. The expression in \Cref{eq:relationsmatrixform} gives a matrix with entries in $\mathsf S_E$. Denote by $I_E \subset \mathsf S_E$ the ideal generated by the entries of $\mathsf M_1^{ij},\ldots,\mathsf M_s^{ij}$, with $i,j$ ranging in $\set{0,1,\ldots,k}$. 
Combining \Cref{eqn:explicit-Kur} and \Cref{prop:B-quadratic}, we deduce that the extension of the ideal $I_E$ along the natural inclusion $\mathsf S_E \subset \widehat{\mathsf S}_E$ is precisely the Kuranishi ideal $\BI_E\subset \widehat{\mathsf S}_E$.

Explicitly, the ideal $I_E \subset \mathsf S_E$ is generated by the matrix entries of the following expressions.
\begin{equation}
\label{eq:B-relations}
\renewcommand{\arraystretch}{1.1}
\begin{array}{@{}l@{\qquad}l@{}}
\text{\bfseries Block} & \text{\bfseries Relations}\\
\midrule[1.2pt]
\pmb{\mathcal{H}}^2_{01}
& a_1\lambda_1+a_2\lambda_2+a_3\lambda_3-X_1a_1-X_2a_2-X_3a_3 \\
\midrule[0.4pt]
\pmb{\mathcal{H}}^2_{10}
& \lambda_1 b-bX_1,\quad
  \lambda_2 b-bX_2,\quad
  \lambda_3 b-bX_3 \\
\midrule[0.4pt]
\pmb{\mathcal{H}}^2_{11}
& a_1b+[X_2,X_3],\quad
  a_2b+[X_3,X_1],\quad
  a_3b+[X_1,X_2] \\
\midrule[0.4pt]
\pmb{\mathcal{H}}^2_{\alpha\alpha},\ \alpha>1
& [X_{\alpha,2},X_{\alpha,3}],\quad
  [X_{\alpha,3},X_{\alpha,1}],\quad
  [X_{\alpha,1},X_{\alpha,2}] \\
\bottomrule[0.6pt]
\end{array}
\end{equation}
We derive explicitly the relations associated with the blocks $\pmb{\mathcal{H}}^2_{10}$ and $\pmb{\mathcal{H}}^2_{11}$. These give the only generators of the Kuranishi ideal which are explicitly required to complete the main argument (cf.~\Cref{lem:B-offdiagonal-traces}).

We start with the block $\pmb{\mathcal{H}}^2_{10}$. The relevant part of the Ext quiver, namely the one containing the vertices $0,1$, is given by the following picture.
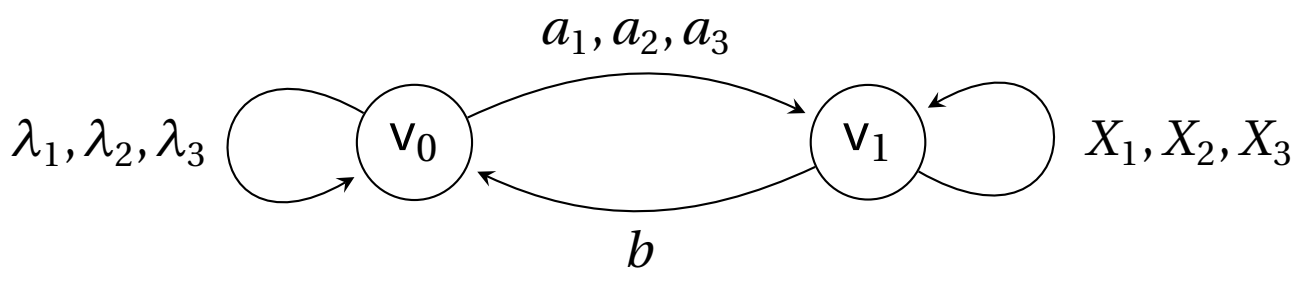
\begin{figure}[H]
\centering
\begin{tikzpicture}[>=stealth,shorten >=2pt,looseness=.5,auto]
  \matrix [matrix of math nodes,
           column sep={3cm,between origins},
           row sep={3cm,between origins},
           nodes={circle, draw, minimum size=5.5mm}]
{ 
|(A)| \mathsf v_0 & |(B)| \mathsf v_1 \\         
};
\tikzstyle{every node}=[font=\small\itshape]

\node [anchor=west,right] at (-0.8,0.7) {$a_1,a_2,a_3$};    \node [anchor=west,right] at (-0.25,-0.7) {$b$}; 
\node [anchor=west,right] at (-4.3,0) {$\lambda_1,\lambda_2,\lambda_3$}; 
\node [anchor=west,right] at (2.8,0) {$X_1,X_2,X_3$}; 

\draw[->] (A) to [bend left=25,looseness=1] (B) node [midway] {};
\draw[->] (B) to [bend left=25,looseness=1] (A) node [midway,below] {};
\draw[->] (B) to [out=-30,in=30,looseness=9] (B) node {};
\draw[->] (A) to [out=150,in=210,looseness=9] (A) node {};
\end{tikzpicture}
\end{figure}
The length-$2$ paths from $\mathsf v_1$ to $\mathsf v_0$ are
\[
bX_i, \quad \lambda_ib, \qquad \text{ for }i\in\Set{1,2,3}.
\]
By Table~\ref{eq:B-representatives}, we obtain
\[
    H_2(bX_i)=-h^2\boldit{e}_i,\qquad H_2(\lambda_ib)=h^2\boldit{e}_i,\qquad \text{for } i\in\{1,2,3\}.
\]
As seen in Table~\ref{eq:B-representatives}, the sequence $\epsilon_1\defeq h^2\boldit{e}_1,\epsilon_2\defeq h^2\boldit{e}_2,\epsilon_3\defeq h^2\boldit{e}_3$ provides a basis of $\mathcal{H}^2_{10}$. Thus the only nonzero coefficients appearing in the relations in \Cref{eq:relationsmatrixform} for the block $\pmb{\mathcal{H}}^2_{10}$ are given by
\begin{equation}
    u_r(b,X_r)=-1,\qquad u_r(\lambda_r,b)=1,\qquad \text{for } r\in\{1,2,3\}. 
\end{equation}
Thus, by \Cref{eq:relationsmatrixform}, we get the relations in the second row of Table \ref{eq:B-relations}, namely
\[
\mathsf M^{10}_1 = \lambda_1 b-bX_1,\quad \mathsf M^{10}_2 = \lambda_2 b-bX_2,\quad \mathsf M^{10}_3 = \lambda_3 b-bX_3.
\]

We now compute relations associated with the block $\pmb{\mathcal{H}}^2_{11}$. The paths of length $2$ from $\mathsf v_1$ to $\mathsf v_1$ are
\begin{align*}
    a_ib\qquad &\text{ for }i\in\{1,2,3\},\\
    X_iX_j\qquad &\text{ for }i,j\in\{1,2,3\}.
\end{align*}
By Table~\ref{eq:B-representatives}, we obtain
\begin{align*}
    H_2(a_r b)&=h^2(\boldit{e}_r,\boldit{e}_r^{\mathsf t}) & \text{for }r\in\{1,2,3\},\\
    H_2(X_i X_j)&=\sum_{1\leq l\leq 3}\varepsilon_{ijl} h^2(\boldit{e}_l,\boldit{e}_l^{\mathsf t}) & \text{for }i,j,l\in\{1,2,3\},
\end{align*}
where $\varepsilon_{ijl}$ is the Levi--Civita symbol (in particular, the displayed sum has only one nonzero term, corresponding to $l$ making $(i,j,l)$ a cyclic triple). The sequence $\epsilon_1\defeq h^2(\boldit{e}_1,\boldit{e}_1^{\mathsf t}),\epsilon_2\defeq h^2(\boldit{e}_2,\boldit{e}_2^{\mathsf t}),\epsilon_3\defeq h^2(\boldit{e}_3,\boldit{e}_3^{\mathsf t})$ provides a basis of $\mathcal{H}^2_{11}$. The nonzero coefficients appearing in \Cref{eq:relationsmatrixform} for the block $\pmb{\mathcal{H}}^2_{11}$ are given by
\begin{align*}
u_r(a_r,b)&=1, &\text{for }r\in\set{1,2,3},\\
u_l(X_i,X_j)&=\varepsilon_{ijl}, & \text{for }i,j,l\in\{1,2,3\}.
\end{align*}
Thus we obtain the relations in the third row of Table \ref{eq:B-relations}, namely
\[
\mathsf M^{11}_1 = a_1b+[X_2,X_3],\quad \mathsf M^{11}_2 = a_2b+[X_3,X_1],\quad \mathsf M^{11}_3 = a_3b+[X_1,X_2].
\]

\begin{lemma}
\label{lem:B-offdiagonal-traces}
The trace of every oriented cycle containing an off-diagonal arrow
belongs to $I_E$.
\end{lemma}
\begin{proof}
If $p\notin\Supp Q$, there are no cycles containing an off-diagonal arrow. Assume therefore that $q_1=p$. Every oriented cycle in $\mathbf{Q}_E$ containing an off-diagonal arrow must be based at $\mathsf v_0$ or $\mathsf v_1$. By cyclicity of the trace, without loss of generality, we can work with a cycle based at $\mathsf v_0$. 
    
Let $\Lambda_{0}$ (resp. $\Lambda_{1}$) denote the set of all words (possibly empty) in the loops at $\mathsf v_0$ (resp.~$\mathsf v_1$).  Any oriented cycle based at $\mathsf v_0$ and containing an off-diagonal arrow can be written in the form
    \[
        s_1 b t_1 a_{n_1} s_2 b t_2 a_{n_2} \cdots s_\ell b t_\ell a_{n_\ell} s_{\ell+1},
    \]
    for some $\ell \ge 1$, with $s_{\ell+1}, s_i \in \Lambda_{0}$, $t_i \in \Lambda_{1}$, and $n_i \in \{1,2,3\}$ for all $1 \leq i \leq \ell$.
    
    Working with coefficients modulo $I_E$, by Table \ref{eq:B-relations}, we have
    \[
       bX_i=\lambda_i b \mod I_E,
    \]
    for all $i\in\{1,2,3\}$. Consequently, for every $t\in\Lambda_1$, there exists $s\in\Lambda_0$ such that
    \[
       bt=s b \mod I_E.
    \]
Since words in $\Lambda_0$ can be taken out of the trace, as they are scalars, we obtain
    \[
        \tr(s_1 b t_1 a_{n_1} s_2 b t_2 a_{n_2} \cdots s_\ell b t_\ell a_{n_\ell} s_{\ell+1})=s_0\tr(b a_{n_1} b a_{n_2} \cdots b a_{n_\ell})\mod I_E,
    \]
    for some $s_0\in\Lambda_0$. Since $b a_{n_1},\dots b a_{n_{\ell-1}}$ are also scalars, we obtain
    \[
        s_0\tr(b a_{n_1} b a_{n_2} \cdots b a_{n_\ell})=s_0 ba_{n_1} \cdots ba_{n_{\ell-1}}\tr(ba_{n_\ell})\mod I_E.
    \]
    By Table \ref{eq:B-relations}, $a_{n_\ell}b=[X_j,X_k]\mod I_E$, for a suitable choice of $j$ and $k$. Thus
    \[
         \tr(ba_{n_\ell})=\tr(a_{n_\ell}b)=\tr [X_j,X_k] \mod I_E.
    \]
    This is enough to conclude since the trace of a commutator is zero.
\end{proof}

\subsection{Conclusion of the argument}
\label{sec:B-conclusion}
The natural choice of dgla's controlling the deformations of $R$ and $Q[-1]$ is
\[
L_R=\pmb{L}_{00},\qquad L_{Q}=\bigoplus_{\alpha=1}^k\pmb{L}_{\alpha\alpha}.
\]
The splitting chosen for $L$, which was compatible with its block decomposition, induces, by restriction, splittings for $L_R$ and $L_{Q}$. Repeating the construction of \Cref{sec:kuranishiideal}, we obtain the Kuranishi ideals
\[
\BI_R\subset\widehat{\mathsf S}_R=\mathbb{C}\llbracket \pmb{\mathcal{H}}^{1,\ast}_{00} \rrbracket,\qquad \BI_{Q}\subset\widehat{\mathsf S}_Q=\mathbb{C}\llbracket \oplus_\alpha\pmb{\mathcal{H}}^{1,\ast}_{\alpha\alpha} \rrbracket.
\]
Since there are no relations in Table \ref{eq:B-relations} associated with the blocks $\pmb{\mathcal{H}}^2_{00}$, we conclude that $\BI_{R}$ is the zero ideal. As for $\BI_Q$, it is the extension along $\mathsf S_Q \subset \widehat{\mathsf S}_Q$ of the ideal $I_Q\subset \mathsf S_Q$ generated by the entries of the matrices
\[
[X_{\alpha,2},X_{\alpha,3}],\ [X_{\alpha,3},X_{\alpha,1}],\ [X_{\alpha,1},X_{\alpha,2}],\quad \text{for }\alpha\in\{1,\dots k\}.
\]
Let us set
\[
\mathsf S_{R,Q}\defeq \mathsf S_R\otimes\mathsf S_Q.
\]
Setting the off-diagonal arrows to zero defines a
$G_E$-equivariant surjection
\begin{equation}
\label{eqn:GE-equiv-surj}
\begin{tikzcd}
  \mathsf S_E\arrow[two heads]{r} & \mathsf S_{R,Q},
\end{tikzcd}
\end{equation}
which sends the ideal $I_E$ onto the ideal 
\[
I_{R,Q}\defeq \mathsf S_R\otimes I_Q+I_R\otimes \mathsf S_Q \subset \mathsf S_{R,Q}.
\]
\begin{remark}
\label{rmk:relationbetweenideals}
Reductivity of $G_E$ implies that the induced map $I_E^{G_E}\rightarrow I_{R,Q}^{G_E}$ is surjective.
\end{remark}
We obtain an induced morphism
\begin{equation}
\label{eq:morphismexplained}
\begin{tikzcd}
\mathsf S_E^{G_E}/I_E^{G_E} \arrow{r}{\gamma} & \mathsf S_{R,Q}^{G_E}/I_{R,Q}^{G_E}.
\end{tikzcd}
\end{equation}

\begin{theorem}
\label{thm:B-completed-local-isomorphism}
For every closed point $s=(R,[Q])\in N\times\Sym^n(\BP^3)$, with $p = \Sing R$ in $\Supp Q$ and $E=R\oplus Q[-1]$, the canonical map
\[
\begin{tikzcd}
\widehat{\OO}_{\BM,[E]} \arrow{r}{\widehat u_s^{\#}} & \widehat{\OO}_{N \times  \Sym^{n}(\BP^3),s}
\end{tikzcd}
\]
is an isomorphism. Consequently, $u\colon N \times  \Sym^{n}(\BP^3) \to \BM$ is \'etale.
\end{theorem}

\begin{proof}
Since, by construction, the inclusion $L_Q\oplus L_R\subset L$ controls direct sums of deformations and preserves the chosen splittings, the identifications of \cite[Lemma 2.4]{Chen-Pertusi-Zhao} and \cite[Lemma 3.2]{zbMATH07585348} match the morphism $\widehat u_s^\#$ with the completion of the morphism $\gamma$ in \eqref{eq:morphismexplained}. Therefore, it is enough to prove that $\gamma$ is an isomorphism. Consider the commutative diagram with exact rows
\begin{equation}\label{iso-offdiag}
\begin{tikzcd}
0\arrow[r]&I_E^{G_E}\arrow[r]\arrow[d,"\alpha"]
&\mathsf S_E^{G_E}\arrow[r]\arrow[d,"\beta"]
&\mathsf S_E^{G_E}/I_E^{G_E}\arrow[r]\arrow[d,"\gamma"]&0\\
0\arrow[r]&I_{R,Q}^{G_E}\arrow[r]
&\mathsf S_{R,Q}^{G_E}\arrow[r]
&\mathsf S_{R,Q}^{G_E}/I_{R,Q}^{G_E}\arrow[r]&0
\end{tikzcd}
\end{equation}
where $\beta\colon\mathsf S_E^{G_E} \to \mathsf S_{R,Q}^{G_E}$ is the morphism setting the off-diagonal variables to zero, thus inducing $\gamma$. 

In order to show that $\gamma$ is an isomorphism, we apply the snake lemma. Surjectivity of $\gamma$ readily follows from surjectivity of $\beta$, which in turn follows from the surjectivity of \eqref{eqn:GE-equiv-surj} and the reductivity of $G_E$.

To show injectivity of $\gamma$ or equivalently $\ker\gamma=0$, it is enough to show that $\coker\alpha=0$ and $\ker\alpha=\ker\beta$. The map $\alpha$ is surjective (cf.~\Cref{rmk:relationbetweenideals}), i.e.~$\coker\alpha=0$. On the other hand, $\ker\alpha=\ker\beta$ is a straightforward consequence of the containment of ideals
\[
\ker\beta\subset I_E^{G_E},
\]
which follows from \Cref{lem:B-offdiagonal-traces}.
\end{proof}

\begin{remark}
\label{rmk:p-not-in-Q}
Suppose $p = \Sing R$ does not belong to $\Supp Q$. This case is much easier. Indeed, the statement of \Cref{lem:B-offdiagonal-traces} is trivially true, since there are no oriented cycles containing an off-diagonal arrow. The proof of \Cref{thm:B-completed-local-isomorphism} then works with no changes.
\end{remark}

\bibliographystyle{amsplain-nodash}
\bibliography{bib}

\bigskip

\medskip

\noindent
\emph{Charles Almeida}, \texttt{charlesalmeida@mat.ufmg.br} \\
\textsc{Universidade Federal de Minas Gerais (UFMG)  \\ Instituto de Ci\^encias Exatas (ICEx) \\ Departamento de Matem\'atica, 
Av. Ant\^onio Carlos, 6627, \\ 30123-970 Belo Horizonte-MG, Brazil}

\bigskip
\noindent
\emph{Nicolò Bignami}, \texttt{nbignami@sissa.it} \\
\textsc{SISSA, Via Bonomea 265, 34136 Trieste, Italy}

\bigskip
\noindent
\emph{Marcos Jardim}, \texttt{jardim@unicamp.br} \\
\textsc{Universidade Estadual de Campinas (UNICAMP)  \\ Instituto de Matem\'atica, Estat\'{\i}stica e Computação Cient\'{\i}fica (IMECC) \\ Departamento de Matem\'atica, 
Rua S\'ergio Buarque de Holanda, 651, \\ 13083-970 Campinas-SP, Brazil}

\bigskip
\noindent
\emph{Andrea T. Ricolfi}, \texttt{aricolfi@sissa.it} \\
\textsc{SISSA, Via Bonomea 265, 34136 Trieste, Italy}

\end{document}